\documentclass[12pt,a4paper]{book}
\usepackage{t1enc}
\usepackage[latin1]{inputenc}
\usepackage[english]{babel}
\usepackage{amsfonts}
\usepackage{amsmath}
\usepackage{amssymb}
\usepackage{amsthm}
\usepackage{amscd}
\usepackage[all]{xypic}
\usepackage{latexsym}
\usepackage{makeidx}
\makeindex
\date{}
\newcommand{\al}{\alpha}
\newcommand{\be}{\beta}
\newcommand{\cd}{\cdot}

\newcommand{\De}{\Delta}
\newcommand{\de}{\delta}
\newcommand{\e}{\epsilon}
\newcommand{\fr}{\frac }
\newcommand{\ga}{\gamma}
\newcommand{\Ga}{\Gamma}

\newcommand{\lam}{\lambda}
\newcommand{\lan}{\langle}

\newcommand{\C}{\mathbb{C}}

\newcommand{\Pro}{\mathbb{P}}
\newcommand{\Q}{\mathbb{Q}}
\newcommand{\R}{\mathbb{R}}
\newcommand{\Z}{\mathbb{Z}}
\newcommand{\mcal}{\mathcal}
\newcommand{\mf}{\mathfrak}
\newcommand{\mb}{\mbox}

\newcommand{\nf}{\normalfont}
\newcommand{\om}{\omega}
\newcommand{\Om}{\Omega}
\newcommand{\op}{\oplus}
\newcommand{\ot}{\otimes}

\newcommand{\ran}{\rangle}
\newcommand{\ra}{\rightarrow}
\newcommand{\si}{\sigma}
\newcommand{\Si}{\Sigma}
\newcommand{\ti}{\tilde}
\newcommand{\ul}{\underline}
\newcommand{\vphi}{\varphi}

\newtheorem{notation}{Notation}[section]
\newtheorem{com}[notation]{Comments}
\newtheorem{conj}[notation]{Conjecture}
\newtheorem{conv}[notation]{Convention}
\newtheorem{defn}[notation]{Definition}
\newtheorem{teo}[notation]{Theorem}
\newtheorem{lem}[notation]{Lemma}
\newtheorem{cor}[notation]{Corollary}
\newtheorem{prop}[notation]{Proposition}
\newtheorem{rem}[notation]{Remark}
\newtheorem{ex}[notation]{Example}

\begin{document}

\frontmatter

\begin{titlepage}

\title{ \ \ \ \ \ \ Orbifold Cohomology of \\ 
        \ \ \ \ \ \ $ADE$-singularities}

\author{\ \ \ \ \ \ Fabio Perroni\\
\vspace{3cm}\\
\ \ \ \ \ \ Advisors: Prof. Barbara Fantechi  and\\ \ \ \ \ \ \  \ \ \ \ \ \ \  Prof. Lothar G\"ottsche 
\vspace{5cm}\\
\ \ \ \ \ \ Submitted in partial fulfillment of the\\
\ \ \ \ \ \ requirements for the degree\\
\ \ \ \ \ \ of Doctor of Philosophy\\
\ \ \ \ \ \ in the International School for Advanced Studies\\ 
\ \ \ \ \ \ Via Beirut 2-4, 34014 Trieste, Italy}

\date{\ \ \ \ \  October 17, 2005}

\maketitle

\newpage

\vspace{2cm}

\noindent Fabio Perroni \\
Institut für Mathematik\\Universität Zürich\\Winterthurerstrasse 190\\CH-8057 Zürich\\
E-Mail: fabio.perroni@math.unizh.ch

\end{titlepage}

\tableofcontents

\chapter{Introduction}

Orbifolds arise in geometry in two different ways and as a consequence they can be given 
 two different descriptions.  On one hand, as a topological space $Y$  (respectively an algebraic variety) which is the union of open 
subsets $U$ of the form $U=\ti{U}/G$, where $\ti{U}$ is a smooth manifold (respectively a smooth variety) 
and $G$ is a finite group acting on it by  diffeomorphisms (respectively algebraically). 
The data $(U,\ti{U},G)$, as $U$ varies, define an orbifold structure $[Y]$ on $Y$. 
This way of thinking about orbifolds is in some sense 
concrete and geometric, but,  to define a good notion of maps between such objects requires some work. 
On the other hand, from moduli problems, orbifolds arise in a more abstract  but stil  natural way. 
Orbifolds arising from moduli problems are called smooth algebraic stacks in the sense of Deligne-Mumford \cite{DM}.
If $[Y]$ denotes such a stack, the coarse moduli  space $Y$ associated to it is,  locally in the \'{e}tale topology,
the quotient of a smooth variety by the action of a finite group. This variety parametrizes objects modulo isomorphisms. 

Gromov-Witten invariants for proper orbifolds have been defined first by Chen and Ruan \cite{CR1} for almost complex orbifolds.
Their construction arises from orbifold string theory introduced by Dixon, Harvey, Vafa and Witten \cite{D1},
\cite{D2}, for orbifolds which are global quotients. In the algebraic category, orbifold Gromov-Witten theory
has been developed by Abramovich, Corti, Graber and Vistoli \cite{ACV}, \cite{AGV}. The definition of  Gromov-Witten invariants for
orbifolds is similar to that for almost complex  manifolds and smooth varieties. 
The main change is that one replaces stable maps by twisted stable maps; these are orbifold morphisms
from nodal curves with orbifold structure to a target orbifold. These morphisms have to satisfy
some stability conditions that guarantee the existence of a compact  moduli space. 
Similar to the smooth case, the genus zero orbifold Gromov-Witten invariants can be used to define an orbifold quantum cohomology.
The degree zero part of the orbifold quantum cohomology, i.e. using only  degree zero maps, 
is the orbifold cohomology ring of $[Y]$.  

The orbifold cohomology of $[Y]$ has been defined by Chen and Ruan \cite{CR} for almost complex orbifolds. 
This has been extended to a noncommutative ring by Fantechi and G\"ottsche \cite{FG} in the case where the
orbifold is a global quotient. Abramovich, Graber and Vistoli defined the orbifold cohomology in the algebraic case \cite{AGV}.

The underlying vector space of the orbifold cohomology is the direct sum of the cohomology of $Y$ 
and the cohomology of the twisted sectors. These are orbifolds that parametrize points of $[Y]$ together with nontrivial 
automorphisms. Note that, if $U=\ti{U}/G$ is an open subset of $Y$ and $\ti{y} \in \ti{U}$, by an automorphism of $\ti{y}$ we mean
an element of the stabilizer of $\ti{y}$ in $G$. Twisted stable maps, being orbifold morphisms, 
carry informations about automorphisms of points. So, twisted sectors come naturally in the picture of orbifold cohomology.

We can see here a difference between the smooth case and the orbifold case. In the smooth case, the degree zero part 
of the quantum cohomology is the ordinary cohomology of the target space. In the orbifold case, the degree zero part contains 
the ordinary cohomology of $Y$ as a subring. 

Let  $[Y]$ be a global quotient orbifold, i.e. its underlying space $Y$ is of the form $Y=X/G$. 
Originating in Physics \cite{D1}, \cite{D2}, \cite{Z}, cohomological invariants of the orbifold $[Y]$ 
have been defined and studied  \cite{BB}. These invariants can be recovered from
the additive structure of $H^*_{orb}([Y])$. However, the orbifold cohomology ring contains much more information, 
notably the orbifold cup product.

\vspace{0.5cm}

Assume that $[Y]$ is a Gorenstein orbifold (see Definition \ref{gorenstein def}) and that the coarse moduli space $Y$ 
admits a crepant resolution $\rho :Z\ra Y$: 
this means  that $Z$ is a smooth variety, $\rho$ is an isomorphism outside the singular locus and $\rho^*(K_Y)=K_Z$
(see also Definition \ref{crepant def}).  
Here $K_Y$ and $K_Z$ are the canonical bundle of $Y$ and $Z$ respectively.

Motivated from Physics, Ruan stated the following \textit{cohomological crepant resolution conjecture} \cite{Ru}, see also Conjecture 
\ref{conj}.
Let $\rho :Z \ra Y$ be a crepant resolution, assume further that we have chosen an integral basis $\be_1,...,\be_n$ of
the kernel $\mb{Ker}~ \rho_*$ of the group homomorphism $\rho_*:H_2(Z,\Q)\ra H_2(Y,\Q)$. 
Then, we assign a formal variable $q_l$ to each class $\be_l$.
Using contributions from genus zero Gromov-Witten invariants of $Z$
whose homology classes belong to $\mb{Ker}~ \rho_*$, we can deform the ordinary cup product on 
$H^*(Z,\C)$ in order to  obtain  new ring structures. This deformation will depend on the complex parameters
$q_1,...,q_n$.  Giving these parameters the same value $-1$, we get the so called \textit{quantum corrected cohomology
ring}  $H^*_{\rho}(Z,\C)$. Then the conjecture states that, under suitable assumptions, $H^*_{\rho}(Z,\C)$
is isomorphic to $H^*_{orb}([Y])$. 

There are several motivations for studying this conjecture. It is related to the problem of 
understanding the behaviour of quantum cohomology under birational transformations. 
Another motivation comes from mirror symmetry as most of the known Calabi-Yau $3$-folds are crepant resolution 
of Calabi-Yau orbifolds.

The conjecture has been proved in the following cases: for $Y=S^{(2)}$,  the symmetric product of a compact complex surface $S$, 
and $Z=S^{[2]}$  the Hilbert scheme of two points in $S$, in \cite{LQ};
for $Y=S^{(r)}$,  the symmetric product of a complex projective surface $S$ with trivial canonical 
bundle $K_S \cong \mcal{O}_S$, and $Z=S^{[r]}$  the Hilbert scheme of $r$ points in $S$, in \cite{FG}, Theorem 3.10;
for $Y=V/G$  the quotient of  a  symplectic vector space of finite dimension by a finite group of symplectit automorphisms,
and $Z$  a crepant resolution, when it does exists, in \cite{GK}, Theorem 1.2.
Notice that in the last two cases, the  crepant resolution $Z$ has a holomorphic symplectic form, so the Gromov-Witten invariants
vanish and there are no quantum corrections. This means that the orbifold cohomology ring is isomorphic to
the ordinary cohomology ring of $Z$.

\vspace{0.5cm}

In this thesis, we study Ruan's conjecture for orbifolds with transversal $ADE$ singularities,
which are in a sense the simplest ones to which the conjecture applies. Indeed, assume that 
in the orbifold $[Y]$ the open subset of points with trivial automorphisms group is dense.
If $[Y]$ is Gorenstein, and even under the much more general assumptions that no chart
$(U,\ti{U}, G)$ contains a pseudoreflection, then the orbifold structure $[Y]$ can be recovered
by the singular variety or analytic space $Y$. The simplest possible Gorenstein singularities are 
the $ADE$ surface singularities (or rational double points also called  Du Val singularities). 

An orbifold $[Y]$ has transversal $ADE$ singularities if, \'etale locally, the coarse moduli space
$Y$ is isomorphic to a product $R\times \C^k$, where $R$ is a germ of an $ADE$ singularity.
In general for any Gorenstein orbifold $[Y]$, there exists a closed subset $W \subset Y$
of codimension $\geq 3$ such that $Y\backslash W$ has transversal $ADE$ singularities.
So the case we study is in a natural sense the simplest, ignoring higher codimension
phenomena. On the other hand, transversal $ADE$ singularities occur natural in many contexts:
transversal $A_1$ in $S^{(2)}$ for $S$ a surface; transversal $A_n$ in many complete intersections
in weighted projective spaces; they also are, at least locally, examples of symplectic orbifolds.

In this thesis, we lay down the structure to deal with general $ADE$ singularities.
After that, we concentrate on the transversal $A_n$ case, with a further mild technical 
assumption we address Ruan's conjecture by computing explicitly both the orbifold cohomology
and the quantum corrections. The former is achieved in general, for the later we have an explicit conjecture
(see Conjecture \ref{mia congettura}) which is only verified under additional, and 
somewhat unnatural, technical assumptions. We include a sketch of an argument
which would prove the conjecture in general, assuming some technical results which we cannot so far prove.
The conjecture is proven fully in the transversal $A_1$ case. 
Finally, we construct an explicit isomorphism between the orbifold cohomology ring $H^*_{orb}([Y])$
and the quantum corrected cohomology ring  $H^*_{\rho}(Z)$ in the transversal $A_1$ case,
verifying Ruan's conjecture. In the transversal $A_2$ case,
the quantum corrected $3$-point function (see Chapter 2.4.2) can not be evaluated in $q_1=q_2=-1$,
so that $H^*_{\rho}(Z,\C)$ is not defined. Thus Ruan's conjecture has to be slightly modified.
However we have found that, if  $q_1=q_2$ is a third root of the 
unit different to $1$, then the resulting ring $H^*(Z)(q_1,q_2)$
is isomorphic to the orbifold cohomology ring $H^*_{orb}([Y])$. Also in this case we give an 
explicit isomorphism between them (see Theorem \ref{a2}). 
Thus proving a slightly modified version of Ruan's conjecture. We expect that
in the  $A_n$ case, the same modification of Ruan's conjecture holds, i.e. if $q_1=...=q_n$ is  an $(n+1)$-th root of unity
such that $q_1 \cd \cd \cd q_n \not= 1$,  the resulting ring $H^*(Z)(q_1,...,q_n)$
is isomorphic to the orbifold cohomology ring $H^*_{orb}([Y])$.

\vspace{1.2cm}

\noindent \textbf{The structure of the thesis is the following.}

\vspace{0.2cm}

In Chapter 1 we collect some basic definitions on orbifolds, morphisms of orbifolds
and orbifold vector bundles. In Chapter 2  we first review the definition of orbifold 
cohomology ring for a complex orbifold, then we state the 
\textit{cohomological crepant resolution conjecture} as given by Ruan in  \cite{Ru}. 
In Chapter 3 we define orbifolds with transversal $ADE$-singularities, 
see Definition \ref{definition ADE arbifolds}. Then we give a description of the twisted sectors in general.
Finally we specialize to orbifolds with transversal $A_n$-singularities and, under the technical assumption of trivial monodromy,
we compute the orbifold cohomology ring. 
In Chapter 4 we study the crepant resolution. We first  show that any variety with transversal $ADE$-singularities $Y$  
has a unique crepant resolution $\rho :Z\ra Y$, Proposition \ref{existence and unicity of res}. 
Then we restrict our attention to the case of transversal $A_n$-singularities and trivial monodromy and we 
give an explicit description of the cohomology ring of $Z$. 
Chapter 5 contains the computations of the Gromov-Witten invariants of $Z$ in the $A_n$ case. We also give  
a description of the quantum corrected cohomology ring of $Z$. In Chapter 6 
we prove Ruan's conjecture in the $A_1$ case and, in the $A_2$ case with minor modifications.

\chapter{Acknowledgements}

First of all, I would like to thank very much my Advisors, 
Professor Barbara Fantechi and Professor Lothar G\"ottsche, for having introduced me to the 
topic of this thesis, for all they taught me over the past years and for all the discussions
which helped me to overcome many problems.  

\vspace{0.2cm}

I would like to thank the members of the Mathematical Physics Sector of SISSA for giving me the opportunity to 
enlarge my scientific knowledge and for the stimulating atmosphere that I found here. 
In particular, I would like to thank Professor Ugo Bruzzo which helped me mathematically and morally
through a couple of difficult times. I would also like to thank Professor Mudumbai Narasimhan for 
introducing me to the fascinating world of Algebraic Geometry, his way of thinking about  Geometry 
has been influential and stimulating for me.

I am very grateful to Professor Christian Okonek for inviting me to give a talk on the subject
of this thesis at the University of Z\"urich. His interest in this work and his suggestions 
have been very important for me.

I am grateful to Jim Bryan and Tom Graber for sharing and discussing with me the results of their paper \cite{BG}.
I would like to thank Professor Yongbin Ruan, I appreciated very much his interest and enthusiasm in my work.
I also owe thanks to several people for useful discussions and helpful advice.

\mainmatter

\chapter{Orbifolds}

We collect here some basic definitions. First of all we define  smooth (or $C^{\infty}$) orbifolds
and complex holomorphic orbifold using local charts, then we give some examples of orbifolds to explain
the definition. We give the definition of a morphism between orbifolds  and we show that the category of orbifolds 
is a $2$-category, which means that given two orbifolds the set of morphisms between them has a structure of a groupoid. 
Finally we define orbifold vector bundles over a given orbifold. As an example of orbifold vector bundle,
we will recall the definition of the tangent bundle to any given orbifold.

\section{Definition of orbifold}

The notion of orbifold\index{orbifold} was first introduced by Satake in \cite{Sa} 
under the name of \textit{$V$-manifold}. A characterization of orbifolds, as defined by 
Satake, in terms of their sheaves is given in \cite{MP}. A different definition
was given by Chen and Ruan in \cite{CR1}. The second definition is slightly more general than the first,
indeed it leads to the notion of non-reduced orbifold while an orbifold 
as defined by Satake is called reduced. In \cite{CR1}, section 4.1, the authors prove 
that the two definitions are equivalent, once restricted to reduced orbifolds.

In this paper we will follow \cite{Sa} and \cite{MP} with minor modifications in order to include 
non-reduced orbifolds.

Let $Y$ be a paracompact Hausdorff topological space. 
A \textit{uniformizing system}\index{uniformizing system ! for an open subset}
for an open subset $U\subset Y$ is a collection of the following objects:
\begin{description}
\item{$\ti{U}$}\quad  a connected open subset of $\R^d$;
\item{$G$}\quad  a finite group of $C^{\infty}$-automorphisms of $\ti{U}$
        such that: the fixed-point set of each element of the group is either the whole
        space or of codimension at least $2$, the multiplication in $G$ is given by $g_1 \cd g_2 = g_1 \circ g_2$
        where $\circ$ is the composition;
\item{$\chi$}\quad  a continuous map $\ti{U} \ra U$ that induces an homeomorphism from $\ti{U}/G$ to $U$,
        where $\ti{U}/G$ is the quotient space with the quotient topology. Here $G$  acts  on $\ti{U}$ on the left.
\end{description}

\noindent We will call the subgroup of $G$ which consists of elements fixing the whole space
the \textit{kernel} of the action\index{kernel of the action}, and it will be denoted by $\rm{Ker}(G)$.

\vspace{0.5cm}

\begin{notation}\nf
Given an open subset $U$ of $Y$, a uniformizing system for $U$
will be denoted by $(\ti{U}, G_U, \chi_U)$. If the dependence on $U$ is clear from the context,
it will also be denoted by $(\ti{U}, G, \chi)$.
\end{notation}

\vspace{0.2cm}

\begin{defn} The \textbf{dimension}\index{uniformizing system ! dimension of} of an uniformizing system 
$(\ti{U}, G, \chi)$ is the dimension of $\ti{U}$ as a real manifold.
\end{defn}

\vspace{0.5cm}

Let $(\ti{U}, G, \chi)$ and $(\ti{U}', G', \chi')$ be uniformizing systems for $U$ and $U'$ respectively, and let $U \subset U'$.
An \textit{embedding} between such uniformizing systems\index{uniformizing system ! embedding of}
is a pair $(\vphi, \lam)$, where $\vphi :\ti{U} \ra \ti{U}'$
is a smooth embedding such that $\chi' \circ \vphi= \chi$ and 
$\lam :G \ra G'$ is a group homomorphism such that $\vphi \circ g = \lam (g) \circ \vphi$
for all $g\in G$. Furthermore, $\lam$ induces an isomorphism from $\rm{Ker}(G)$
to $\rm{Ker}(G')$.

\vspace{0.5cm}

\begin{defn} \label{atlas}
An \textbf{orbifold atlas}\index{orbifold atlas} on $Y$ is a family $\mcal{U}$ of uniformizing 
systems for open sets in $Y$ satisfying the following conditions:
\begin{enumerate}
\item The family ${\chi (\ti{U})}$ is an open covering of $Y$,  for $(\ti{U}, G, \chi)\in \mcal{U}$.
\item Let $(\ti{U}, G, \chi),(\ti{U}', G', \chi')\in \mcal{U}$ be uniformizing systems for $U$ and $U'$ respectively, 
and let $U \subset U'$. Then there exists an embedding $(\vphi, \lam):(\ti{U}, G, \chi)\ra (\ti{U}', G', \chi')$.
\item Let $(\ti{U}, G, \chi),(\ti{U}', G', \chi')\in \mcal{U}$ be uniformizing systems for $U$ and $U'$ respectively.
Then, for any point $y\in U\cap U'$, there exists an open neighbourhood $U'' \subset U\cap U'$
of $y$ and a uniformizing system $(\ti{U}'', G'', \chi'')$ for $U''$ which belong to the family $\mcal{U}$.
\end{enumerate}

Two such atlases are said to be equivalent if they have a common refinement,
where an atlas $\mcal{U}$ is said to refine $\mcal{V}$ if for every chart in $\mcal{U}$
there exists an embedding into some chart in $\mcal{V}$.

A \textbf{smooth orbifold}\index{orbifold ! smooth} structure on a paracompact Hausdorff topological space $Y$ 
is an equivalence class of orbifold atlases on $Y$.
\end{defn}

\vspace{0.5cm}

\begin{notation}\nf We denote by $[Y]$ the smooth orbifold structure on the topological space $Y$. We will 
call this simply an orbifold.
\end{notation}

\vspace{0.2cm}

\begin{defn} The orbifold $[Y]$ is said of \textbf{real dimension}\index{orbifold ! dimension} $d$ if all the uniformizing 
systems of an atlas have dimension $d$. The real dimension of $[Y]$ will be denoted by $\rm{dim}_{\R}[Y]$ or, by abuse of 
notation by $\rm{dim}[Y]$ .
\end{defn}

\vspace{0.2cm}

\begin{rem}\label{rem1}\nf
Every orbifold atlas for $Y$ is contained in a unique maximal one, and two orbifold atlases
are equivalent if and only if they are contained in the same maximal one.
Therefore we shall often tacitly work with a maximal atlas.
\end{rem}

\vspace{0.2cm}

\begin{prop}[\nf{\cite{MP} Proposition A.1}]\label{pretransition} \nf
Let $(\vphi, \lam)$ and $(\psi, \mu)$ be two embeddings from $(\ti{U}, G, \chi)$ to $(\ti{U}', G', \chi')$.
Then, there exists a $g' \in G$ such that 
        \begin{align*}
        \psi = g' \circ \vphi \quad \rm{and} \quad \mu = g' \cd \lam \cd {g'}^{-1}.
        \end{align*} 
Moreover, this $g'$ is unique up to composition by an element of $\rm{Ker}(G')$.
\end{prop}

\vspace{0.2cm}

\begin{rem}\label{rem2}\nf Notice that, for any uniformizing system $(\ti{U}, G, \chi)$ and $g\in G$,
there exists a $g' \in G'$ such that $\vphi \circ g = g' \circ \vphi$. Moreover this $g'$
is unique up to an element in $\rm{Ker}(G')$. 
Thus, in the definition of an embedding $(\vphi, \lam)$, the existence of $\lam$ is required
to guarantee a continuity for the kernels of the actions.
\end{rem}

\begin{rem}\label{rem3}\nf
The previous remark implies that, for any embedding $(\vphi, \lam)$, the group homomorphism $\lam $ is injective.
\end{rem}

\vspace{0.2cm}

\begin{lem}[\nf{\cite{MP}, Lemma A.2}]\label{rem4} \nf
Let $(\vphi, \lam): (\ti{U}, G, \chi) \ra (\ti{U}', G', \chi')$ be an embedding. 
If $g' \in G'$ is such that $\vphi(\ti{U}) \cap (g' \circ \vphi)(\ti{U})\not= \emptyset $, then
$g'$ belongs to the image of $\lam$. 
\end{lem}

\vspace{0.2cm}

\begin{rem}\label{rem5} \nf
Let $(\ti{U}, G, \chi)$ be a uniformizing system for the open subset $U$ of $Y$. Let $U' \subset U$
be an open subset. Then we have an \textit{induced uniformizing system}\index{uniformizing system ! induced}
$(\ti{U}', G', \chi')$ for $U'$,  where $\ti{U}'$ is a connected component of $\chi^{-1}(U')$
and $G'$ is the maximal subgroup of $G$ that acts on $\ti{U}'$.
Clearly there is an embedding of $(\ti{U}', G', \chi')$ in  $(\ti{U}, G, \chi)$.

It follows that, for a given orbifold $[Y]$, we can choose an orbifold atlas $\mcal{U}$ 
arbitrarily fine.
\end{rem}

\begin{rem}\label{rem6} \nf
Let $[Y]$ be an orbifold and let $(\ti{U}_1, G_1, \chi_1)$, $(\ti{U}_2, G_2, \chi_2)$
be uniformizing systems for open subsets $U_1$, $U_2$ of $Y$ in the same orbifold structure
$[Y]$. For any point $y\in U_1 \cap U_2$, there is an open neighbourhood $U_{12}$ of $y$ 
such that $U_{12}\subset U_1 \cap U_2$, a uniformizing system $(\ti{U}_{12}, G_{12}, \chi_{12})$ for $U_{12}$,
compatible with $[Y]$, and embeddings $(\vphi_i, \lam_i) :(\ti{U}_{12}, G_{12}, \chi_{12}) \ra (\ti{U}_i, G_i, \chi_i)$ for 
$i\in \{1,2\}$.
So, we get the isomorphism:  
        \begin{align*}
        \vphi_{12}:=\vphi_2 \circ \vphi_1^{-1} : \vphi_1(\ti{U}_{12}) \ra \vphi_2(\ti{U}_{12}).
        \end{align*}

Let now $U_1$, $U_2$ and $U_3$ be open subsets of $Y$ such that $U_1 \cap U_2 \cap U_3 \not= \emptyset$, 
and assume that there are uniformizing systems $(\ti{U}_1, G_1, \chi_1)$, $(\ti{U}_2, G_2, \chi_2)$ and 
$(\ti{U}_3, G_3, \chi_3)$ for $U_1$, $U_2$ and $U_3$ respectively. Then, from Proposition \ref{pretransition},
there exists $g\in G_3$ such that 
        \begin{align*}
        \vphi_{23} \circ \vphi_{12} = g\circ \vphi_{13},
        \end{align*}
where the equation holds if we restrict the functions to some open subsets of the domains.
\end{rem}

\vspace{0.2cm}

\begin{defn} 
A \textbf{reduced}\index{orbifold ! reduced} orbifold is an orbifold structure $[Y]$ on $Y$ such that there
exists an orbifold atlas $\mcal{U}$ for $[Y]$ with the following propertiy:
for any uniformizing system $(\ti{U}, G, \chi) \in \mcal{U}$, $\rm{Ker}(G)$ is the trivial group.
\end{defn}

\vspace{0.2cm}

\begin{defn} \label{u.s. at y def}
Let $[Y]$ be a smooth orbifold and $y\in Y$ be a point. 
A \textbf{uniformizing system for $[Y]$ at $y$}\index{uniformizing system ! at a point}
is given by an open neighbourhood $U_y$ of $y$ in $Y$ and a uniformizing system $(\ti{U}, G, \chi)$ for $U_y$
in the orbifold structure $[Y]$ such that, $\ti{U}\subset \R^d$ is 
a ball centered in the origin $0\in \R^n$, $G$ acts trivially on $0$ and $\chi^{-1}(y)=0$.
\end{defn}

\vspace{0.2cm}

\begin{notation}\label{u.s. at y}\normalfont
For a given orbifold $[Y]$ and a point $y\in Y$, a uniformizing system at $y$
will be denoted by $(\ti{U}_y, G_y, \chi_y)$ and  $\chi(\ti{U}_y)$ by $U_y$.
The group $G_y$ will be also called the local group at $y$\index{local group}.
\end{notation}

\vspace{0.5cm}

We now define a complex orbifold. We will use the same notation as in the smooth case.

\vspace{0.2cm}

Let $Y$ be a paracompact Hausdorff topological space. A \textit{complex uniformizing system}\index{uniformizing system ! complex}
for an open subset $U$ of $Y$ is a triple $(\ti{U}, G, \chi)$, where
$\ti{U} \subset \C^d$ is a connected open subset, $G$ is a finite group of holomorphic 
automorphisms of $\ti{U}$ and $\chi$ is a continuous map satisfying the same properties 
required in 
the smooth case.

\vspace{0.2cm}

\begin{defn} The \textbf{complex  dimension}\index{uniformizing system ! complex dimension of} of a complex uniformizing system 
$(\ti{U}, G, \chi)$ is the dimension of $\ti{U}$ as a complex manifold.
\end{defn}

\vspace{0.2cm}

Let $(\ti{U}, G, \chi)$ and $(\ti{U}', G', \chi')$ be complex uniformizing systems for $U$ and $U'$ respectively, 
and let $U \subset U'$.
A \textit{complex embedding} between such uniformizing systems\index{uniformizing system ! complex embedding of}
is a pair $(\vphi, \lam)$ satisfying the same properties stated in the smooth case
but  where $\vphi :\ti{U} \ra \ti{U}'$ is holomorphic.

\begin{defn}
A \textbf{complex orbifold atlas}\index{orbifold !  atlas ! complex} on $Y$ is a family $\mcal{U}$ of complex uniformizing 
systems for open sets in $Y$ satisfying the conditions 1., 2. and 3. of Definition \ref{atlas}
where we replace embeddings with complex embeddings.

Two such atlases are said to be equivalent if they have a common refinement,
where an atlas $\mcal{U}$ is said to refine $\mcal{V}$ if for every chart in $\mcal{U}$
there exists a complex embedding into some chart in $\mcal{V}$.

A \textbf{complex orbifold}\index{orbifold ! complex} structure on a paracompact Hausdorff topological space $Y$ 
is an equivalence class of complex orbifold atlases on $Y$.
\end{defn}

\vspace{0.2cm}

\begin{defn}\label{orbifold complex dimension}
 The complex orbifold $[Y]$ is said of \textbf{complex dimension}\index{orbifold ! complex dimension} $d$ if 
all the uniformizing systems of a complex  atlas have complex dimension $d$. The complex 
dimension of $[Y]$ will be denoted by $\rm{dim}_{\C}[Y]$ or, by abuse of 
notation by $\rm{dim}[Y]$ .
\end{defn}

\vspace{0.2cm}

\begin{rem} \nf
Proposition \ref{pretransition}, Lemma \ref{rem4} and Remarks \ref{rem1}, \ref{rem2},
\ref{rem3}, \ref{rem5} and \ref{rem6}, holds in the complex case. The notions of reduced complex orbifold
and of uniformizing systems at a point are defined for complex orbifolds in the same way of the smooth case. 
\end{rem}

\vspace{0.2cm}

\begin{lem}[\nf{\textit{Linearization lemma},\cite{C} Theorem 4.}]\index{Linearization lemma} \nf
Let $[Y]$ be a complex orbifold and  let $y\in Y$ be a point. Then we can choose a local uniformizing
system  $(\ti{U}_y, G_y, \chi_y)$ at $y$ such that  $G_y$ acts linearly on $\ti{U}_y$.
\end{lem}

\section{Examples of orbifolds} 

In this section we give some examples of orbifolds, principally of complex orbifolds.

\begin{ex}\label{smooth manifolds as orbifolds}[\nf{\textit{Smooth and complex manifolds}}] \nf 
Let $Y$ be a smooth (resp. complex) manifold. Let $(U,\al)$ be a chart, 
where $U\subset Y$ is open and $\al :U \ra \R^d$ (resp. $\C^d$) is an homeomorphism with the image.
Then $(\ti{U} = \al (U), G=\{1\}, \chi = \al^{-1})$ is a smooth (resp. complex) uniformizing  system for $U$.
Let $(U', \al')$ be another chart in the same smooth (resp. complex) atlas of $(U,\al)$.
Then the previous construction gives smooth (resp. complex) uniformizing system for $U\cap U'$ and for $U'$.
Moreover we have smooth (resp. complex) embeddings corresponding to $U\cap U' \subset U, U'$.
It follows that a smooth (resp. complex) manifold has a natural orbifold structure.
\end{ex} 

\begin{ex}[\nf{\textit{Global quotients}}]\label{global quotient} \nf
Let $X$ be a smooth (resp. complex) manifold with the
action of a finite group of diffeomorphisms (resp. biholomorphic transformations) $G$. 
Let $Y$ be the quotient space with the quotient topology, $Y=X/G$.
Clearly $\mcal{U}= \{ (U=X, G, \chi)\}$ is a smooth (resp. complex) orbifold atlas for $Y$, where 
$\chi : X \ra Y$ is the quotient map. So the class 
of $\mcal{U}$ define a smooth (resp. complex) orbifold structure over $Y$. This orbifold structure
will be denoted by $[X/G]$ and is called \textit{global quotient}\index{global quotient}. 
\end{ex}

\begin{ex}[\nf{\textit{Weighted projective space} $[\Pro(a_0,...,a_d)]$}]\label{wps} \nf
Let $a_0,...,a_d$ be positive integers.
Then  $\Pro(a_0,...,a_d)$ is, by definition, the quotient
        \begin{align*}
        \Pro(a_0,...,a_d)= ( \C^{d+1}-\{0\} )/ \C^*
        \end{align*}
of $\C^{d+1}-\{0\}$ under the equivalence relation 
        \begin{align*}
        \left( x_{0},...,x_{d} \right) \thicksim ( {\lam}^{a_{0}} x_{0},...,{\lam}^{a_{a}} x_{d} )
        \qquad \mbox{for} \quad \lam \in \C^*.
        \end{align*}
For any $\left( x_{0},...,x_{d} \right) \in \C^{d+1}-\{0\}$ we will denote with 
$[ x_{0},...,x_{d}]\in \Pro(a_0,...,a_d)$ the equivalence class of $\left( x_{0},...,x_{d} \right)$.

Let $X_i =\{x_i \not= 0\} \subset \C^{d+1}-\{0\}$ and $U_i =X_i/ \C^*$.
For any $i=0,...,d$, let $\ti{U}_i=\C^d$ with coordinates $(z_0,...,\hat{z_i},...,z_d)$
and $G_i= \mu_{a_i}$, the group of $a_i$-th rooth of unity, acting on $\ti{U}_i$ as follows
        \begin{align*}
        \e \cd (z_0,...,\hat{z_i},...,z_d)= (\e^{-a_0}z_0,..., \e^{-a_d}z_d) \qquad \mbox{for} \quad \e \in  \mu_{a_i}.
        \end{align*}
Then  the map $\ti{U}_i \ra \Pro(a_0,...,a_d)$ defined by 
        \begin{align*}
        (z_0,...,\hat{z_i},...,z_d) \mapsto [z_0,...,z_{i-1},1,z_{i+1},...,z_d]
        \end{align*}
induces an isomorphism between $\ti{U}_i/\mu_{a_i}$ and $U_i$.

In order to give an orbifold atlas on $\Pro(a_0,...,a_d)$ it is enough to construct,
for any $[x]\in \Pro(a_0,...,a_d)$, a uniformizing system at $[x]$, say $(\ti{U}_{[x]}, G_{[x]}, \chi_{[x]})$,
and embeddings $(\ti{U}_{[x]}, G_{[x]}, \chi_{[x]}) \ra (\ti{U}_i, \mu_{a_i}, \chi_i)$
whenever $U_{[x]} \subset U_i$. 

Let $[x] \in U_i$, then take any $z\in \ti{U}_i$ such that $[z] =[x]$. 
Define $G_{[x]}$ to be the stabilizer in $\mu_{a_i}$ of $z$. For $\ti{U}_{[x]}$ take a small ball in $\ti{U}_i$ such that $G_{[x]}$
acts on it and, for any $g \in \mu_{a_i}-G_{[x]}$, $g(\ti{U}_{[x]}) \cap \ti{U}_{[x]} =\emptyset$.
Let $U_{[x]}$ be the quotient of $\ti{U}_{[x]}$ by $G_{[x]}$ and $\chi_{[x]}:\ti{U}_{[x]} \ra U_{[x]}$
be the quotient map. By construction 
there is an embedding $(\ti{U}_{[x]}, G_{[x]}, \chi_{[x]}) \ra (\ti{U}_i,\mu_{a_i}, \chi_i)$. 

If $U_{[x]} \subset U_j$ for some $j\not= i$ we have to construct an embedding 
        \begin{align*}
        (\ti{U}_{[x]}, G_{[x]}, \chi_{[x]}) \ra (\ti{U}_j,\mu_{a_j}, \chi_j).
        \end{align*}
Notice that, if $U_i \cap U_j \not= \emptyset$, for any $z \in \ti{U}_i$ such that $[z] \in U_i \cap U_j$
there is a biholomorphic map with domain a suitable neighbourhood of $z$ in $\ti{U}_i$ and values  in $\ti{U}_j$
defined as follow
        \begin{align*}
        \ti{\phi}_{ij}: (z_0,...,\hat{z_i},...,z_d) \mapsto \left( \fr{z_0}{z_j^{a_0/a_j}},...,\fr{1}{z_j^{a_i/a_j}},...,
        \hat{z_j},...,\fr{z_d}{z_j^{a_d/a_j}}\right)
        \end{align*}
where $z_j^{1/a_j}$ is a choosed $a_j$-th  rooth of $z_j$. Then using this we obtain the required embedding.
Notice that, if $z'=\ti{\phi}_{ij}(z)$, then  $[z'] =[x]$ and
the stabilizers of $z'$ in $\mu_{a_j}$ and  of $z$ in $\mu_{a_i}$
are isomorphic.
\end{ex}

\begin{ex} [\nf{\textit{Hypersurfaces in weighted projective spaces}}]\label{hwps} \nf 
We use the same notation of the previous example.
A polynomial $F$ in the variables $x_0,...,x_d$ is \textit{weighted homogeneous}\index{weighted homogeneous} of degree $D$, 
if the following holds
        \begin{align*}
        F({\lam}^{a_{0}} x_{0},...,{\lam}^{a_{d}} x_{d})=\lam^D F(x_0,...,x_d),  
        \end{align*}
for any $\lam \in \C^*$.
In this case, the equation $F=0$ defines a closed subvariety $Y$ of  $\Pro(a_0,...,a_d)$.

Assume that  the affine variety $\{F=0\}$ is smooth in $\C^{d+1}-\{0\}$. Than the same construction 
of the previous example can be used to give an orbifold structure over $Y$.
\end{ex}

\vspace{0.5cm}

\begin{rem} \nf
There is a more algebraic construction of $\Pro(a_0,...,a_d)$ which goes as follows (\cite{Do}). 
Define $S(a_0,...,a_d)$ to be the polynomial ring $\C [x_0,...,x_d]$ graded by the condition $deg(x_i)=a_i$, for $i\in \{0,...,d \}$.
Then 
        \begin{align*}
        \Pro(a_0,...,a_d) = \mbox{Proj}(S(a_0,...,a_d)).
        \end{align*}
Here we use the definition given in \cite{H} for the Proj of a graded $\C$-algebra.

The condition for $F$ to be weighted homogeneous of degree $D$ means that $F$ is homogeneous of degree $D$ in  $S(a_0,...,a_d)$.

Let $Y$ to be a variety as in the last example. Then the canonical sheaf $K_Y$ of $Y$
is given as follows
        \begin{align*}
        K_Y =\mcal{O}_Y(D-\sum_{i=0}^d a_i),
        \end{align*}
where, for any integer $m$, $\mcal{O}_Y(m)$ is the sheaf associated to the graded module $\left(S(a_0,...,a_d)/(F) \right)(m)$,
\cite{Do}. See Notation \ref{canonical sheaf def} for the definition of $K_Y$.

This construction has been extensively used to give examples of orbifold with trivial canonical bundle, i.e. Calabi-Yau orbifold.
\end{rem}

\begin{rem} \nf
In \cite{BCS}, the authors introduces the notion of  \textit{toric Deligne-Mumford stack},
this is an orbifold structure over a simplicial toric variety.
Note that a simplicial toric variety is an algebraic variety with quotient 
singularities (\cite{Fu2}, Section 2.2).

A toric Deligne-Mumford stack corresponds to a combinatorial object called a \textit{stacky fan}.
A stacky fan $\mathbf{\Sigma}$ is a triple consisting of a finitely generated abelian group
$N$, a simplicial fan $\Si$ in $\Q \ot_{\Z} N$ with $d$ rays (see \cite{Fu2}, Section 2.2,
for the definition of a simplicial fan), and a map $\be :\Z^d \ra N$
where the image of the standard basis of $\Z^d$ generates the rays in $\Si$.

A rational simplicial fan $\Si$ in some lattice $N\cong \Z^d$ (where rational means that any simplex of the fan
is generated by vectors in the lattice) gives a canonical stacky fan $\mathbf{\Si} = (N, \Si, \be)$
where $\be$ is the map defined by the minimal lattice points on the rays.
Hence, there is a natural toric Deligne-Mumford stack associated to every simplicial toric variety.

Notice that weighted projective spaces and hypersurfaces in weighted projective spaces are
toric varieties and the orbifold structures constructed in Examples \ref{wps} and  \ref{hwps}
are examples of toric Deligne-Mumford stacks.
\end{rem}

\section{Morphisms of orbifolds}

In this section we review the notion of morphism between two given orbifolds and of 
\textit{natural transformation}\index{natural transformation} 
between two morphisms of the same orbifolds. We will see that orbifolds and morphisms form a $2$-category:
$1$-morphisms are morphisms and $2$-morphisms are natural transformations.

Satake in the paper \cite{Sa} defined  \textit{$C^{\infty}$-map} between orbifolds. 
With this definition, smooth orbifolds and $C^{\infty}$-maps 
form a category. It turns out that, given a $C^{\infty}$-map from the orbifold $[X]$ to the orbifold $[Y]$,
and an orbifold vector bundle $[E]$ over $[Y]$,
it is not possible, in general, to pull-back $[E]$ using this map.
This point is explained, for example, in \cite{CR1} Section 4.4, \cite{MP} Section 2
and in \cite{PRY} Section 2.

In order to be able to pull-back vector bundles, the correct definition of a morphism between two
orbifolds is that of \textit{strong map} defined by I. Moerdijk and D. A. Pronk in \cite{MP} (Section 5),
or \textit{good map} defined by W. Chen and Y. Ruan in \cite{CR1}, Definition 4.4.1.
As proved in \cite{LU}, Proposition 5.1.7, the two definitions are equivalent.

\vspace{0.5cm}

\begin{defn} \label{compatible system}
Let $[X]$ and $[Y]$ be two orbifolds, $f:X\ra Y$ be a continuous map. 
A \textbf{compatible system}\index{compatible system} for $f$
is given by the following objects:
        \begin{enumerate}
        \item two atlases $\mcal{V}$ and $\mcal{U}$ for $[X]$ and $[Y]$ respectively;
        \item a correspondence that associates to any uniformizing system $(\ti{V}, H, \rho )$ in $\mcal{V}$, an uniformizing system 
                $(\ti{U}, G, \chi )$ in $\mcal{U}$ and a smooth function $f_{\ti{V}} :\ti{V} \ra \ti{U}$,
                such that 
                \begin{eqnarray*}
                \chi \circ f_{\ti{V}} & = & f \circ \rho;
                \end{eqnarray*}         
        \item a correspondence that associates to any embedding $(\psi , \mu) :(\ti{V}, H, \rho ) \ra  (\ti{V}', H', \rho' )$
                an embedding $(\vphi, \lam ):(\ti{U}, G, \chi ) \ra (\ti{U}', G', \chi' )$
                between the corresponding uniformizing systems such that
                \begin{eqnarray*}
                f_{\ti{V}'} \circ \psi & = & \vphi \circ f_{\ti{V}}.
                \end{eqnarray*}                 
        \end{enumerate}
Moreover we require that the assignment of the above objects is functorial with respect to the composition of embeddings
$(\psi , \mu) :(\ti{V}, H, \rho ) \ra  (\ti{V}', H', \rho' )$ and $(\psi' , \mu') :(\ti{V}', H', \rho' ) \ra  (\ti{V}'', H'', \rho'' )$.
\end{defn}

\vspace{0.2cm}

\begin{rem}\label{f_H} \nf
From the definition we have  that for any uniformizing system $(\ti{V}, H, \rho )$ in $\mcal{V}$, let 
$(\ti{U}, G, \chi )$ be the corresponding uniformizing system of $\mcal{U}$, then 
there is a group homomorphism $f_H: H \ra G$ such that
        \begin{align*}
        f_{\ti{V}} \circ h = f_H(h) \circ f_{\ti{V}} \qquad \mbox{for any $h\in H$}.
        \end{align*}
\end{rem}

\begin{notation} \nf A compatible system for $f$ will be denoted by $\ti{f} : \mcal{V} \ra \mcal{U}$,
where $\mcal{V}$ and $\mcal{U}$ are atlases of $[X]$ and $[Y]$ respectively. For any uniformizing system
$(\ti{V}, H, \rho)\in \mcal{V}$, we will denote by $\ti{f}(\ti{V}, H, \rho)$ the corresponding uniformizing system in 
$\mcal{U}$. For any embedding $(\psi, \mu )$ in $\mcal{V}$ we will denote by $\ti{f}(\psi, \mu )$ the corresponding embedding
in $\mcal{U}$.
\end{notation}

\vspace{0.2cm}

\begin{lem}[\nf{\cite{CR1} Remark 4.4.7}] 
Let $\ti{f}_i : \mcal{V}_i \ra \mcal{U}_i$, for $i\in \{1,2\}$, be two compatible systems of the same function $f:X\ra Y$.
Then there exists a common refinement $\mcal{V}$ of both $\mcal{V}_1$ and $\mcal{V}_2$, a refinement $\mcal{U}$ of both 
$\mcal{U}_1$ and $\mcal{U}_2$, and  compatible systems $\ti{f}_i' :\mcal{V} \ra \mcal{U}$, $i\in \{1,2\}$, such that
$\ti{f}_i'$ is induced by $\ti{f}_i$  for any $i\in \{1,2\}$. 
\end{lem}

\vspace{0.2cm}

\noindent The previous Lemma means that given two compatible systems we can always assume that they are defined over the same atlases.

\vspace{0.2cm}

\begin{defn}
Two compatible systems $\ti{f}_i : \mcal{V}_i \ra \mcal{U}_i$, for $i=1,2$, of the same map $f$, are said  \textbf{equivalent}
if  they coincide when defined over the same atlases.

A \textbf{morphism}\index{morphism between orbifolds} $[f]$ from the orbifold $[X]$ to the orbifold $[Y]$ is
given by a continuous map $f:X\ra Y$ and  an equivalence class of compatible systems for $f$.
\end{defn}

\vspace{0.5cm}

We now review the notion of natural transformation between two morphisms. 

\begin{defn}
Let $\ti{f}_1$ and $\ti{f}_2$ be compatible systems for $f:X\ra Y$. Assume that they are defined 
over the same atlases $\mcal{V}$ and $\mcal{U}$
of $[X]$ and $[Y]$ respectively.
A \textbf{natural transformation}\index{natural trnsformation} from $\ti{f}_1$ to $\ti{f}_2$ is 
a correspondence that associates to any uniformizing system $(\ti{V}, H, \rho) \in \mcal{V}$
an automorphism $\de_{\ti{V}}$ of $\ti{f}_1(\ti{V}, H, \rho) = \ti{f}_2(\ti{V}, H, \rho)$ 
such that 
        \begin{align*}
        (f_2)_{\ti{V}} = \de_{\ti{V}} \circ (f_1)_{\ti{V}},
        \end{align*}
and such that, for any embedding $(\psi, \mu): (\ti{V}, H, \rho) \ra (\ti{V}', H', \rho')$ in $\mcal{V}$, 
the following diagram commutes
        \begin{eqnarray*}
        \begin{CD}
         \ti{f}_1(\ti{V}, H, \rho) @>\ti{f}_1(\psi, \mu) >> \ti{f}_1(\ti{V}', H', \rho')\\
        @V\de_{\ti{V}} VV       @VV\de_{\ti{V}'}V \\
        \ti{f}_2(\ti{V}, H, \rho) @>\ti{f}_2(\psi, \mu) >> \ti{f}_2(\ti{V}', H', \rho').
        \end{CD}
        \end{eqnarray*}
\end{defn}

\vspace{0.2cm}

\begin{rem} \nf Let $\ti{f}_1, \ti{f}_2 : \mcal{V} \ra \mcal{U}$ be compatible systems for $f:X\ra Y$,
and let $\{ \de_{\ti{V}}~:~ (\ti{V},H, \rho) \in \mcal{V} \}$ be a natural transformation between them.
Let $ \mcal{V}'$ and $\mcal{U}'$  be refinement of $ \mcal{V}$ and $\mcal{U}$ respectively.
Then there is a natural transformation $\{ \de_{\ti{V}}'~:~ (\ti{V}',H', \rho') \in \mcal{V}' \}$
between the compatible systems $\ti{f}_1', \ti{f}_2' : \mcal{V}' \ra \mcal{U}'$ that are induced by $\ti{f}_1$ and $\ti{f}_2$.
We will say that $\{ \de_{\ti{V}}\}$ and $\{ \de_{\ti{V}}'\}$ are equivalent.
\end{rem}

\vspace{0.2cm}

\begin{defn}
Let $[f]_1$ and $[f]_2$ be morphisms from $[X]$ to $ [Y]$. A \textbf{natural transformation}\index{natural transformation}
between   $[f]_1$ and $[f]_2$ is an equivalence class of  natural transformations between two compatible systems 
that represents $[f]_1$ and $[f]_2$.
\end{defn}

\vspace{0.2cm}

\begin{notation} \nf We will denote a natural transformation between the morphisms $[f_1]$ and $[f_2]$
by $[f_1] \Rightarrow [f_2]$. We will use the same symbol for natural transformations between compatible systems.
\end{notation}

\vspace{0.2cm}

\begin{rem}\label{orbifold as groupoid}\nf \textit{(Orbifolds as groupoids)}
Let  $[Y]$ be an orbifold. Any atlas $\mcal{U}$ of $[Y]$ determines a groupoid which represents $[Y]$.
This is shown in \cite{MP2} Theorem 4.1.1; for a detailed study of the relations between orbifolds and groupoids
see \cite{MP}, \cite{Mo}.

By abuse of notation, we will denote also by $\mcal{U}$ the groupoid associated to the covering $\mcal{U}$. 

Let $[X]$ and $[Y]$ be orbifolds and $f:X\ra Y$ a continuous map. Then a compatible system 
$\ti{f} : \mcal{V} \ra \mcal{U}$ induces a morphism of the corresponding groupoids.
Conversely, for any pair of groupoids $\mcal{V}$ and $\mcal{U}$  that represent $[X]$ and $[Y]$
respectively, and a groupoid morphism $\ti{f} : \mcal{V} \ra \mcal{U}$, there is an induced compatible system for $f$. 
This is the content of Proposition 5.1.7 in \cite{LU}. Moreover, natural transformations between compatible systems
correspond to natural transformations between the associated morphisms of groupoids (see \cite{Mo} Section 2.2
for the definition of natural transformation of morphisms of groupoids).
So, the notion of compatible system is the same of \textit{strong map} introduced in \cite{MP}, Section 5.
\end{rem}

\begin{rem} \label{orbifold as stack} \nf \textit{(Orbifolds as stacks)}
There is another way to define orbifolds, that is as 
\textit{smooth Deligne-Mumford stacks} (\cite{DM}, \cite{E}, \cite{V}).
The way to pass from our definition to stacks is through groupoids. Indeed, given an atlas $\mcal{U}$,
we can construct the corresponding groupoid $\mcal{U}$, then there is a procedure that associates to the groupoid
$\mcal{U}$ a stack (\cite{V}, Appendix).
\end{rem}

\vspace{0.2cm}

We  now show that our definition of morphisms between orbifolds corresponds to the definition of morphisms between stacks.

\begin{prop}\label{proof of orbifolds=stacks}
The $2$-category of orbifolds with orbifold morphisms is equivalent to the $2$-category of smooth Deligne-Mumford 
stacks with morphisms of stacks.
\end{prop}

\noindent \textbf{Proof.} Let $(\mbox{Orb})$ denote the category of orbifolds and orbifold morphisms,
and let $(\mbox{Stacks})$ denote the category of stacks and stack morphisms. From Remark \ref{orbifold as stack}
it follows that we have a correspondence that associates to any object $[Y]$ in $(\mbox{Orb})$
an object $\mcal{Y}$ in $(\mbox{Stacks})$. We now construct a correspondence that associates to any morphism
in $(\mbox{Orb})$ a morphism in $(\mbox{Stacks})$ in such a way that this gives a functor $(\mbox{Orb}) \ra (\mbox{Stacks})$
such that, for any pair of orbifolds $[X]$ and $[Y]$, the map
        \begin{equation}\label{morphismi orb > morphismi stacks}
        \mbox{Mor}_{\mbox{(Orb)}}([X], [Y]) \ra \mbox{Mor}_{\mbox{(Stacks)}}(\mcal{X},\mcal{Y})
        \end{equation}
is bijective. Note that this would be an equivalence, indeed any smooth Deligne-Mumford stack is represented by a groupoid
(\cite{V} Appendix), then apply Remark \ref{orbifold as groupoid}.

Let $[f] :[X] \ra [Y]$ be a morphism  of orbifolds.  We construct a morphism $F : \mcal{X} \ra \mcal{Y}$
of stacks in the following way: 
        \begin{enumerate}
        \item we give open coverings $\{[\ti{V}_i/H_i]\}_{i\in I}$ of $\mcal{X}$ and $\{[\ti{U_j}/G_j]\}_{j\in J}$ of $\mcal{Y}$; 
        \item we define morphisms $F_i : [\ti{V}_i/H_i] \ra  [\ti{U}_{\nu(i)}/G_{\nu(i)}]$, 
                where $\nu : I \ra J$ is a function; 
        \item for any $i, i' \in I$, we give a natural transformation $\de_{ii'} : F_{ii'} \Rightarrow  F_{i'i}$, 
                such that $\de_{i'i''}\circ \de_{ii'} = \de_{ii''}$ on triple intersections, 
                where $F_{ii'}$ denote the restriction of $F_i$ to 
                $[\ti{V}_i/H_i] \times_{\mcal{X}} [\ti{V}_{i'}/H_{i'}]$.
        \end{enumerate}
Let $\ti{f}: \mcal{V} \ra \mcal{U}$ be a compatible system representing $[f]$.
Then define  
        \begin{align*}
        \{[\ti{V}_i/H_i]\}_{i\in I} := \{[\ti{V}/H]~ : ~ (\ti{V}, H, \rho) \in \mcal{V} \}
        \end{align*}
and 
        \begin{align*}
        \{[\ti{U_j}/G_j]\}_{j\in J} :=\{[\ti{U}/G]~ : ~ (\ti{U},U , \chi) \in \mcal{U} \}.
        \end{align*}
The function $\nu :I \ra J$ is the correspondence given in point 2 of Definition \ref{compatible system}.
We now define  a morphism $F_i : [\ti{V}_i/H_i] \ra  [\ti{U}_{\nu(i)}/G_{\nu(i)}]$ for any $i\in I$.
Let $f_{\ti{V}_i} : \ti{V}_i \ra \ti{U}_{\nu (i)}$ and $f_{H_i}: H_i \ra G_{\nu (i)}$ given as in Definition \ref{compatible system}
and Remark \ref{f_H}. An object in $[\ti{V}_i/H_i]$ over the base $B$ is represented by the following  diagram:
        \begin{eqnarray*}
        \begin{CD}
        &P& @>\al >> \ti{V}_i \\
        &@VVV& \\       
        &B&
        \end{CD}
        \end{eqnarray*}
where $P\ra B$ is a principal $H_i$-bundle, and $\al : P\ra \ti{V}_i$ is an $H_i$-equivariant smooth function.
This object can be also given by an open covering $\{B_k\}_K$ of $B$ and the following data
        \begin{eqnarray*}
        &h_{lm}& \in H_i, ~ \mbox{for any $l,m$ such that} B_l \cap B_m \not= \emptyset, ~ \mbox{such that}~ h_{mn}\cd h_{lm}=h_{ln}\\
        &\al_k &: B_k \ra \ti{V}_i,~ \mbox{smooth functions such that}~ \al_l =h_{lm} \circ \al_m
        \end{eqnarray*}
modulo the equivalence relation that identifies $h_{lm}$ with $h_l \cd h_{lm} \cd h_m^{-1}$ and $\al_k$ with $h_k \circ \al_k$
for $h_k \in H_i$, for any $k\in K$.
So, using this description of the objects, we define $F_i$ on the objects to be  given by the following correspondence
        $$
        F_i : \{ h_{lm},~\al_k \}/\ti \mapsto \{ g_{lm} =f_{H_i}(h_{lm}),~ \be_k =f_{\ti{V}_i}\circ \al_k\}/\ti.
        $$
On morphisms, $F_i$ is defined in an analogous way.

Let $i, i' \in I$. Then, the natural transformation $\de_{ii'} : F_{ii'} \Rightarrow  F_{i'i}$ is constructed as follows.
We can cover $[\ti{V}_i/H_i] \times_{\mcal{X}} [\ti{V}_{i'}/H_{i'}]$ using open embeddings that 
are induced by embeddings of uniformizing systems as follows
        \begin{align*}
        \Psi: [\ti{V}_{lm}/H_{lm}]\ra  [\ti{V}_i/H_i] \times_{\mcal{X}} [\ti{V}_{i'}/H_{i'}], 
        \end{align*}
where $(\ti{V}_{lm}, H_{lm}, \rho_{lm}) \in \mcal{V}$.
Point 3 of Definition \ref{compatible system} gives natural transformations from the restriction of
$F_{ii'}$ on $[\ti{V}_{lm}/H_{lm}]$ to  the restriction of $F_{i'i}$ on $[\ti{V}_{lm}/H_{lm}]$. 
The functoriality assumption in Definition \ref{compatible system}
ensures that these natural transformations patch together to give $\de_{ii'}$. Finally, the cocycle condition
on natural transformations follows from the functoriality condition given in Definition \ref{compatible system}.

\vspace{0.2cm}

To show that the map  (\ref{morphismi orb > morphismi stacks})
is bijective, we will construct an inverse. Let
        \begin{align*}
        F: \mcal{X} \ra \mcal{Y}
        \end{align*}
be a morphism of stacks. Let $\ti{x}\in X$ be a point, and let $(\ti{V}_x, H_x, \rho_x )$ be a uniformizing system at $x$. 
We assume that $\ti{V}_x$ is simply connected. We can assume that the restriction of $F$ to
$[\ti{V}_x/ H_x ]$ take values in a global quotient stack $[\ti{U}/G]$. So, we have 
        \begin{align*}
        F_x: [\ti{V}_x/ H_x ] \ra [\ti{U}/G].
        \end{align*}
The following object 
        \begin{eqnarray}\label{ob}
        \begin{CD}
        &\ti{V}_x \times H_x& @>>> \ti{V}_x \\
        & @VVV &\\
        & \ti{V}_x &
        \end{CD}
        \end{eqnarray}
has  group of automorphisms $H_x$, where the horizontal map $\ti{V}_x \times H_x \ra \ti{V}_x$ is given by $(\xi, k)  \mapsto k(\xi)$
and the vertical map is the projection.
Indeed, for any $h\in H_x$, the following diagram is cartesian
        \begin{eqnarray*}
        \begin{CD}
        \ti{V}_x \times H_x @>>> \ti{V}_x \times H_x \\
        @VVV @VVV \\
        \ti{V}_x @>>> \ti{V}_x
        \end{CD}
        \end{eqnarray*}
where the vertical arrows are projections, the  lower horizontal is $\xi \mapsto h(\xi)$ and the upper horizontal
$(\xi, k) \mapsto ( h(\xi), k \cd h^{-1})$.
So we get a group homomorphism $f_{H_x} : H_x \ra G$. Now notice that 
we have a morphism $\ti{V}_x \ra [\ti{U}/G]$ which is the composition of $\ti{V}_x \ra [\ti{V}_x/ H_x ]$
with $F_x$. This correspond to the following diagram
        \begin{eqnarray}\label{imob}
        \begin{CD}
        &Q& @>\be>> \ti{U} \\
        &@VVV&\\
        &\ti{V}_x &
        \end{CD}
        \end{eqnarray}
where $Q \ra \ti{V}_x$ is a principal $G$-bundle and $Q \ra \ti{U}$ is $G$-equivariant.
But since $\ti{V}_x$ is simply connected, $Q \ra \ti{V}_x$ has a section $s$. So, if  we define $f_{\ti{V}}=\be \circ s$,
we get a map  $f_{\ti{V}}:\ti{V}_x \ra \ti{U}$.
Let $h \in H_x$. If we see $h$ as an automorphism  of (\ref{ob}), then $F_x (h)$ is an automorphism of (\ref{imob}).
It follows that $f_{\ti{V}} \circ h = f_H(h) \circ f_{\ti{V}}$.

So, we have constructed a map $\mbox{Mor}_{\mbox{(Stacks)}}(\mcal{X},\mcal{Y}) \ra \mbox{Mor}_{\mbox{(Orb)}}([X], [Y])$.
We now show that it is an inverse of (\ref{morphismi orb > morphismi stacks}). Consider an object in $[\ti{V}_x/ H_x ](B)$. Let it be
described by the data $\{ h_{lm}, ~ \al_k \}$ with respect to the covering $\{B_k \}_{k\in K}$. 
Let us denote the image of this object under $F_x$ by $\{ g_{lm}, ~ \be_k \}$.
Then we have the following commutative diagram 
        \begin{eqnarray*}
        \begin{CD}      
        B_k @> \al_k >> \ti{V}_x  @>f_{\ti{V}}>> \ti{U}\\
        @VidVV @VVV @VVV \\
        B_k @>>> [\ti{V}_x/ H_x ] @>F_x >> [\ti{U}/G].
        \end{CD}
        \end{eqnarray*}
This shows that $\be_k = f_{\ti{V}} \circ \al_k$. Now, for any $l,m \in K$, the equation
$\al_l = h_{lm} \circ \al_m$ means that $h_{lm}$ is a morphism between two objects over 
the same base $B_{lm}$. Then, if we apply the functor $F_x$, we have that $F_x(h_{lm})$
is a morphism between the corresponding objects. This shows that $ g_{lm}= f_{H_x}(h_{lm})$. \qed

\section{Orbifold vector bundles}

We shall now  review the definition of orbifold vector bundles. We will follow the definition given in 
\cite{MP}. Then, following \cite{Sa}, we will review the construction of the 
\textit{orbifold tangent bundle}\index{orbifold tangent bundle}.
The orbifold cotangent  bundle and its exterior powers are constructed in the same way. Then we review the notion of 
\textit{differential form}\index{differential form} and also of the de Rham cohomology for orbifolds.
Finally we recall a Theorem, due to Satake, which states that the de Rham cohomology of $[Y]$ is isomorphic to 
the singular cohomology of the underlying topological space.

\begin{defn} \label{orbifold bundle definition}
Let $[Y]$ be an orbifold and let $\mcal{U}$ be the maximal atlas. A 
\textbf{smooth orbifold vector bundle}\index{orbifold vector bundle ! smooth}
$[E]$ on $[Y]$ is given by:
        \begin{enumerate}
        \item for any uniformizing system $(\ti{U}, G, \chi) \in \mcal{U}$ a (ordinary) vector bundle $E_{\ti{U}}$
                on $\ti{U}$;
        \item for each embedding $(\vphi, \lam) :(\ti{U}, G, \chi)\ra (\ti{U}', G', \chi')$ an isomorphism 
                $E(\vphi, \lam) : E_{\ti{U}} \ra \vphi^* E_{\ti{U}'}$, moreover we require that these isomorphisms 
                are functorial in $(\vphi, \lam)$.
        \end{enumerate}

Let $[E]$ be an orbifold vector bundle on $[Y]$. A \textbf{section}\index{orbifold vector bundle ! smooth section of}
$[s]$ of $[E]$ is given by a (ordinary) section $s_{\ti{U}}$ of $E_{\ti{U}}$, for any uniformizing system 
$(\ti{U}, G, \chi) \in \mcal{U}$,
such that, for each embedding $(\vphi, \lam) :(\ti{U}, G, \chi)\ra (\ti{U}', G', \chi')$, 
        $$
        E(\vphi, \lam)(s_{\ti{U}}) = s_{\ti{U}'}.
        $$      
\end{defn}

\vspace{0.2cm}

\begin{rem} \nf
To define an orbifold vector bundle $[E]$ on $[Y]$ (up to isomorphism), it is enough to specify 
the bundles $E_{\ti{U}}$ and vector bundle maps $E(\vphi, \lam)$ for all uniformizing systems $(\ti{U}, G, \chi)$
in some atlas  with the property that the images $U= \chi (\ti{U}) \subset Y$ form a basis for the topology
on $Y$. This is Remark 1, in Section 2, \cite{MP}.
\end{rem}

\begin{rem} \nf
Let $[E]$ be an orbifold vector bundle on $[Y]$, and let $(\ti{U}, G, \chi)$ be a uniformizing system.
Then each $g\in G$ defines an embedding $(g, Ad_g) :(\ti{U}, G, \chi) \ra (\ti{U}, G, \chi)$, where $Ad_g :G\ra G$
is given by $g' \mapsto g\circ g' \circ g^{-1}$. So, we have an isomorphism
        \begin{align*}
        E(g, Ad_g): E_{\ti{U}} \ra g^* E_{\ti{U}}.
        \end{align*}
This defines an action of $G$ on $E_{\ti{U}}$. Thus we see that $E_{\ti{U}}$ is a $G$-\textit{equivariant vector bundle} on $\ti{U}$.
\end{rem}

\vspace{0.2cm}

\begin{ex}[\nf{\textit{Orbifold tangent and cotangent bundle}}]\index{orbifold tangent bundle} \nf 
Let $[Y]$ be a smooth  orbifold, and let $\mcal{U}$ be the maximal compatible atlas. 
The \textit{orbifold tangent bundle} of $[Y]$ is the orbifold vector bundle $T_{[Y]}$ defined as follows:
        \begin{enumerate}
        \item $(T_{[Y]})_{\ti{U}} = T_{\ti{U}}$ is the tangent bundle of $\ti{U}$, for any uniformizing system $(\ti{U}, G, \chi)$;
        \item for any embedding $(\vphi, \lam):(\ti{U}, G, \chi)\ra (\ti{U}', G', \chi')$, 
                $T_{[Y]}(\vphi, \lam)= T\vphi$ is the tangent morphism of $\vphi$.
        \end{enumerate}

In the same way we define the \textit{orbifold cotangent bundle}\index{orbifold cotangent bundle} $T^{\vee}_{[Y]}$.
Then we can form the $p$-th exterior product $\wedge^p T^{\vee}_{[Y]}$.
\end{ex}

\vspace{0.5cm}

\begin{defn}
Let $[Y]$ be a smooth orbifold. A \textbf{differential $p$-form}\index{orbifold differential form} on $[Y]$ is a section
of $\wedge^p T^{\vee}_{[Y]}$. 

The space of differential $p$-forms over $[Y]$ will be denoted by  $\Om^p_{[Y]}$. 
\end{defn}

\vspace{0.2cm}

We can define the exterior differential $d:\Om^p_{[Y]}\ra \Om^{p+1}_{[Y]}$ and the wedge product on differential
forms (see \cite{Sa}, Section 3, for more details).
So, we define the $p$-th de Rham cohomology group of $[Y]$ in the usual way:
        \begin{align*}
        H^p_{dR}([Y], \R)=\fr{\rm{Ker}(d:\Om^p_{[Y]}\ra \Om^{p+1}_{[Y]})}{\rm{Im}(d:\Om^{p-1}_{[Y]}\ra \Om^{p}_{[Y]})}.
        \end{align*}

The following result holds, see \cite{Sa} Theorem 1, \cite{Kaw} pag. 78.

\begin{teo}
For any $p\geq 0$, there is a natural isomorphism between the $p$-th singular cohomology 
group $H^p(Y,\R)$ of the topological space $Y$ and the $p$-th 
de Rham cohomology group $H^p_{dR}([Y], \R)$ of the orbifold $[Y]$.
Moreover, under this isomorphism the exterior
product in $\Om^*_{[Y]}$ corresponds to the cup product in $H^p(Y,\R)$.
\end{teo}

\vspace{0.2cm}

We can define connections on orbifold bundles. For an orbifold vector bundle 
with a linear connection we have characteristic forms by Weil homomorphism.
The cohomology class of a characteristic form is independent of the choice 
of the connection. So we have  \textit{Euler classes} for oriented 
orbifold vector bundles,  \textit{Chern classes} for complex orbifold vector bundles,
and \textit{Pontrjagin classes} for real orbifold vector bundles.
Moreover one can see that these characteristic classes are defined over
the rational numbers.

\vspace{0.5cm}

Integration over compact orbifolds is defined as follows. 
First of all, assume that $[Y]=[V/G]$ is a global quotient orbifold. 
Let $\om$ be a  differential $p$-form  on $[Y]$. By definition,
$\om$ is a $G$-equivariant $p$-form $\tilde{\om}$ on $V$. 
Then the integration of $\om$ on $[Y]$ is defined by 
        \begin{align*}
        \int_{[Y]}^{orb} \om := \fr{1}{|G|} \int_V  \ti{\om},
        \end{align*}
where $|G|$ is the order of the group $G$. We use the convention such that, if the degree of the differential form
is different from the dimension of the manifold, then the integral is zero.

Let $[Y]$ be a compact orbifold. Fix a $C^{\infty}$ partition of unity $\{\rho_i \}$
subordinated to  the covering $\{U_i \}$, where each $U_i$ is an uniformized
open set in $Y$. For any $p$-form $\om$ on $[Y]$, the integration of $\om$ over $[Y]$
is defined by
        \begin{align*}
        \int_{[Y]}^{orb} \om := \sum_i  \int_{[U_i]}^{orb} \rho_i \cdot \om_{|U_i},
        \end{align*}
where $[U_i]$ has the global quotient structure. This definition
is independent of the choice of the partition of unity (\cite{Sa}, Section 8).

\vspace{0.2cm}

The following Theorem holds, which is a form of Poincar\'e Theorem for orbifolds, see also \cite{Sa}, Theorem 3.

\begin{teo}\index{Poincar\'e duality for orbifolds}
The bilinear form 
        \begin{eqnarray*}
        H^p_{dR}([Y], \R) & \times & H^{n-p}_{dR}([Y], \R) \ra \R \\
        (\om, \eta) & \mapsto & \int_{[Y]}^{orb} \om \wedge \eta
        \end{eqnarray*}
is non degenerate.
\end{teo}

\vspace{0.5cm}

We now recall the definition of \textit{complex orbifold vector bundle}\index{orbifold vector bundle ! complex} 
and of \textit{holomorphic section}\index{orbifold vector bundle ! holomorphic section of} over a complex orbifold.

\begin{defn}
Let $[Y]$ be a complex orbifold and let $\mcal{U}$ be the maximal atlas. A 
\textbf{complex orbifold vector bundle} $[E]$ on $[Y]$ is given by:
        \begin{enumerate}
        \item an ordinary complex vector bundle $E_{\ti{U}}$ on $\ti{U}$, for any uniformizing system 
                $(\ti{U}, G, \chi) \in \mcal{U}$;
        \item an isomorphism $E(\vphi, \lam) : E_{\ti{U}} \ra \vphi^* E_{\ti{U}'}$ of complex vector bundles, 
                for any embedding $(\vphi, \lam) :(\ti{U}, G, \chi)\ra (\ti{U}', G', \chi')$ such that
                these isomorphisms are functorial in $(\vphi, \lam)$.
        \end{enumerate}

Let $[E]$ be an orbifold vector bundle on $[Y]$. A \textbf{holomorphic section}
$[s]$ of $[E]$ is given by a ordinary holomorphic section $s_{\ti{U}}$ of $E_{\ti{U}}$, 
for any uniformizing system $(\ti{U}, G, \chi) \in \mcal{U}$,
such that, for each embedding $(\vphi, \lam) :(\ti{U}, G, \chi)\ra (\ti{U}', G', \chi')$, 
        $$
        E(\vphi, \lam)(s_{\ti{U}}) = s_{\ti{U}'}.
        $$      
\end{defn}

\chapter{Orbifold cohomology}

The notion of \textit{orbifold cohomology ring} was introduced by Chen and Ruan 
for an almost complex orbifold \cite{CR}. This has been extended by Fantechi and G\"ottsche  to a noncommutative ring \cite{FG},
in the case where the orbifold is a global quotient. Abramovich, Graber and Vistoli 
gave the definition of orbifold cohomology ring in the algebraic case, that is for a smooth Deligne-Mumford stack \cite{AGV}.

This chapter can be divided in two parts.  In the first part, Sections 1,2,3, 
we review the definition of orbifold cohomology ring for a complex orbifold. We follow closely the paper \cite{CR}. 
In the second part, Section 4 we state the \textit{cohomological crepant resolution conjecture}
given by Ruan \cite{Ru}. The aim of this paper is to verify this conjecture for a certain class of orbifolds 
which will be defined in the next chapter.

\vspace{0.2cm}

\begin{notation} \nf In this chapter all orbifolds will be complex orbifolds. So,  morphisms  will be holomorphic
and  orbifold vector bundles will be complex. 
Furthermore, we assume that $[Y]$ has complex dimension $\rm{dim}_{\C}[Y]=d$, see Definition \ref{orbifold complex dimension}.

We will use the same notation \ref{u.s. at y}, so, for any $y\in Y$,
$(\ti{U}_y, G_y, \chi_y)$ will denote a uniformizing system at $y$.
\end{notation}

\section{Inertia orbifold}

Let $[Y]$ be an orbifold. Let us consider the set
        $$
        Y_1=\{ (y, (g)_y):~ y\in Y,~ (g)_y \subset G_y ~ \mbox{is a conjugacy class} \},
        $$
as usual we will denote by $(g)_y$ the conjugacy class of $g\in G_y$. We want to define an orbifold $[Y_1]$
in such a way that there is a morphism $[\tau ]:[Y_1] \ra [Y]$ with  continuous function $\tau: Y_1 \ra Y$ is given by:
$(y, (g)_y) \mapsto y$. This orbifold will be called  the \textit{inertia orbifold}.

\vspace{0.5cm}

We introduce an equivalence relation on the set $Y_1$. Let $y \in Y$ and let  $(\ti{U}_y, G_y, \chi_y)$ 
be an uniformizing system at $y$ (\ref{u.s. at y}).
For any $y' \in U_y =\chi_y(\ti{U}_y)$, let $(\ti{U}_{y'}, G_{y'}, \chi_{y'})$ be a uniformizing system 
at $y'$ such that $U_{y'} \subset U_y$. Then we have an embedding 
$(\vphi , \lam ) :(\ti{U}_{y'}, G_{y'}, \chi_{y'}) \ra (\ti{U}_y, G_y, \chi_y)$ with
$\lam :G_{y'} \ra G_y$ injective, Remark \ref{rem3}. So, for any conjugacy class $(g')_{y'} \subset G_{y'}$, 
we can associate the conjugacy class of $\lam (g')$ in $G_y$, i.e. $(\lam (g') )_y \subset G_y$.
Notice that, the class $(\lam (g') )_y \subset G_y$ does not depend on the chosen embedding, Proposition \ref{pretransition}.
In this situation, we say that $(y', (g')_{y'})$ is \textit{equivalent} to $(y, (\lam (g') )_y)$.
This generate an equivalence relation on $Y_1$.

\vspace{0.2cm}

\begin{notation}\label{T} \nf We will denote by $T$ the set of equivalence classes of $Y_1$ just described.  
Elements of $T$  will be denoted by $(g)$. For any $(g)\in T$, the equivalence class $(g)$
 will be denoted by $Y_{(g)}$. The class of $(y, (1)_y) \in Y_1$ will be denoted by
$(1)$, so $Y_{(1)}=Y$.
\end{notation}

\vspace{0.5cm}

\begin{lem}[\nf{ \cite{CR} Lemma 3.1.1}]\label{inertia}
For any  $(g)\in T$, the set $Y_{(g)}$ has a topology and an orbifold structure $[Y_{(g)}]$ which is given as follows: 
for any point $(y, (g)_y)\in (Y_1)_{(g)}$, an uniformizing system at $(y, (g)_y)$ is given by 
        \begin{align*}
        (\ti{U}_y^{g}, C(g), {\chi_y}\mid),
        \end{align*}
where $g$ is a representative of $(g)_y$, $\ti{U}_y^{g}$ is the fixed-point set of $g$ in $\ti{U}$,
$C(g) \subset G_y$ is the centralizer of $g$ in $G_y$ and ${\chi_y}\mid$ denotes the restriction of $\chi_y$
to $\ti{U}_y^{g}$.

There is an orbifold morphism
        \begin{align*}
        [\tau_{(g)}]: [Y_{(g)}] \ra [Y]
        \end{align*}
which locally is given by the inclusion $\ti{U}_y^{g} \ra \ti{U}_y$ and the group injection $C(g) \ra G_y$. 

The orbifold structure $[Y_1]$ is defined to be the disjoint union of  $[Y_{(g)}]$, as $(g)$ varies in $T$, so that
        \begin{align*}
        [Y_1]= \bigsqcup_{(g)\in T} [Y_{(g)}].
        \end{align*}
There is a morphism  $[\tau ]: [Y_1] \ra [Y]$ defined by requiring that its  restriction to $[Y_{(g)}]$ is $[\tau_{(g)}]$.

Moreover, if  $[Y]$   is a complex orbifold, so is $[Y_1]$ and $[\tau]: [Y_1] \ra [Y]$ is holomorphic.
\end{lem}

\vspace{0.5cm}

\begin{defn}
Let $[Y]$ be an orbifold. The \textbf{inertia orbifold}\index{inertia orbifold} of $[Y]$
is the orbifold $[Y_1]$ descrbed in Lemma \ref{inertia}.

For any $(g)\in T$, $(g) \not= (1)$, the orbifold $[Y_{(g)}]$ is called a \textbf{twisted sector}\index{twisted sector}.
On the other hand the orbifold $[Y_{(1)}]$ is the \textbf{nontwisted sector}\index{nontwisted sector}.
\end{defn}

\vspace{0.5cm}

\begin{prop}\label{involuzione I}
Let $[Y]$ be an orbifold. Then, there is a morphism 
        \begin{align*}
        [I] : [Y_1] \ra [Y_1]
        \end{align*}
with continuous function $I: Y_1 \ra Y_1$ given  by $(y, (g)_y) \mapsto (y, (g^{-1})_y)$. It 
is an involution, i.e. $[I] \circ [I] = [id]$.
\end{prop}

\noindent \textbf{Proof.} For any $y\in Y$ and an uniformizing system at $y$, $(\ti{U}_y, G_y, \chi_y)$,
we get the following uniformizing system for $[Y_1]$:
        \begin{align*}
        \left( \bigsqcup_{(g)_y \subset G_y} \ti{U}_y^g,\quad  \op_{(g)_y \subset G_y} C(g),\quad
         \bigsqcup_{(g)_y \subset G_y} (\chi_y)\mid{\ti{U}_y^g}\right)
        \end{align*}
with the property that, if $g$ represents $(g)_y$, then $g^{-1}$ represents $(g^{-1})_y$.
Then, the restriction of $[I]$ to this uniformizing system is given by
        \begin{eqnarray*}
        \ti{U}_y^g & \ra & \ti{U}_y^{g^{-1}}, \quad  \ti x  \mapsto  \ti x \\
        C(g) & \ra & G(g^{-1}), \quad h \mapsto h. 
        \end{eqnarray*} \qed

\section{Orbifold cohomology group}

As vector space, $H^*_{orb}([Y])$ is the cohomology of the inertia orbifold $[Y_1]$. 
The grading takes the \textit{degree shifting}\index{degree shifting} 
of the elements of the local groups into account.

\vspace{0.5cm}

We review now the definition of  degree shifting. Let $y\in Y$ be any point and let 
$(\ti{U}_y,G_y, \chi_y)$ be a uniformizing system at $y$ (Definition \ref{u.s. at y def}). 
The origin $0$ of $\ti{U}_y$ is fixed by the action of $G_y$, so we have an action of $G_y$ on the tangent space of $\ti{U}_y$ at $0$. 
We represent this action by a group homomorphism $R_y: G_y \ra GL(d,\C)$, where $d=\rm{dim}_{\C}Y$. For every $g\in G_y$,
$R_y(g)$ can be written as a diagonal matrix (\cite{Se}, Chapter 2, Proposition 3):
        \begin{align*}
        R_y(g)=\rm{diag}\left( \rm{exp}(2\pi i m_{1,g}/m_g),...,\rm{exp}(2\pi i m_{d,g}/m_g)\right),
        \end{align*}
where $m_g$ is the order of $R_y(g)$, and $0\leq m_{i,g} < m_g$ is an integer. Since this matrix depends 
only on the conjugacy class $(g)_y$ of $g\in G_y$, we define a function $\iota :Y_1 \ra \Q$ by
        \begin{align*}
        \iota (y,(g)_y)= \sum_{i=1}^d \frac{m_{i,g}}{m_g}.
        \end{align*}

\begin{lem}
For any $(g)\in T$ (see Notation \ref{T}), the function $\iota :Y_{(g)} \ra \Q$ is  constant on each connected component.
\end{lem}

\noindent \textbf{Proof.} From the definition of $Y_{(g)}$ (see Notation \ref{T}), it follows that it is enough to prove
the following statement. Let $(y,(g)_y), (y',(g')_{y'}) \in Y_{(g)}$  such that  there exists an embedding 
$(\vphi, \lam) : (\ti{U}_{y'},G_{y'},\chi_{y'}) \ra (\ti{U}_{y},G_{y},\chi_{y})$ which sends the origin $0_{y'}\in \ti{U}_{y'}$ 
to the origin $0_y\in \ti{U}_{y}$, then $\iota (y,(g)_y) = \iota (y',(g')_{y'})$.

This follows from the fact that $\vphi$ is $\lam$-equivariant, so the tangent map 
$T\vphi :T_{0_{y'}}\ti{U}_{y'} \ra T_{0_y}\ti{U}_{y}$ is a morphism of representations. \qed

\vspace{0.5cm}

\begin{defn}\label{degree shifting definition}
For any $(g)\in T$, the \textbf{degree shifting number\index{degree shifting number} of} $(g)$ is the locally constant function
        \begin{align*}
        \iota (y,(g)_y) : Y_{(g)} \ra \Q.       
        \end{align*}
If $Y_{(g)}$ is connected, we identify $\iota (y,(g)_y)$ with its value, and we will denote it also by $\iota_{(g)}$.
\end{defn}

\vspace{0.5cm}

\begin{rem} \nf In the literature, the degree shifting number is also called the \textit{age}\index{age}.
\end{rem}

\begin{rem}\label{degree shifting as eigenvalue} \nf 
Note that, for any $(y, (g)_y) \in Y_{(g)}$, $\rm{exp} (2\pi i \iota (y,(g)_y))$ is the eigenvalue of the linear map 
        \begin{align*}
        \rm{det} R_y(g): \wedge^d \C^d \ra \wedge^d \C^d
        \end{align*}
\end{rem}

\vspace{0.5cm}

\begin{defn}
For any integer $p$, the degree $p$ \textbf{orbifold cohomology group}\index{orbifold cohomology group} of $[Y]$, 
$H^p_{orb}([Y])$,  is defined as follows
        \begin{align*}
        H^p_{orb}([Y])= \oplus_{(g)\in T} H^{p-2\iota_{(g)}}(Y_{(g)}),
        \end{align*}
where $H^*(Y_{(g)})$ is the singular cohomology of $Y_{(g)}$ with complex coefficients.
The \textbf{total orbifold cohomology group} of $[Y]$ is 
        \begin{align*}
        H^*_{orb}([Y])=\op_p H^p_{orb}([Y]).
        \end{align*}
\end{defn}

\vspace{0.5cm}

\begin{rem} \nf Note that $H^*_{orb}([Y])$ is a priori rationally graded. It is integrally graded if and only if all the degree 
shifting numbers are  half-integers. 
\end{rem}

\begin{rem} \nf From Remark \ref{degree shifting as eigenvalue} it follows that the degree shiftings are integers
if and only if the induced action of $G_y$ on $\wedge^d T_{0_y} \ti{U}_{0_y}$ is trivial. This means that 
the canonical sheaf of the singular variety $Y$ is locally free
which implies (and indeed is equivalent to say) that $Y$ is 
Gorenstein, see Definition \ref{gorenstein def}.
\end{rem}

\vspace{0.2cm}

\begin{rem} \nf The orbifold cohomology group $H^*_{orb}([Y])= \op_p H^p_{orb}([Y])$ can be split into 
even and odd parts, \cite{FG} Definition 1.8. By definition
        \begin{align*}
        H^{\rm{even}}_{orb}([Y]) = \oplus_{(g)\in T} H^{\rm{even}}(Y_{(g)}),
        \end{align*}
and analogously for the odd part. In general  this decomposition is not related with the even/odd 
decomposition given by the orbifold cohomology grading. On the other hand, for Gorenstein orbifolds,
the two gradings coincide.
\end{rem}

\vspace{0.5cm}

On the vector space $H^*_{orb}([Y])$ there is a complex valued pairing $\lan , \ran_{orb}$
which we will call the \textit{Poincar\'e duality}\index{Poincar\'e duality}. 

\begin{defn}
Let  $[I]:[Y_1] \ra [Y_1]$ be the holomorphic morphism defined in Proposition \ref{involuzione I}.
Then $[I]$ sends $[Y_{(g)}]$ to $[Y_{(g^{-1})}]$. 
The \textbf{Poincar\'e duality pairing} is the following pairing
        \begin{align*}
        \lan , \ran_{orb}:H^p_{orb}([Y])\times H^{2n-p}_{orb}([Y])\ra \C, \quad  \mbox{for}~ 0\leq p \leq 2n, 
        \end{align*}
defined as the direct sum of 
        \begin{align*}
        \lan , \ran_{orb}^{(g)}:H^{p-2\iota_{(g)}}(Y_{(g)})\times H^{2n -p-2\iota_{(g^{-1})}}(Y_{(g^{-1})})\ra \C,
        \end{align*}
where 
        $$
        \lan \al, \be \ran_{orb}^{(g)}=\int_{[Y_{(g)}]}^{orb} \al \cdot I^*(\be).
        $$
\end{defn}

\vspace{0.5cm}

We recall here Proposition 3.3.1 of \cite{CR}.

\begin{prop}
The Poincar\'e duality pairing is nondegenerate.
\end{prop}

\vspace{1cm}

\section{Orbifold cup product}

In this section we review the definition of \textit{orbifold cup product}. It is an
associative product on the total orbifold cohomology group, the resulting ring is the \textit{orbifold cohomology ring}.

\vspace{0.5cm}

For any positive integer $k$, consider the following set
        \begin{align*}
        Y_k=\{(y,(\ul{g})_y): ~y\in Y,~ \ul{g}=(g_1,...,g_k),~g_i\in G_y \},
        \end{align*}
where $(\ul{g})_y$ is the conjugacy class of $\ul{g}$. Here, two $k$-tuples $(g_1^{(1)},...,g_k^{(1)})$ and $(g_1^{(2)},...,g_k^{(2)})$
are conjugate if there exists  $g\in G_y$ such that $g_i^{(2)}=g\cd g_i^{(1)} \cd g^{-1}$ for all $i=1,...,k$.

\vspace{0.5cm}

We can define an equivalence relation on the set $Y_k$ as we did for $Y_1$, see Notation \ref{T}. 
The resulting set of equivalence classes will be denoted by $T_k$.

\vspace{0.2cm}

\begin{lem}[\nf{\cite{CR} Lemma 4.1.1}]\label{multisector}
For any $(\ul{g})\in T_k$, let $(Y_k)_{(\ul{g})}$ be the corresponding equivalence class of $Y_k$.
Then, there is a topology and an orbifold structure $[(Y_k)_{(\ul{g})}]$ over $(Y_k)_{(\ul{g})}$ such that, for any
point $(y, (\ul{g})_y)\in (Y_k)_{(\ul{g})}$, an uniformizing system at $(y, (\ul{g})_y)$ is given by 
        \begin{align*}
        (\ti{U}_y^{\ul{g}}, C(\ul{g}), {\chi_y}\mid),
        \end{align*}
where $\ul{g}=(g_1,...,g_k)$ is a representative of $(\ul{g})_y$, $\ti{U}_y^{\ul{g}}= \ti{U}_y^{g_1}\cap ... \cap  \ti{U}_y^{g_k}$,
and $C(\ul{g})=C(g_1)\cap ...\cap C(g_k)$.

We obtain an orbifold structure $[Y_k]$ defined as follows: 
        \begin{align*}
        [Y_k]= \bigsqcup_{(\ul{g})\in T_k} [Y_{(\ul{g})}].
        \end{align*}

For any $i=1,...,k$, we have morphisms
        \begin{align*}
        [e_i]:  [Y_k] \ra [Y_1]
        \end{align*}
that locally are given by: the topological inclusions $\ti{U}_y^{\ul{g}} \ra \ti{U}_y^{g_i}$ and 
the group injections $C(\ul{g}) \ra C(g_i)$.

Moreover, if  $[Y]$   is a complex orbifold, so is $[Y_k]$ and each $[e_i]: [Y_k] \ra [Y_1]$ is holomorphic.
\end{lem}

\vspace{0.2cm}

\begin{defn}
Let $[Y]$ be an orbifold. For any positive integer $k$, the \textbf{$k$-multisector}\index{multisector} of $[Y]$ is the orbifold
$[Y_k]$ constructed in Lemma \ref{multisector}.
\end{defn}

\vspace{0.2cm}

Consider the map $\rm{o}:T_k \ra T$ induced by $(g_1,...,g_k) \mapsto g_1\cd...\cd g_k$. The set $T_k^{\rm{o}}=\rm{o}^{-1}(1)$
is the subset of $T_k$ consisting of equivalence classes $(\ul{g})$ such that $\ul{g}=(g_1,...,g_k)$ satisfies the condition
$g_1\cd...\cd g_k=1$. 

\begin{notation} \nf 
We will denote by $Y_k^{\rm{o}}$ the set
        \begin{align*}
        Y_k^{\rm{o}} = \bigsqcup_{(\ul{g})\in T_k^{\rm{o}}Y_{(\ul{g})}},
        \end{align*}
and similarly, by $[Y_k^{\rm{o}}]$ the orbifold 
        \begin{align*}
        [Y_k^{\rm{o}}]=\bigsqcup_{(\ul{g})\in T_k^{\rm{o}}}[Y_{(\ul{g})}].
        \end{align*}
\end{notation}

\vspace{0.5cm}

The definition of the orbifold cup product requires the construction of an \textit{obstruction bundle}\index{obstruction bundle}
$[E_{(\ul{g})}]$ over each component $[Y_{(\ul{g})}]$ of $[Y_3^{\rm{o}}]$. We now review the definition of $[E_{(\ul{g})}]$.

\vspace{0.2cm}

Let $Y_{(\ul{g})}$ be a component of $Y_3^{\rm{o}}$, $(y,(\ul{g})_y) \in Y_{(\ul{g})}$ be a point and 
$(\ti{U}_y^{\ul{g}}, C(\ul{g}), {\chi_y}\mid )$ be a uniformizing system of $[Y_{(\ul{g})}]$ at $(y,(\ul{g})_y)$. 

If $\ul{g}=(g_1,g_2,g_3)$ is a  representative  of $(\ul{g})_y$, then $g_1 \cd g_2 \cd g_3 =1$. So we have
a morphism
        \begin{eqnarray*}
        \pi_1(S^2-\{0,1,\infty\}) &\ra & G_y \\
        \ga_i &\mapsto& g_i.
        \end{eqnarray*}
Here, $S^2$ is the unit sphere in $\R^3$, so $\pi_1(S^2-\{0,1,\infty\})$ is the free group generated
by three elements $\ga_1,\ga_2,\ga_3$ with the unique relation $\ga_1 \cd \ga_2 \cd \ga_3 =1$.
Geometrically, we can represent $\ga_1,\ga_2$ and $\ga_3$ as loops around $0,1$ and $\infty$ respectively.

There is a compact Riemann surface $\Si$ and a projection $\ti{\pi} :\Si \ra S^2$ which is a Galois covering, with Galois group 
the subgroup $G$ of $G_y$ generated by $g_1,g_2,g_3$, branched over $0,1,\infty$, see \cite{FG} Appendix.
In particular, the group $G$ acts on the vector space $H^1(\Si,\mcal{O}_{\Si})$.

\vspace{0.2cm}

Then we define a vector bundle over $\ti{U}_y^{\ul{g}}$ as follows:
        \begin{equation} \label{local obstruction for orbifold}
        (E_{(\ul{g})})_{\ti{U}_y^{\ul{g}}}:= \left( H^1(\Si, \mcal{O}_{\Si} ) \ot (T_{\ti{U}_y})\mid {\ti{U}_y^{\ul{g}}} \right)^G,
        \end{equation}
where $()^G$ means the  $G$-invariant part. Note that $H^1(\Si,\mcal{O}_{\Si})\ot (T_{\ti{U}_y})\mid {\ti{U}_y^{\ul{g}}}$
is a $G$-bundle, so $(E_{(\ul{g})})_{\ti{U}_y^{\ul{g}}}$ is a vector bundle. See \cite{FG}, page 201, for further details.

\vspace{0.2cm}

Let $(\psi, \mu):(\ti{U}_{y'}^{\ul{g}'}, C(\ul{g}'), {\chi_{y'}}\mid) \ra (\ti{U}_y^{\ul{g}}, C(\ul{g}), {\chi_y}\mid)$
be an embedding compatible with $[Y_{(\ul{g})}]$. Then we can suppose that it is induced by an embedding
$(\vphi, \lam): (\ti{U}_{y'}, G_{y'}, \chi_{y'} ) \ra (\ti{U}_y, G_y, \chi_y ) $ compatible with $[Y]$.
Then, if $G'$ denotes the subgroup of $G_{y'}$ generated by $\ul{g}'$, $\lam$ induces an isomorphism 
$G' \ra G$, where $G$ is the subgroup of $G_y$ generated by $\ul{g}$. This induces an isomorphism $\Si ' \ra \Si $
which is compatible with the group actions, where $\Si'\ra S^2 $ is the Galois covering associated to $\ul{g}'$. 
So, we obtain an isomorphism
        \begin{equation}\label{transition for obstruction}
        (E_{(\ul{g}')})_{\ti{U}_{y'}^{\ul{g}'}} \ra \psi^* (E_{ ( \ul{g} ) } )_{\ti{U}_y^{\ul{g}}}.
        \end{equation}

\vspace{0.2cm}

We have the following result.

\begin{prop}[\nf{\cite{CR} }] \label{obstruction for orbifold coho, construction}
The vector bundles given by (\ref{local obstruction for orbifold}), for any uniformizing system, and the isomorphisms
(\ref{transition for obstruction}) associated to any embedding of local uniformizing systems compatible with $[Y_{ (\ul{g}) } ]$,
satisfies the conditions of Definition \ref{orbifold bundle definition}. So these data defines an orbifold vector bundle 
on $[Y_3^{\rm{o}}]$.
\end{prop}

\vspace{0.2cm}

\begin{defn}\label{obstruction for orbifold coho}
The \textbf{obstruction bundle}\index{obstruction bundle  for orbifold cohomology}
for the orbifold cohomology of $[Y]$ is the orbifold vector bundle defined in 
Proposition \ref{obstruction for orbifold coho, construction}. 
\end{defn}

\vspace{0.2cm}

\begin{notation} \nf 
We will denote the obstruction bundle for the orbifold cohomology of $[Y]$ by
$[E]$. The restriction of  $[E]$ to the component $[Y_{(\ul{g})}]$, for any $(\ul{g}) \in T_3^{\rm{o}} $, will be denoted by
$[ E_{ ( \ul{g} ) } ]$.
\end{notation}

\vspace{0.5cm}

\begin{defn}\label{ocpdefn}
Let $[Y]$ be a complex orbifold such that $Y$ is compact.
For $\al, \be, \ga \in H^*_{orb}([Y])$, the \textbf{$3$-point function}\index{$3$-point function  for orbifold cohomology}
is defined as follows
        \begin{align*}
        \lan \al, \be, \ga \ran_{orb}= \sum_{(\ul{g})\in T_3^{\rm{o}}} \int_{[Y_{(\ul{g})}]}^{orb} [e_1]^* \al 
        \wedge [e_2]^* \be \wedge [e_3]^* \ga \wedge c_{top}([E_{(\ul{g})}]),
        \end{align*}
where  $[e_i]:[Y_{(\ul{g})}] \ra [Y_1]$, for $i=1,2,3$, is the morphism defined in Lemma \ref{multisector}.

The \textbf{orbifold cup product}\index{orbifold cup product} of $[Y]$ is the product on $H^*_{orb}([Y])$ 
        \begin{eqnarray*}
        \cup_{orb} : H^*_{orb}([Y]) \times H^*_{orb}([Y]) &\ra&  H^*_{orb}([Y])\\
        (\al, \be) & \mapsto &  \al \cup_{orb} \be
        \end{eqnarray*}
where $\al \cup_{orb} \be$ is defined by the following relation
        $$
        \lan \al \cup_{orb} \be , \ga \ran_{orb}=\lan\al, \be, \ga \ran_{orb} \quad \mbox{for any} ~ \ga \in H^*_{orb}([Y]).
        $$
\end{defn}

\vspace{0.2cm}

\begin{rem} \nf
The orbifold cup product can be defined also for an \textit{almost complex orbifold} that is not compact.
We will not recall the definition in this general case, it can be found in \cite{CR}, Definition 4.1.2.
\end{rem}

\vspace{0.5cm}

We will report in the following Theorem the most important properties of the orbifold cup product.
This is the main result of the paper \cite{CR}, see Theorem 4.1.5. Even if the Theorem holds for almost complex orbifolds 
which are not necessarily compact, we will present the result under stronger assumptions.

\begin{teo}\label{properties}
Let $[Y]$ be a complex orbifold such that the underlying topological space $Y$ is compact. Assume that $[Y]$ has complex dimension
$\rm{dim}_{\C} [Y] =d$. The orbifold cup product preserves the orbifold grading, i.e.
        $$
        \cup_{orb} :  H^p_{orb}([Y]) \times H^q_{orb}([Y]) \ra  H^{p+q}_{orb}([Y]),
        $$
for any $0\leq p,q \leq 2d$ such that $p+q \leq 2d$, and has the following properties.
        \begin{description}
        \item[\nf Associativity.] The  orbifold cup product is associative and has a unity $e_{[Y]}$.
                Moreover, $e_{[Y]} \in H^0_{orb}([Y])=H^0(Y)$ and it coincides with the unity of the usual cup product 
                of $Y$. 
        \item[\nf Poincar\'e duality.] For any $(\al, \be) \in H^p_{orb}([Y]) \times H^{2d -p}_{orb}([Y])$, with $0\leq p \leq 2d$,
                we have
                $$
                \int^{orb}_{[Y]} \al \cup_{orb} \be = \lan \al, \be \ran_{orb}.
                $$
        \item[\nf Deformation invariance.] The orbifold cup product $\cup_{orb}$ is invariant under deformations 
                of the complex structure of $[Y]$.
        \item[\nf Supercommutativity.] If $[Y]$ is Gorenstein, the total orbifold cohomology is integrally graded, and
                we have supercommutativity
                $$
                \al \cup_{orb} \be = (-1)^{\rm{deg}~ \al \cd \rm{deg} ~ \be} \be \cup_{orb} \al.
                $$
        \item[\nf  Compatibility with the usual cup product.] The restriction of $\cup_{orb}$ to the cohomology 
                of the nontwisted sector, i.e. $H^*(Y)$, is equal to the usual cup product of $Y$.
        \end{description}
\end{teo}

\section{Ruan's conjecture}

In this section we recall the statement of the \textit{cohomological crepant resolution conjecture}, given by Ruan in \cite{Ru}.
The conjecture gives a precise relation between the orbifold cohomology ring of a complex orbifold $[Y]$ and 
the cohomology ring of a crepant resolution of singularities of $Y$, when such a resolution  exists.

\begin{notation} \nf In this section $Y$ will be a complex algebraic variety. For an orbifold $[Y]$, 
we mean a complex orbifold structure over the topological space $Y$, where $Y$ has  the strong (or complex) topology,
see \cite{Mu} Chapter I, Section 10.
\end{notation}

\subsection{Crepant resolutions}

We first recall the definition of Gorenstein variety, Gorenstein orbifold, and  crepant resolution.
For more details see \cite{C3folds} and \cite{YPG}.

\begin{notation}\label{canonical sheaf def} \nf 
For any normal variety $Y$, we will denote by  $Y_0$  the smooth locus of $Y$ and by $l: Y_0 \ra Y$  the inclusion.
Then $K_Y$ will denote the sheaf $l_* K_{Y_0}$, where $ K_{Y_0}$ is the canonical sheaf of $Y_0$.
\end{notation}

\begin{defn}\label{gorenstein def}
$Y$ is \textbf{Gorenstein}\index{Gorenstein variety} if it is Cohen-Macaulay and $K_Y$ is locally free.

A \textbf{Gorenstein orbifold}\index{Gorenstein orbifold} is a complex orbifold structure $[Y]$ over a Gorenstein variety $Y$.
\end{defn}

\begin{defn}\label{crepant def}
Let $Y$ be a Gorenstein variety. A resolution of singularities $\rho:Z \ra Y$ is crepant if $\rho^*(K_Y)\cong K_Z$.
\end{defn}

\vspace{0.2cm}

\begin{rem}\label{crepant rem} \nf Crepant resolutions of Gorenstein varieties with quotient  singularities  
are known to exists in dimension $d=2,3$. 

In dimension $d=2$ the following stronger result holds, (\cite{BPV},  Chapter III, Theorem 6.2). 
\begin{teo}\label{surface resolution}
Every normal surface $Y$ admits a unique crepant resolution of singularities.
\end{teo} 

In dimension $d=3$ the existence of a crepant resolution is proven in \cite{Ro}, Main Theorem, pag. 493. 
However, in this case, the uniqueness result does not hold.

In dimension $d \geq 4$ crepant resolution exists only in rather special cases.
\end{rem}

\vspace{0.5cm}

\begin{ex} \nf (\textit{Hilbert scheme of points on surfaces}). 
An important class of examples for which a crepant resolution always exists, is the symmetric product of a compact complex surface $S$. 

Let $S$ be a compact complex surface. The $r$-th symmetric product of $S$ is the quotient of $S\times ...\times S$,
$r$-times, by the symmetric group $\mathfrak{S}_r$ acting by permutation. We will denote this quotient by $S^{(r)}$. 
Note that $S^{(r)}$ is a variety of dimension $2r$ with quotient singularities.
Now, let $S^{[r]}$ be the Hilbert scheme parametrizing $r$-points in $S$. 
Then, there is a morphism $\rho :S^{[r]} \ra S^{(r)}$ which is a crepant resolution. See \cite{Be}, Section 6, for a review.
\end{ex}

\begin{ex}\label{symplectic} \nf (\textit{Symplectic quotient singularities}). 
Let $V$ be a finite dimensional complex vector space equipped 
with a non degenerate symplectic form. Let $G \subset Sp(V)$ be a finite subgroup of the group $Sp(V)$ 
of symplectic automorphisms. The quotient $V/G$ has a natural structure of irreducible affine algebraic variety
with coordinate ring $\C[V/G]=\C[V]^G$, the subalgebra of $G$-invariants polynomials on $V$. Moreover 
the variety $V/G$ is  Gorenstein (\cite{Be1}, Proposition 2.4).

There are strong necessary conditions on $G$ in order for $V/G$ to have a crepant resolution (\cite{Ve}, Theorem 1.2).
Moreover there are examples of $G$ that do not match these conditions, \cite{GK} Theorem 1.1.
\end{ex}

\subsection{Quantum corrections}

We review the definition of \textit{quantum corrected cohomology ring} given by Ruan in \cite{Ru}.

\begin{notation}\label{3.4.2} \nf In this section $Y$ will be a Gorenstein projective algebraic variety.
So, for any crepant resolution  $\rho :Z \ra Y$, $Z$ will be a nonsingular projective algebraic variety.
\end{notation}

\vspace{0.5cm}

Let $[Y]$ be a Gorenstein orbifold and let $\rho :Z \ra Y$ be a crepant resolution. Consider the group homomorphism
        \begin{equation}\label{rho}
        \rho_* : H_2(Z, \Q) \ra H_2(Y, \Q)
        \end{equation}
induced by $\rho$. Choose a basis $\be_1,...,\be_n$ of $\mbox{Ker} ~ \rho_* \subset  H_2(Z, \Q)$ which consists of 
homology classes of rational curves.
We will call  $\be_1,...,\be_n$ an \textit{integral basis}\index{integral basis} of $\mbox{Ker} ~ \rho_*$.
Then,  the homology class of any effective curve that is contracted by $\rho$ can be written in a unique way as 
$\Ga=\sum_{l=1}^n a_l \be_l$, for $a_l\geq 0$. 

For each $\be_l$, we assign a formal variable $q_l$. Then $\Ga$ corresponds to $q_1^{a_1}\cd \cd \cd q_n^{a_n}$.
Define the \textit{quantum corrected $3$-point function}\index{quantum corrected $3$-point function} as 
        \begin{equation}\label{qc3}
        \lan \ga_1, \ga_2, \ga_3 \ran_{qc}(q_1,...,q_n):= \sum_{a_1,...,a_n\geq 0} 
        \Psi_{\Ga}^Z(\ga_1, \ga_2, \ga_3 )q_1^{a_1}\cd \cd \cd q_n^{a_n},
        \end{equation}
where $\ga_1,\ga_2,\ga_3\in H^*(Z)$ are cohomology classes, $\Ga =\sum_{l=1}^n a_l \be_l $,
and  $\Psi_{\Ga}^Z(\ga_1, \ga_2, \ga_3 )$ is the genus zero Gromov-Witten invariant, see (\ref{GWI})
for the definition.

\vspace{0.2cm} 

\begin{notation} \nf 
We assume that the quantum corrected $3$-point function is represented by an analytic function in the variables $q_1,...,q_n$
in some region of the complex space $\C^n$. We will denote this function by $\lan \ga_1, \ga_2, \ga_3 \ran_{qc}(q_1,...,q_n)$.
In the following, when  we  valuate this function on particular values of the $q_i$'s,  we will implicitly
understand that the analytic function is defined on such values.
\end{notation}

\vspace{0.2cm}

We now define a family of rings $H^*_{\rho}(Z)(q_1,...,q_n)$ depending on the parameters $q_1,...,q_n$,
where  $q_1,...,q_n$ belong to the domain of definition of the quantum corrected $3$-point function.

\begin{defn}\label{quantum corrected cup product}
The \textit{quantum corrected triple intersection}\index{quantum corrected triple intersection}  
$\lan \ga_1, \ga_2, \ga_3 \ran_{qc}(q_1,...,q_n)$ is defined as follows
        \begin{align*}
        \lan \ga_1, \ga_2, \ga_3 \ran_{\rho}(q_1,...,q_n)=\lan \ga_1, \ga_2, \ga_3 \ran +\lan \ga_1, \ga_2, \ga_3 \ran_{qc}(q_1,...,q_n),
        \end{align*}
where $\lan \ga_1, \ga_2, \ga_3 \ran :=\int_Z \ga_1 \cup \ga_2 \cup \ga_3$.
The \textit{quantum corrected cup product}\index{quantum corrected cup product}  $\ga_1 \ast_{\rho} \ga_2$ is defined by the equation
        \begin{align*}
        \lan \ga_1 \ast_{\rho} \ga_2, \ga \ran =\lan \ga_1 ,\ga_2, \ga \ran_{\rho}(q_1,...,q_n)
        \end{align*}
for arbitrary $\ga \in H^*(Z)$, where $\lan \ga_1 , \ga_2 \ran := \int_Z \ga_1 \cup \ga_2$. 
\end{defn}

\vspace{0.2cm}

\begin{rem} \nf  Note that the quantum corrected cup product is a family of products on $H^*(Z)$ depending on the parameters
$q_1,...,q_n$. These parameters belong to the domain of definition of the quantum corrected $3$-point function
$\lan \ga_1, \ga_2, \ga_3 \ran_{qc}(q_1,...,q_n)$.
\end{rem}

\begin{rem} \nf Our definition of quantum corrected triple intersection and of quantum corrected cup product
is slightly more general than that given by Ruan in \cite{Ru}. We can recover the definition given by Ruan by 
giving to the parameters the value $q_1=...=q_n =-1$, provided that this point belongs to the
domain of the quantum corrected $3$-point function.
\end{rem}

\vspace{0.5cm}

\begin{teo}\label{qccohomology}
For any $(q_1,...,q_n)$ belonging to the domain of the quantum corrected $3$-point function, the quantum corrected cup product
$\ast_{\rho}$ satisfies the following properties.
        \begin{description}
        \item[\nf{Associativity}.] It is an associative product on $H^*(Z)$, moreover it has an
                identity which coincides with the identity of the usual cup product of $Z$.
        \item[\nf{Supercommutativity}.] It is supercommutative, that is
                \begin{align*}
                \ga_1 \ast_{\rho} \ga_2 = (-1)^{\rm{deg}~\ga_1 \cd \rm{deg}~\ga_2} \ga_2\ast_{\rho}\ga_1 
                \end{align*}
                for any $\ga_1,\ga_2 \in H^*(Z)$.
        \item[\nf{Homogeneity}.] For any $\ga_1,\ga_2 \in H^*(Z)$, the following equality holds
                \begin{align*}
                deg~(\ga_1\ast_{\rho} \ga_2 )= \rm{deg}~\ga_1 + \rm{deg}~\ga_2 .
                \end{align*}
        \end{description}
For any $(q_1,...,q_n)$ as before, the resulting ring structure on $H^*(Z)$ given by $\ast_{\rho}$,
will be denoted by $H^*_{\rho}(Z)(q_1,...,q_n)$.
\end{teo}

\noindent \textbf{Proof.} Notice that the  definition of the quantum corrected cup product, $\ast_{\rho}$, is analogous to the definition
of the small quantum product for a smooth projective algebraic variety $Z$, as given for example in \cite{CK}, Definition 8.1.1.
(see also \cite{CK}, Proposition 8.1.6.). The only difference is in the set to which the effective curves belong.
Indeed, for the small quantum product, the quantum corrected $3$-point function is defined as a sum over all the effective curves in 
$Z$. In our definition, we take into account only effective curves that are contracted by $\rho$.

Let $B \subset H_2(Z, \Z)/\mbox{tor}$ be the set of homology classes $\be$ of effective curves in $Z$ such that 
$\rho_* (\be) =0$, where $\rho_*$ is defined in (\ref{rho}). 
Then, from Lemma \ref{B} below it follows that  the proof of associativity and supercommutativity is 
the same proof of Theorem 8.1.4. in \cite{CK} with $H_2(Z, \Z)$ replaced by $B$.

To prove the homogeneity property, notice that, for any $\be \in B$, the following equality holds
        $$
        \int_{\be} c_1(K_Z) =0.
        $$
Then, apply Proposition 8.1.5., \cite{CK}. \qed

\begin{lem}\label{B}
Let $B \subset H_2(Z, \Z)/\mbox{tor}$ be the set of homology classes $\be$ of effective curves in $Z$ such that 
$\rho_* (\be) =0$, where $\rho_*$ is defined in (\ref{rho}). Then $B$ satisfies the following properties:
        \begin{itemize}
        \item $B$ is a semigroup under addition and contains the zero of $ H_2(Z, \Z)/\mbox{tor}$;
        \item for any $\be \in B$, if $\be = \al_1 + \al_2$ with $\al_1, \al_2 \in H_2(Z, \Z)/\mbox{tor}$,
                then $\al_1, \al_2 \in B$.
        \end{itemize}
\end{lem}

\noindent \textbf{Proof.} The first condition is clear. So, we prove the second.

Let $\be = \al_1 + \al_2 \in B$, with $\al_1, \al_2 \in H_2(Z, \Z)/\mbox{tor}$. Then 
        $$
        0 = \rho_*(\be) = \rho_*(\al_1) + \rho_*( \al_2).
        $$
But, $\rho_*(\al_1)$ and $\rho_*( \al_2)$ are homology classes of effective curves in $Y$, and since $Y$ is 
assumed to be projective (Notation \ref{3.4.2}), $\rho_*(\al_1)=\rho_*( \al_2)=0$.

To see this, take any very ample line bundle $L$ on $Y$. Then $\int_{\rho_*(\al_i)} c_1(L) \geq 0$ for any $i=1,2$.
Since $\int_{\rho_*(\be)} c_1(L) =0$, it follows that $\int_{\rho_*(\al_1)} c_1(L)=\int_{\rho_*(\al_2)} c_1(L)  = 0$.
This implies that $\rho_*(\al_1) =\rho_*(\al_2)=0$. \qed

\subsection{The conjecture}

We review the statement of the \textit{cohomological crepant resolution conjecture}.

We use the same notations of the previous section.

\vspace{0.5cm}

The \textit{quantum corrected cohomology ring} of $Z$ is the ring  obtained from  $H^*(Z)(q_1,...,q_n)$
by giving to the variables $q_1,...,q_n$ the same value $-1$. It will be denoted by $H^*_{\rho}(Z)$, 
so 
        $$
        H^*_{\rho}(Z) = H^*(Z)(-1,...,-1).
        $$

\vspace{0.5cm}

\begin{conj}[\nf{\textit{Ruan}, \cite{Ru}}]\label{conj} The rings $H^*_{\rho}(Z)$ and $H^*_{orb}([Y])$ 
are isomorphic.
\end{conj}

\vspace{0.5cm}

\begin{rem} \nf The conjecture has been proved in the following cases:
        \begin{enumerate}
        \item for $Y=S^{(2)}$, be the symmetric product of a compact complex surface $S$, and $Z=S^{[2]}$
                be the Hilbert scheme of two points in $S$, in \cite{LQ};
        \item for $Y=S^{(r)}$, be the symmetric product of a complex projective surface $S$ with trivial canonical 
                bundle $K_S \cong \mcal{O}_S$, and $Z=S^{[r]}$ be the Hilbert scheme of $r$ points in $S$,
                in \cite{FG}, Theorem 3.10;
        \item for $Y=V/G$ be the quotient of  a  symplectic vector space of finite dimension by a finite group of symplectic automorphisms,
                and $Z$ is a crepant resolution, when it does exists, in \cite{GK}, Theorem 1.2.
        \end{enumerate}
Note that in cases 2 and 3, the  crepant resolution $Z$ has a holomorphic symplectic form, so the Gromov-Witten invariants
vanish and there are no quantum corrections. This means that the orbifold cohomology ring is isomorphic to
the ordinary cohomology ring of $Z$.
\end{rem}

\begin{rem} \nf The aim of this thesis is to verify this conjecture for orbifolds with transversal 
$ADE$-singularities, see Definition \ref{definition ADE arbifolds}. 
We will prove the conjecture in the case of transversal $A_1$-singularities in Chapter 6.1 and
for transversal $A_2$-singularities with minor modifications. Indeed, in this case,
the quantum corrected $3$-point function can not be evaluated in $q_1=q_2=-1$. However, we will show that by
giving to $q_1$ and $q_2$ other values (which we have found), the resulting ring $H^*_{\rho}(Z)(q_1,q_2)$
is isomorphic to the orbifold cohomology ring $H^*_{orb}([Y])$. This is the content of Theorem \ref{a2}.
\end{rem}

\chapter{Orbifolds with transversal $ADE$-singularities}

In this Chapter we define orbifolds with transversal $ADE$-singularities.
We describe the inertia orbifold in terms of the \textit{monodromy}\index{monodromy}, which we introduce in Section 4.
Finally, we compute the orbifold cohomology ring of orbifolds with transversal $A$-singularities and trivial 
monodromy.

Orbifolds with transversal $ADE$-singularities are generalizations of orbifolds associated to 
quotient surface singularities which are Gorenstein, also called \textit{rational double points}\index{rational double points}. 
So, in the first section we recall the definition of such surface singularities and collect some properties. 

\section{Rational double points}

Let $G\subset SL(2,\C)$ be a finite subgroup, $G\not= \{1 \}$. The inclusion $G \subset SL(2,\C)$ induces an action of $G$ on $\C^2$.
This action has $0\in \C^2$ as the only fixed point, and is free on $\C^2 \backslash \{0 \}$. The quotient  $\C^2 / G$
has a structure of algebraic variety whose ring of regular functions is $\C [u,v]^G$, the subring of 
the polynomial ring  $\C [u,v]$ consisting of functions which are invariants under the action of $G$.
Moreover, such quotient can be represented as an hypersurface $R$ in $\C^3$ passing through the origin $0\in \C^3$
and with $0\in \C^3$ as the only singular point (\cite{DV}, Chapter 5, Section 40).

\begin{defn}\label{rdp definition}
A \textbf{rational double point}\index{rational double point} is the germ of a surface singularity $R \subset \C^3$
which can be obtained as a quotient $\C^2 / G$ of $\C^2$ by a finite subgroup $G$ of $SL(2,\C)$.
\end{defn}

\vspace{0.2cm}

\begin{rem}\nf Rational double points are Gorenstein, see Example \ref{symplectic}.
\end{rem}

\vspace{0.5cm}

The finite subgroups of $SL(2,\C)$ can be classified in the following way. There is a group homomorphism 
$SL(2,\C) \ra PGL(2,\C)$ which is onto and two-to-one. Then, the finite subgroups of $SL(2,\C)$ are  inverse images of
finite subgroups of  $PGL(2,\C)$. 

By identifying the sphere $S^2$ with the complex projective line $\Pro^1$, we can see that
the symmetry groups of the five regular polyhedra are finite subgroups of $PGL(2,\C)$. These groups are called the 
\textit{polyhedral groups}\index{polyhedral groups}. Now, since the cube and the dodecahedron are duals respectively to 
the octahedron and icosahedron, their symmety groups are isomorphic (\cite{DV}, Chapter 2, Section 8). 
So,  the symmety groups of the regular polyhedra provides three finite subgroups of $PGL(2,\C)$.

The classification of finite subgroups of $PGL(2,\C)$ is given by the following theorem. 
The proof and other details can be found in \cite{DV}, Chapter 2, Section 10.

\begin{teo}\label{G}
Any finite subgroup of $PGL(2,\C)$ is conjugate to one of the following subgroups: 
the symmetry group of the tetrahedron, $\mb{E}_6$, of order $12$;
the symmetry group of the octahedron, $\mb{E}_7$, of order $24$; the symmetry group of the icosahedron, $\mb{E}_8$, of order $60$;
the dihedral group, $\mb{D}_n$ for $n\geq 4$, of order $4(n-2)$; the cyclic group, $\mb{A}_n$, of order $n+1$.
\end{teo}

\vspace{0.5cm}

To this classification of the finite subgroups of $SL(2,\C)$, corresponds a classification of rational double points.
We report now the classification of rational double points as hypersurfaces in $\C^3$ with coordinate $(x,y,z)$.
In the left column we report the group, while in the right the equation of the  corresponding  singularity.
        \begin{equation}
        \begin{matrix} \label{rdp}
        \mb{A}_n: &xy - z^{n+1} & \mb{for}~n\geq 1 \\
        \mb{D}_n: & x^2 + y^2 z +z^{n-1}  & \mb{for}~ n\geq 4\\
        \mb{E}_6: & x^2 + y^3 +z^4 &  \\
        \mb{E}_7: &x^2 + y^3 + yz^3  & \\
        \mb{E}_8: & x^2 + y^3 + z^5& .
        \end{matrix}
        \end{equation}
This is proved in \cite{DV}, Chapter 5, Section 39. 

\vspace{0.2cm}

\begin{rem}\nf It can be proved that rational double points are the only rational surface singularities.
For more details on this and on other characterizations of rational double points, see \cite{Dur}.
\end{rem}

\vspace{0.5cm}

\begin{rem}\label{resolution graph} \nf \textit{(Resolution graph).} 
As pointed out in Remark \ref{crepant rem}, for any rational double point $R$, 
there exists a unique crepant resolution
$\rho :\ti{R} \ra R$ (\cite{BPV}, Chapter III, Theorem 6.2). The exceptional locus of $\rho$ is the 
union of rational curves $C_1,..., C_n$ with selfintersection numbers equal to $-2$. 
Moreover, it is possible to associate a graph to the collection of these curves in the following way:
there is a vertex for any irreducible component of the exceptional locus; two vertices are joined by an edge if
and only if the corresponding components have non zero intersection. The list of the graphs obtained by resolving 
rational double points is given in \cite{Dur}, Table 1, and in \cite{BPV}, Chapter III, Proposition 3.6.
Each of this graph is called \textit{resolution graph}\index{resolution graph} of the corresponding rational double point.
\end{rem}

\begin{rem}\label{McKay} \nf \textit{(McKay correspondence).} McKay observed that the resolution graph of a rational double point
can be recovered from the representation theory of the corresponding subgroup $G \subset SL(2,\C)$ \cite{McKay}, see also
\cite{Cra}.

Let $G \subset SL(2,\C)$ be a finite subgroup and $Q=\C^2$ be the representation induced by the inclusion $G \subset SL(2,\C)$.
Let $\rho_0,...,\rho_m$ be the irreducible representations of $G$, with  $\rho_0$ being the trivial one. Then, for
any $j=1,...,m$ we can decompose $Q\ot \rho_j$ as follows
        $$
        Q\ot \rho_j = \op_{i=0}^m a_{ij} \rho_i, \quad a_{ij} = \rm{ dim}_{\C} Hom_{G}(\rho_i, Q\ot \rho_j ).
        $$
The \textit{McKay graph} of $G\subset SL(2,\C)$ is the graph with one vertex for any irreducible representation,
two vertices are joined by $a_{ij}$ arrows. This graph is denoted by $\ti{\Ga}_G$.

The following theorem holds, \cite{McKay}, see also \cite{Cra} Theorem 1.19.

\begin{teo}\label{McKay theo}
The McKay graph $\ti{\Ga}_G$ is an extended Dynkin graph of $\ti{A} \ti{D}\ti{E}$ type.
Moreover the subgraph $\Ga_G$ consisting of nontrivial representations is the resolution graph 
of the corresponding rational double point.
\end{teo}
\end{rem}

\section{Definition}

\begin{conv} \nf All the varieties are defined over the field $\C$ of complex numbers.
By an open subset of a variety we mean open in the strong topology, (see \cite{Mu}, Chapter I, Section 10 for the definition). 
We will specify when we use a Zariski-open subset.
\end{conv}

\begin{notation}\nf From now on, $R$  will denote a surface in $\C^3$ defined by one of the equations
(\ref{rdp}), i.e. a surface with a rational double point at the origin $0\in \C^3$. 
The crepant resolution of $R$ will be denoted by $\ti{R}$.
\end{notation}

\vspace{0.5cm}

Let $Y$ be a projective variety over $\C$. We say that $Y$ has \textit{transversal $ADE$-singularities} if 
the singular locus $S$ of $Y$ is connected, smooth, and  the couple $(S,Y)$ is locally 
isomorphic to $(\C^k\times \{0\}, \C^k \times R)$.

\vspace{0.2cm}

\begin{rem}\nf Let $Y$ be a $3$-fold with canonical singularities. Then, with the exception of at most 
a finite number of points, every point in $Y$ has an open neighbourhood which is nonsingular or
isomorphic to $\C \times R$, \cite{C3folds}, Corollary 1.14.
\end{rem}

\vspace{0.2cm}

The following Proposition is a particular case of the fact that every complex variety with quotient singularities
has a unique reduced  orbifold structure, \cite{S}, Theorem 1.3.

\begin{prop}\label{adeorb}
Let $Y$ be a variety with transversal $ADE$-singularities. Then there is a unique reduced complex holomorphic orbifold structure $[Y]$
on $Y$.
\end{prop}

\noindent \textbf{Proof.} The surface $R$ is isomorphic to the quotient of $\C^2$ by the action of a
finite subgroup $G$ of $SL(2,\C)$. 

For each point $y\in Y$, if $y\in Y-S$, we can take $U_y = \ti{U}_y =  C^{k+2}$ and $G_y=\{1\}$, otherwise
$y$ has a uniformized neighbourhood $U_y\cong \C^k \times R$ with uniformizing system $(\ti{U}_y\cong\C^k\times \C^2,G_y \cong G, \chi_y)$,
where $G_y$ acts trivially on the first factor of $\C^k\times \C^2$, while the action on the second factor is induced by the inclusion
$G_y \subset SL(2,\C)$. Then one can see that these charts patch together to give an orbifold structure  $[Y]$ over $Y$.

The uniqueness of the orbifold structure follows from \cite{S}, Theorem 1.3.  \qed

\vspace{0.5cm}

\begin{defn}\label{definition ADE arbifolds}
An \textbf{orbifold with transversal $ADE$-singularities}\index{orbifold with transversal $ADE$-singularities}
is the reduced orbifold $[Y]$ associated to a variety $Y$ with transversal $ADE$-singularities.
\end{defn}

\vspace{0.5cm}

\begin{rem}\nf An orbifold with transversal $ADE$-singularities is Gorenstein. This follows from the fact 
that rational double points are Gorenstein.
\end{rem}

\section{Inertia orbifold and monodromy}

We describe the inertia orbifold for orbifolds with transversal $ADE$-singularities. We will use the same notation introduced 
in Section 1 of Chapter 2, in particular Notation \ref{T}.

\vspace{0.2cm}

\begin{notation}\nf By a \textit{topological covering}\index{topological covering} we mean a covering space
as defined in \cite{Massey}, Chapter 5, Section 2, with the difference that we don't require the connectedness. 

Let $p: \ti{X} \ra X$ be a topological covering. For any point $x\in X$, the fundamental group $\pi_1 (X,x)$
of $X$ in $x$ acts on the fiber $p^{-1}(x)$ as defined in \cite{Massey} Chapter 5 Section 7.
We will call this action the \textit{monodromy} of the covering\index{monodromy ! of a topological covering}.
\end{notation} 

\vspace{0.2cm}

\begin{lem}\label{ADE twisted sectors}
Let $[Y]$ be an orbifold with transversal $ADE$-singularities, and  let 
        \begin{align*}
        \ti{Y}_1 :=  \bigsqcup_{(g)\in T,(g)\not= (1) } Y_{(g)}.
        \end{align*}
Then, the restriction of $\tau : Y_1 \ra Y$ to $\ti{Y}_1$ is a topological covering
        \begin{align*}
        \ti{\tau} :\ti{Y}_1 \ra S
        \end{align*}
and the connected components of $\ti{Y}_1$ are the topological spaces $Y_{(g)}$ of the twisted sectors
$[Y_{(g)}]$.
\end{lem}

\noindent \textbf{Proof.} Consider the following open cover of $S$, $\{ U_y \cap S :~  y\in S\}$. Then
        \begin{align*}
        \ti{\tau}^{-1} (U_y \cap S) = \bigsqcup_{(g)_y \subset G_y, (g)_y\not= (1)_y} \ti{U}_y^g,
        \end{align*}
moreover, the restriction of $\ti{\tau}$ to $\ti{U}_y^g$ is an homeomorphism for every $(g)_y \subset G_y$
such that $(g)_y\not= (1)_y$. \qed

\vspace{0.7cm}

\noindent We describe now the monodromy of $\ti{\tau} :\ti{Y}_1 \ra S$ explicitly. 

\vspace{0.2cm}

\begin{notation} \nf For any $y\in S$, the fiber $\tau^{-1}(y)$ is the set of conjugacy classes of  $G_y$. We will denote it by $T_y$.
\end{notation} 

\vspace{0.2cm}

Let $y\in S$ be a fixed point of $S$, $(\ti{U}_y, G_y, \chi_y)$ be a uniformizing system
at $y$ and  $U_y = \chi_y (\ti{U}_y)$. 

For any $y' \in  U_y  \cap S$ such that $U_{y'} \subset U_y$ there is an embedding 
$(\vphi , \lam ) : (\ti{U}_{y'}, G_{y'}, \chi_{y'}) \ra (\ti{U}_y, G_y, \chi_y)$, moreover $\lam : G_{y'} \ra G_y$
is an isomorphism. So, using $\lam$, we get a bijective correspondence $T_{y'} \ra T_y$ given by
$(g)_y' \mapsto (\lam (g))_y$, where $(g)_y' \in T_{y'}$ and $g$ is a representative of $(g)_y'$.
If   $(\psi, \mu ) : (\ti{U}_{y'}, G_{y'}, \chi_{y'}) \ra (\ti{U}_y, G_y, \chi_y)$ is another embedding, 
then $\mu = g \cd \lam \cd {g}^{-1}$ for some $g \in G_y$ (Proposition \ref{pretransition}), so the 
correspondence $T_{y'} \ra T_y$ does not depends on the  chosen embedding.

Let $[\al]\in \pi_1(S,y)$ be a class of a loop $\al$ based on $y$. Choose a finite number of points $y=y_0, y_1,...,y_k, y_{k+1}=y$
in $\al$   such that the sets $U_{y_i}$ form a cover of $\al$. The previous argument gives bijective maps
$T_{y_{i+1}} \ra T_{y_i}$ for any $i\in \{0,...,k\}$. The composition of these maps is a bijection $T_y \ra T_y$, i.e. an
element of the group $\rm{Aut}(T_y)$ of automorphisms of $T_y$.
The bijection associated to $\al$ depends  only on the homotopy class of $\al$ and  
the resulting map 
        \begin{align}\label{descrizione monodromia}
        \pi_1(S,y) \ra \rm{Aut}(T_y)
        \end{align}
is a group homomorphism.

Let $y_1 \in S$ be another point and $\ga$ a curve from $y_1$ to $y$. Let 
        \begin{align*}
        \pi_{\ga} :\pi_1 (S,y) \ra \pi_1 (S,y_1)
        \end{align*}
be the homomorphism: $[\al] \mapsto [\ga^{-1}\cd \al \cd \ga]$. In the same way we defined the homomorphism (\ref{descrizione monodromia}),
we can define a correspondence: $T_{\ga} : T_y \ra T_{y_1}$. We have a group homomorphism
        \begin{align*}
        A_{\ga} :\mb{Aut}(T_y) \ra \mb{Aut}(T_{y_1})
        \end{align*}
defined as follows: $a \mapsto T_{\ga} \circ a \circ T_{\ga}^{-1}$. It is easy to see that the following diagram is commutative,
        \begin{equation}\label{m}
        \begin{CD}
        \pi_1 (S,y) @>\mf{m}_y>> \mb{Aut}(T_y) \\
        @V\pi_{\ga}VV @VVA_{\ga}V\\
        \pi_1 (S,y_1) @>\mf{m}_{y_1}>> \mb{Aut}(T_{y_1}).
        \end{CD}
        \end{equation}

\vspace{0.5cm}

\begin{defn} 
Let $[Y]$ be an orbifold with transversal $ADE$-singularities, and let $y\in S$. The \textbf{monodromy}\index{monodromy} of $[Y]$
in $y$ is the group homomorphism (\ref{descrizione monodromia}) and we will denote it as follows:
        \begin{align*}
        \mathfrak{m}_y :\pi_1(S,y) \ra \rm{Aut}(T_y).
        \end{align*}
\end{defn}

\vspace{0.5cm}

\begin{lem}
For any $y \in S$, the following map is a bijection
        \begin{eqnarray*}
        T_y / \pi_1 (S,y) & \ra & T  \\
        \left[ {(g)}_y \right] & \mapsto & \left[ ( y, {(g)}_y ) \right],
        \end{eqnarray*}
where $T_y / \pi_1(S,y)$ is the quotient set of $T_y$ by the action of $\pi_1(S,y)$ given by the monodromy,
$T$ is the set defined in Notation \ref{T}.
\end{lem}

\noindent \textbf{Proof.} This follows from the fact that, for any path connected topological covering $p:\ti{X} \ra X$
and any $x\in X$, the action of $\pi_1(X,x)$ on the fiber $p^{-1}(x)$ is transitive. Notice that 
the restriction of $\ti{\tau}$ to $Y_{(g)}$ is a path connected topological covering of $S$ for all $(g) \in T$. \qed

\vspace{0.5cm}

By definition the monodromy of $[Y]$ in $y$ is a representation of $\pi_1(S,y)$ on the set $T_y$
of conjugacy classes of $G_y$. Identifying every conjugacy class with its characteristic function
we get a representation of $\pi_1(S,y)$ on the vector space of class functions on $G_y$, 
where as usual a class function is a complex valued function on $G$ that is constant on each conjugacy class.
According to the following Proposition, the set of  class functions on $G_y$ which are
characters of irreducible representations of $G_y$ is invariant under the action of  $\pi_1(S,y)$.

\begin{prop}\label{monodromy as aut gammag}
The set of characters of irreducible representations of $G_y$ 
is invariant under  the action of the monodromy $\mf{m}_y$ on the vector space of class functions on $G_y$.

Under the identification of any irreducible representation with its character we have that 
the image of $\mf{m}_y$ is a subgroup of the group of automorphisms of the graph $\Ga_{G_{y}}$.
\end{prop}

\noindent \textbf{Proof.} Let $y\in S$, $(\ti{U}_y, G_y, \chi_y)$ be a uniformizing system
at $y$ and $U_y = \chi_y (\ti{U}_y)$.
For any $y' \in  U_y  \cap S$ such that $U_{y'} \subset U_y \cap S$ there is an embedding 
$(\vphi , \lam ) : (\ti{U}_{y'}, G_{y'}, \chi_{y'}) \ra (\ti{U}_y, G_y, \chi_y)$. Moreover $\lam : G_{y'} \ra G_y$
is an isomorphism. We obtain a map which sends a representation of $G_y$ to a representation of $G_{y'}$ as follows.
Let $\mf{r} : G_y \ra GL(V)$ be a linear representation of $G_y$, then $\mf{r} \circ \lam :G_{y'} \ra GL(V)$
is a linear representation of $G_{y'} $. 

This map has the following properties:  1. is bijective; 2. the character of $\mf{r} \circ \lam$ is the composition 
of the character of $\mf{r}$ with $\lam$; 3. sends isomorphic representations of $G_y$ in 
isomorphic representations of $G_{y'}$; 4. sends irreducible representations of $G_y$ in irreducible representations of $G_{y'}$;
5. is compatible with the tensor product of two representations, that is 
$(\mf{r}_1 \ot \mf{r}_2)\circ \lam =\mf{r}_1\circ \lam \ot \mf{r}_2 \circ \lam$; 
6. sends the representation of $G_y$ on $N_{\ti{U}_{y}^{G_{y}} / \ti{U}_{y}, 0}$
in the representation of $G_{y'}$ in $N_{\ti{U}_{y'}^{G_{y'}} / \ti{U}_{y'}, 0}$. Where $N_{\ti{U}_{y}^{G_{y}} / \ti{U}_{y}, 0}$
(resp.  $N_{\ti{U}_{y'}^{G_{y'}} / \ti{U}_{y'}, 0}$) is the fiber at the origin $0\in \ti{U}_y$ (resp. $0\in \ti{U}_{y'}$)
of the normal vector bundle of $\ti{U}_{y}^{G_{y}}$ in $\ti{U}_{y}$ (resp. of $\ti{U}_{y'}^{G_{y'}}$ in $\ti{U}_{y'}$).
The last assertion is proved as follows. The tangent map 
        \begin{align*}
        T_0 \vphi : T_{\ti{U}_{y'}, 0} \ra T_{\ti{U}_y, 0}
        \end{align*}
induces a linear map
        \begin{align*}
        N_{\ti{U}_{y'}^{G_{y'}} / \ti{U}_{y'}, 0} \ra N_{\ti{U}_{y}^{G_{y}} / \ti{U}_{y}, 0}
        \end{align*}
such that the following diagram commutes for any $g' \in G_{y'}$
        \begin{equation}
        \begin{CD}
        N_{\ti{U}_{y'}^{G_{y'}} / \ti{U}_{y'}, 0} @>g'>> N_{\ti{U}_{y'}^{G_{y'}} / \ti{U}_{y'}, 0}  \\
        @VVV @VVV \\
        N_{\ti{U}_{y}^{G_{y}} / \ti{U}_{y}, 0} @> \lam (g') >> N_{\ti{U}_{y}^{G_{y}} / \ti{U}_{y}, 0}.
        \end{CD}
        \end{equation}

\vspace{0.2cm}

If  $(\psi, \mu ) : (\ti{U}_{y'}, G_{y'}, \chi_{y'}) \ra (\ti{U}_y, G_y, \chi_y)$ is another embedding, 
then $\mu = g \cd \lam \cd {g}^{-1}$, for some $g \in G_y$, Proposition \ref{pretransition}. So, for any representation of $G_y$,
the representations of $G_{y'}$ obtained using $\lam$ and $\mu$ are isomorphic. 

The previous arguments shows that we have a map from the set of isomorphism classes of representations of $G_{y}$ to the set 
of isomorphism classes of representations of $G_{y'}$ that satisfies properties 1., 2., 4., 5., 6. as before.
We will denote with  $\mf{n}_y$ the class of $N_{\ti{U}_{y}^{G_{y}} / \ti{U}_{y}, 0}$ and by $\mf{n}_{y'}$
the class  of $N_{\ti{U}_{y'}^{G_{y'}} / \ti{U}_{y'}, 0} $. 

The first assertion follows from property 2. To show the second claim, let $\mf{r_1},...,\mf{r_m}$ be the classes of irreducible 
representations of $G_{y}$, and let $\mf{r_1}',...,\mf{r_m}'$ the classes of irreducible representations of $G_{y'}$. 
Suppose that the correspondence sends $\mf{r_i}$ to $\mf{r_i}'$. Let $\mf{n}_y \ot \mf{r_i} = \sum_j a_{ji}\mf{r_j} $, 
and $\mf{n}_{y'} \ot \mf{r_i}' = \sum_j a_{ji}'\mf{r_j}' $. Then 
        \begin{eqnarray*}
        \lam \circ (\mf{n}_y \ot \mf{r_i}) & = & \sum_j a_{ji}\lam \circ \mf{r_j} \\
        & = & \sum_j a_{ji}\mf{r_j}',
        \end{eqnarray*}
on the other hand, 
        \begin{align*}
        \lam \circ (\mf{n}_y \ot \mf{r_i})= (\lam \circ \mf{n}_y ) \ot (\lam \circ \mf{r_i}) = \mf{n}_{y'} \ot \mf{r_i}'.
        \end{align*}
So, $a_{ji}=a_{ji}'$ for all $i,j$. \qed

\vspace{0.7cm}

\begin{rem}\label{con} \nf For $G= A_n , ~n \geq 1,~ D_n ~ n\geq 4, E_6, E_7,E_8$ (see Theorem \ref{G}), the automorphism group 
of $\Ga_G$ is given as follows:
        \begin{equation*}
        \begin{matrix} 
        G   & & \mb{Aut}(\Ga_G)\\
            & &          \\     
        A_1 & & \{1 \} \\
        A_n & n\geq 2& \Z_2 \\
        D_4 & & \mf{S}_3 \\
        D_n & n\geq 5 & \Z_2 \\
        E_6 & & \Z_2 \\
        E_7 & & \{ 1 \} \\
        E_8 & & \{ 1 \} 
        \end{matrix}
        \end{equation*}
where we have written on the left side the group $G$, and on the right $\mb{Aut}(\Ga_G)$. 
\end{rem}

\vspace{0.2cm}

The previous considerations give constraints on the topology of the spaces $Y_{(g)}$ for $(g)\in T$.
The following Corollary is an easy consequence of Proposition \ref{monodromy as aut gammag} and Remark \ref{con}.

\begin{cor}\label{monodromy and twisted sectors}
Let $[Y]$ be an orbifold with transversal singularities of type $A_1, E_7$ or $E_8$, then, for any $(g)\not= (1)$,
the topological space $Y_{(g)}$ is  isomorphic to $S$. 

Let $[Y]$ be an orbifold with transversal singularities of type $A_n$, for $n\geq 2$, or $D_n$, for $n\geq 5$, then, 
for any $(g)\not= (1)$, the topological space $Y_{(g)}$ is  isomorphic to $S$ if the monodromy is trivial, 
it is a  double covering of $S$ if the monodromy is not trivial.
\end{cor}

\vspace{0.2cm}

\begin{rem}\label{twisted sectors dipendenza locale}\nf 
Notice that, the twisted sectors $[Y_{(g)}]$ of $[Y]$ depends only on a neighbourhood of $S$ in $Y$.
Indeed, let $U\subset Y$ be an open neighbourhood of $S$ in $Y$, $U$ is a variety with transversal $ADE$-singularities.
Then, the twisted sectors $[U_{(g)}]$ of $[U]$ are canonically isomorphic to $[Y_{(g)}]$. So,
        \begin{align*}
        [Y_1] \cong [Y] \bigsqcup_{(g)\in T, (g)\not= (1)} [U_{(g)}].
        \end{align*}
\end{rem}

\vspace{0.5cm}

We describe now the inertia orbifold for an orbifold with transversal $A_n$-singularities.  We first study 
the case in which the monodromy is trivial. 

\vspace{0.2cm}

\begin{rem}\label{gruppi locali monodromia banale}\nf 
Let $[Y]$ be an orbifold with transversal $A_n$-singularities and trivial monodromy. Then we can 
identify all the local groups $G_y$, for $y\in S$, with the group of $(n+1)$-th roots of unity $\mu_{n+1}$
in a canonical way. Indeed, for any $y\in S$, consider the representation of $G_y$ on the normal space 
$N_{\ti{U}_y^{G_y}/\ti{U}_y , 0_y}$ of $\ti{U}_y^{G_y}$ in $\ti{U}_y $ over the point $0_y\in \ti{U}_y $.
It decomposes as direct sum of two representations
        \begin{align*}
        N_{\ti{U}_y^{G_y}/\ti{U}_y , 0_y} \cong \left( N_{\ti{U}_y^{G_y}/\ti{U}_y , 0_y} \right)^{\mf{g}_y}
        \op  \left( N_{\ti{U}_y^{G_y}/\ti{U}_y , 0_y} \right)^{\mf{g}^{-1}_y},
        \end{align*}
where $\left( N_{\ti{U}_y^{G_y}/\ti{U}_y , 0_y} \right)^{\mf{g}_y}$ 
(resp. $\left( N_{\ti{U}_y^{G_y}/\ti{U}_y , 0_y} \right)^{\mf{g}^{-1}_y}$)
is given by the character $\mf{g}_y$ (resp. $\mf{g}^{-1}_y$). Notice that 
        \begin{align*}
        \mf{g}_y,\mf{g}^{-1}_y : G_y \ra \mu_{n+1}
        \end{align*}
are group isomorphisms, so we can identify $G_y$ with $ \mu_{n+1}$ using $\mf{g}_y$ or $\mf{g}^{-1}_y$.
Since the monodromy is trivial, we can choose isomorphisms 
        \begin{align*}
        \mf{g}_y : G_y \ra \mu_{n+1}
        \end{align*}
such that the following diagram commutes for any $y'\in S$ with $U_{y'} \subset U_y$ and any 
embedding $(\phi, \lam): (\ti{U}_{y'}, G_{y'}, \chi_{y'}) \ra (\ti{U}_{y}, G_{y}, \chi_{y})$,
        \begin{equation*}
        \begin{CD}
        G_{y'} @>\lam>> G_y \\
        @V \mf{g}_{y'}VV @VV \mf{g}_{y}V \\
        \mu_{n+1} @>>=> \mu_{n+1}.
        \end{CD}
        \end{equation*}
\end{rem}

\vspace{0.5cm}

\begin{prop}
Let $[Y]$ be an orbifold with transversal $A_n$-singularities and  trivial monodromy.
Then, 
        \begin{align*}
        [Y_1] \cong [Y] \bigsqcup_{g\in \mu_{n+1}, g\not= 1}[S/\mu_{n+1}],
        \end{align*}
where the action of $\mu_{n+1}$ on $S$ is the trivial one. 
\end{prop}

\noindent \textbf{Proof.} Since the monodromy is trivial we can identify all the twisted sectors with $S$, see 
Corollary \ref{monodromy and twisted sectors}, moreover we identify the local groups $G_y$ with $\mu_{n+1}$,
for all $y\in S$, as explained in Remark \ref{gruppi locali monodromia banale}. These identifications 
give the isomorphism. \qed

\vspace{0.7cm}

Assume that the monodromy is not trivial. So, if $y\in Y$, we have a surjective homomorphism
        \begin{align}\label{monodromia non banale}
        \mf{m}_y : \pi_1(S,y) \ra \rm{Aut}(\Ga_{A_n})
        \end{align}
where $\rm{Aut}(\Ga_{A_n})$ is the cyclic group of order two.
Let $U\subset Y$ be a tubular neighbourhood of $S$ in $Y$, then (\ref{monodromia non banale})
gives a topological covering 
        \begin{align}\label{5ottobre}
        p: \ti{U} \ra U
        \end{align}
with monodromy $\mf{m}_y $. The topological space $\ti{U}$ has a structure of variety with transversal $A_n$-singularities
such that the orbifold $[\ti{U}]$ has trivial monodromy. 

The orbifold structure $[\ti{U}]$ can be described as follows. Let $U_y \subset U$ ba an open set
such that $p^{-1}(U_y) = V_1 \sqcup V_2$, where $\sqcup$ means disjoint union and the restrictions 
of $p$ to $V_1$ and $V_2$ are homeomorphisms. Then we choose uniformizing systems for 
$V_1$ and $V_2$ as follows: $(\ti{V}_1, G_1, \chi_1)$ is $(\ti{U}_y, G_y , \chi_y)$, 
$(\ti{V}_2, G_2, \chi_2)$ is $(\ti{U}_y, G_y , \chi_y)$ but the action of $G_y$ is given as follows
        \begin{align*}
        g\cd z := g^{-1}\cd z, \quad \mb{for} ~ g\in G_y, ~z\in \ti{U}_y.
        \end{align*}

These considerations lead to the following description of the twisted sectors $[Y_{(g)}]$.
Let $(g)\in T$ such that $(g)\not=(1)$, then 
        \begin{align*}
        Y_{(g)} \cong \left( \ti{U}_{(g)} \sqcup \ti{U}_{\e \cd (g)} \right)/\rm{Aut}(\Ga_{A_n}),
        \end{align*}
where the element $(g)$ on the left side is meant in the set $T$ relative to $[Y]$, instead, 
$(g)$ on the right side is the class of a representative $g$ of $(g)$ in the set $T$ relative to 
$[\ti{U}]$, the action of $\rm{Aut}(\Ga_{A_n})$ is given by 
        \begin{align}\label{azione twisted sectors caso con monodromia}
        \e \cd (z, (g)_z) = (\e \cd z, \e \cd ( g)_{\e\cd z}).
        \end{align}

\vspace{0.5cm}

The previous considerations are summarised in the following result.

\begin{prop}\label{twisted sectors with monodromy}
Let $[Y]$ be an orbifold with transversal $A_n$-singularities and non trivial monodromy.
Then, 
        \begin{align*}
        [Y_1] = [Y]\bigsqcup_{(g)\in T, (g)\not= (1)}[Y_{(g)}],
        \end{align*}
where $T \cong \Ga_{A_n}/\rm{Aut}(\Ga_{A_n})$ and 
        \begin{align*}
        Y_{(g)} \cong \left( \ti{U}_{(g)} \sqcup \ti{U}_{(\e \cd(g))} \right)/\rm{Aut}(\Ga_{A_n})
        \end{align*}
as topological spaces, where the action of $\rm{Aut}(\Ga_{A_n})$ is defined in (\ref{azione twisted sectors caso con monodromia}).
\end{prop}

\section{Orbifold cohomology ring}

We describe  now the orbifold cohomology ring of an orbifold $[Y]$ with transversal $A_n$-singularities. 
We first study the case where the monodromy is trivial.

\begin{conv} \nf If the monodromy is trivial, we identify  in a consistent way $G_y$ with $\mu_{n+1}$ for all $y\in S$ 
as explained in Reamrk \ref{gruppi locali monodromia banale}. We also identify 
$\mu_{n+1}$ with the group $\Z_{n+1}$ of integers modulo $n+1$ via the morphism:
        \begin{align*}
        a\in  \Z_{n+1}\mapsto \rm{exp}\left(\frac{2\pi i}{n+1}a\right).
        \end{align*}
Then, $g_1, g_2$ and $g_3$ in $\mu_{n+1}$ correspond  respectively to $a_1, a_2$ and $a_3$ in $\Z_{n+1}$.
We denote by $\ul{a}$ the vector $(a_1,a_2,a_3)$. 
\end{conv}

\vspace{0.5cm}

\begin{lem}\label{normal}
For any $y\in S$, let $(\ti{U}_y, G_y, \chi_y)$ be a uniformizing system of $[Y]$ at $y$. 
Then, the normal vector bundles
        \begin{align*}
        N_{\ti{U}_y^{G_y}/ \ti{U}_y} \ra \ti{U}_y^{G_y} \qquad \mb{for}~y\in S,
        \end{align*}
define an orbifold vector bundle $[N]$ over $[S/ \Z_{n+1}]$ of rank two,  
where the action of $\Z_{n+1}$ on $S$ is the trivial one.

If $n\geq 2$ and the monodromy is trivial, $[N]$ is isomorphic to the direct sum of two orbifold 
vector bundles, $[N]^{\mf{g}}$ and $[N]^{\mf{g}^{-1}}$, of rank one, that is 
        \begin{align*}
        [N] \cong  [N]^{\mf{g}} \op [N]^{\mf{g}^{-1}}.
        \end{align*}
\end{lem}

\noindent \textbf{Proof.} Let $y \in S$, and let $(\ti{U}_y, G_y, \chi_y)$ be a uniformizing system for $[Y]$ at $y$.
Then $(\ti{U}_y^{ G_y}, G_y, {\chi_y}\mid)$ is a uniformizing system for $[S/ \Z_{n+1}]$ at $y$, where ${\chi_y}\mid$
denotes the restriction of $\chi_y$ at $\ti{U}_y^{ G_y}$. Moreover, let $y' \in S$ and $(\ti{U}_{y'}, G_{y'}, \chi_{y'})$ 
be a uniformizing system at $y'$ such that $\chi_y(\ti{U}_y) \subset \chi_{y'}(\ti{U}_{y'})$. Then, any embedding  
$(\vphi, \lam) : (\ti{U}_y, G_y, \chi_y) \ra (\ti{U}_{y'}, G_{y'}, \chi_{y'})$ compatible with $[Y]$ induces an embedding 
$(\vphi \mid , \lam ) :(\ti{U}_y^{G_y}, G_y, {\chi_y}\mid) \ra (\ti{U}_{y'}^{G_{y'}}, G_{y'}, {\chi_{y'}}\mid)$ 
compatible with $[S/\Z_{n+1}]$. In this way, the orbifold $[Y]$ induces an orbifold structure on $S$, which coincides
with $[S/ \Z_{n+1}]$.

For any uniformizing system $(\ti{U}_y^{ G_y}, G_y, {\chi_y}\mid)$, we define 
        \begin{align*}
        N_{\ti{U}_y^{ G_y}} :=  N_{\ti{U}_y^{G_y}/ \ti{U}_y}.
        \end{align*}
For any embedding $(\vphi \mid , \lam ) :(\ti{U}_y^{G_y}, G_y, {\chi_y} \mid) \ra (\ti{U}_{y'}^{G_{y'}}, G_{y'}, {\chi_{y'}} \mid)$,
we get an isomorphism $N_{\ti{U}_y^{ G_y}} \ra {\vphi \mid}^* N_{\ti{U}_{y'}^{ G_{y'}}}$ which is induced by the tangent morphism
        \begin{align*}
        T\vphi : T_{ \ti{U}_y} \ra  T_{ \ti{U}_{y'}}.
        \end{align*}
These data define $[N]$.

Assume now that the monodromy is trivial and identify $\Z_{n+1}$ with $G_y$ for all $y\in S$. 
Notice that $\Z_{n+1}$ acts on  $N_{\ti{U}_y^{\Z_{n+1}}}$. So,
        \begin{equation}\label{dec}
        N_{\ti{U}_y^{\Z_{n+1} }} \cong \left( N_{\ti{U}_y^{\Z_{n+1}}} \right)^{\mf{g}} \op \left( N_{\ti{U}_y^{\Z_{n+1}}} 
\right)^{\mf{g}^{-1}},
        \end{equation}
where  $\mf{g} : \Z_{n+1} \ra \C^*$ is the character $a\mapsto \mb{exp}(\fr{2\pi ia}{n+1})$, $\mf{g}^{-1}$ is the dual of $\mf{g}$,
and $\left( N_{\ti{U}_y^{\Z_{n+1}}} \right)^{\mf{g}}$ (resp. $\left( N_{\ti{U}_y^{\Z_{n+1}}} \right)^{\mf{g}^{-1}}$) 
is the subbundle of $N_{\ti{U}_y^{\Z_{n+1}}}$ on which $\Z_{n+1}$ acts by the character $\mf{g}$ (resp. $\mf{g}^{-1}$). 
Since the monodromy is trivial, the line bundles $\left( N_{\ti{U}_y^{\Z_{n+1}}} \right)^{\mf{g}}$ (resp. 
$\left( N_{\ti{U}_y^{\Z_{n+1}}} \right)^{\mf{g}^{-1}}$) define the orbifold vector bundle $[N]^{\mf{g}}$ 
(resp. $[N]^{\mf{g}^{-1}}$). \qed

\vspace{0.5cm}

\begin{lem}\label{LM}\nf Let $n\geq 2$ and assume that the monodromy is trivial. Consider the natural morphism of orbifolds 
        \begin{align}\label{4ottobre}
        [S/\Z_{n+1}] \ra [S],
        \end{align} 
where the orbifold $[S]$ is the variety $S$, see Example \ref{smooth manifolds as orbifolds}. 
Then, there are line bundles $L$, $M$ and $K$ on $S$ which are defined by the following condition:
the pull-back of $L$, $M$ and $K$ under (\ref{4ottobre}) are respectively the line bundles 
$( [N]^{\mf{g}^{-1}} )^{\ot n+1}$, $\left( [N]^{\mf{g}} \right)^{\ot n+1}$ and $[N]^{\mf{g}^{-1}} \ot [N]^{\mf{g}}$. 
\end{lem}

\noindent \textbf{Proof.} We prove the statement about $[N]^{\mf{g}^{-1}}$ and $L$, the others are  similar.
Let $\{(\ti{U}_i, G_i, \chi_i )\}_{i\in I}$ be a set of uniformizing systems for $[Y]$ such that 
$\{( \ti{U}_i \cap S, G_i, {\chi_i}\mid )\}_{i\in I}$ is an orbifold atlas. We can assume that, for any $i, j \in I$,
we have isomorphisms
        \begin{align*}
        g_{ij} : \left( [N]^{\mf{g}^{-1}} \right)_{\ti{U}_{ij} \cap S} \ra \left( [N]^{\mf{g}^{-1}} \right)_{\ti{U}_{ji} \cap S}.
        \end{align*}
On triple intersections $\ti{U}_{ijk} \cap S$, the isomorphisms $g_{jk} \circ g_{ij}$ differ from $g_{ik}$ by
the action of an element of the group. So,
        $$
        g_{jk}^{n+1} \circ g_{ij}^{n+1}  = g_{ik}^{n+1}.
        $$
This proves the claim. \qed

\vspace{0.5cm}

\begin{teo} \label{orbifold product}
Let $[Y]$ be an orbifold with transversal $A_n$-singularities. Assume that the monodromy is trivial.
Then we have the following identification of vector spaces
        $$
        H^p_{orb}([Y])= H^p(Y_{(0)}) \oplus_{a\in \Z_{n+1} , a\not=0} H^{p-2}(Y_{(a)})
        $$
for all $p$. The orbifold cup product is given as follows, for $\al \in H^*( Y_{(a_1)})$ and  $\be \in H^*( Y_{(a_2)})$:
        \begin{enumerate}
        \item $  \al \cup_{orb} \be  =  \al \cup \be \in H^*(Y)$ if $a_1=a_2=0$
        \item $  \al \cup_{orb} \be  =  \al \cup  i^*(\be) \in H^*(Y_{(a_1)})$ if $a_1\not=0$, $a_2=0$
        \item $  \al \cup_{orb} \be  = \frac{1}{n+1}i_*(\al\cup \be)\in H^*(Y)$ if $a_1\not=0$, $a_2=a_1^{-1}$
        \item $  \al \cup_{orb} \be  = \fr{1}{n+1} \al\cup \be \cup c_1(L) \in H^*(Y_{(a_1+a_2)})$ if $a_1\not=0$, $a_2\not=0$, 
        $a_1+a_2<n+1$ in $\Z$
        \item $  \al \cup_{orb} \be  = \fr{1}{n+1}\al\cup \be \cup c_1(M) \in H^*(Y_{(a_1+a_2 -n+1)})$ if  $a_1\not=0$, $a_2\not=0$,
        $a_1+a_2>n+1$ in $\Z$,
\end{enumerate}

where $L$ and $M$ are the line bundles  defined in Lemma \ref{LM}, 
$i:S \ra Y$ is the inclusion of the singular locus in $Y$, and $\cup $ is the ordinary cup product of $Y$.

Note that, since $[Y]$ is Gorenstein, $\cup_{orb}$ is supercommutative, see Theorem \ref{properties}.
\end{teo}

\noindent \textbf{Proof.} The orbifold vector bundles $[E_{(\ul{a})}]$ have rank $0$ if $a_1=0$, $a_2=0$ or $a_3=0$.
This follows e.g. from \cite{FG}, Lemma 1.12,  where a relation between the rank of $[E_{(\ul{a})}]$  
and the degree shifting numbers of $a_1, a_2$ and $a_3$ is proved.

Let $y\in S$ and $\ul{a}= (a_1, a_2, a_3 )\not= (0,0,0)$. By definition (see equation \ref{local obstruction for orbifold})
we have
        \begin{align*}
        (E_{(\ul{a})})_{\ti{U}_y^{\ul{a}}}:= \left( H^1(\Si, \mcal{O}_{\Si} ) \ot (T_{\ti{U}_y})\mid {\ti{U}_y^{\ul{a}}} \right)^{G},
        \end{align*}
where $G\subset \Z_{n+1}$ is the subgroup generated by $a_1,a_2,a_3$, $\Si \ra \Pro^1$
is the Galois cover with Galois group $G$ branched over $0,1,\infty \in  \Pro^1$ and monodromy respectively $a_1,a_2,a_3$. 

We first notice that we can replace $(T_{\ti{U}_y})\mid {\ti{U}_y^{\ul{a}}}$ with  the normal bundle 
$N_{\ti{U}_y^{\ul{a}}/\ti{U}_y }$. Then, we can replace $\Si$ with the Galois cover $C\ra \Pro^1$,
induced by the inclusion $G\subset \Z_{n+1}$, which has  Galois group $\Z_{n+1}$. So   we have:
        \begin{align*}
        (E_{(\ul{a})})_{\ti{U}_y^{\ul{a}}} \cong 
        \left( H^1(C, \mcal{O}_{C})\ot  N_{\ti{U}_y^{\ul{a}}/\ti{U}_y }\right)^{\Z_{n+1}}.
        \end{align*}

Note that $p:C\ra \Pro^1$ is an abelian cover in the sense of \cite{P}, so
        $$
        p_* \mcal{O}_C=\op_{\mf{c} \in \Z_{n+1}^*} (L^{-1})^{\mf{c}}
        $$
where $\Z_{n+1}^*$ is the group of characters of $\Z_{n+1}$ and $\Z_{n+1}$ acts on $(L^{-1})^{\mf{c}}$ via the character $\mf{c}$.
Note that the characters $\mf{g}$ and $\mf{g}^{-1}$ defined in Lemma \ref{normal} are generators of $\Z_{n+1}^*$.

Using the fact $H^1(C, \mcal{O}_C)\cong H^1(\Pro^1,p_* \mcal{O}_C)$ and the decomposition (\ref{dec}), we have
        \begin{eqnarray*}
        & & \left( H^1(C, \mcal{O}_C)\ot N_{\ti{U}_y^{\ul{a}}/\ti{U}_y} \right)^{\Z_{n+1}}  \cong  \\
        & & H^1(\Pro^1, (L^{-1})^{\mf{g}})\ot  \left( N_{\ti{U}_y^{\ul{a}}/\ti{U}_y} \right)^{\mf{g}^{-1}}
        \op H^1(\Pro^1, (L^{-1})^{\mf{g}^{-1}})\ot  \left( N_{\ti{U}_y^{\ul{a}}/\ti{U}_y} \right)^{\mf{g}}.
        \end{eqnarray*}

By Proposition 2.1 and in particular by example 2.1 i) of \cite{P} we have that 
        \[ L^{\mf{g}} =
        \begin{cases}
        \mcal{O}(2) & \text{ if $a_1+a_2<n+1$},\\
        \mcal{O}(1) & \text{ if $a_1+a_2 \geq n+1$}
        \end{cases} \]
and 
                \[ L^{\mf{g}^{-1}}=
        \begin{cases}
        \mcal{O}(1) & \text{ if $a_1+a_2 \leq n+1$},\\
        \mcal{O}(2) & \text{ if $a_1+a_2 >  n+1$}.
        \end{cases} \]

It follows that 
        \[ {(E_{(\ul{a})})_{\ti{U}_y^{\ul{a}}}} =
        \begin{cases}
        \left( N_{\ti{U}_y^{\ul{a}}/\ti{U}_y} \right)^{\mf{g}^{-1}} & \quad \mbox{if} \quad a_1 +a_2 <n+1 \\
        \left( N_{\ti{U}_y^{\ul{a}}/\ti{U}_y} \right)^{\mf{g}}\quad & \quad \mbox{if} \quad a_1 +a_2 >n+1 
        \end{cases} \]
and the obstruction bundle is 
        \[ {[E_{(\ul{a})}]} =
        \begin{cases}
        [N]^{\mf{g}^{-1}}  & \quad \mbox{if} \quad a_1 +a_2 <n+1 \\
        [N]^{\mf{g}}       & \quad \mbox{if} \quad a_1 +a_2 >n+1 
        \end{cases} \] 
\qed

\vspace{0.7cm}

We now study the case in which the monodromy is not trivial. We have seen that the monodromy (\ref{dec})
gives a topological covering $p:\ti{U} \ra U$ of an open neighbourhood $U$ of $S$ in $Y$.
Moreover $\ti{U}$ is a variety with transversal $A_n$-singularities, so we have a natural orbifold structure 
$[\ti{U}]$ on $\ti{U}$. Let $\ti{S}$ be the singular locus of $\ti{U}$. 
The orbifold $[\ti{U}]$ has trivial monodromy, so we can identify all the local groups $G_{\ti{y}}$ with $\Z_{n+1}$, 
where $\ti{y} \in \ti{S}$, see Remark \ref{gruppi locali monodromia banale}. 
Under the identification of $\ti{S}$ with any twisted sector $\ti{U}_{(g)}$, 
we have the following presentation for the orbifold cohomology of $[\ti{U}]$:
        \begin{align*}
        H^*_{orb}([\ti{U}]) \cong H^*(\ti{U}) \op_{a\in \Z_{n+1}, a\not= 0} H^{*-2}(\ti{S}).
        \end{align*}
This is an easy consequence of Theorem \ref{orbifold product}. 

Notice that the group $\rm{Aut}(\Ga_{A_n})$ acts on $H^*_{orb}([\ti{U}])$. Given $\e \in \rm{Aut}(\Ga_{A_n})$
and $\ti{\de} + \ti{\al}_a \in H^*_{orb}([\ti{U}])$, then we have
        \begin{align}\label{z2action}
        \e \cd (\ti{\de} + \ti{\al}_a) := (\e^{-1})^* \ti{\de} + ((\e^{-1})^* \ti{\al})_{-a},
        \end{align}
where $\ti{\de}$ denotes an element of $H^*(\ti{U})$ and $\ti{\al}_a$ denotes the element of 
$\op_{a\in \Z_{n+1}, a\not= 0} H^{*-2}(\ti{S})$ given by the $\ti{\al} \in  H^{*-2}(\ti{S})$
in the $a$-th addendum.

It follows from Proposition \ref{twisted sectors with monodromy} that there is a natural isomorphism 
of vector spaces:
        \begin{align}\label{orbcohowithmon}
        H^*_{orb}([U]) \cong \left( H^*_{orb}([\ti{U}])\right)^{\rm{Aut}(\Ga_{A_n})},
        \end{align}
where the right hand side is the $\rm{Aut}(\Ga_{A_n})$-invariant subspace of $H^*_{orb}([\ti{U}])$.

Even if $\ti{U}$ and $U$ are not compact, it is possible to define
an orbifold cup product $\cup_{orb}$ on $H^*_{orb}([\ti{U}])$ (resp. on $H^*_{orb}([U])$)  such that the resulting ring 
$(H^*_{orb}([\ti{U}]),\cup_{orb})$ (resp. $(H^*_{orb}([U]),\cup_{orb})$ has the properties listed in 
Theorem \ref{properties}, see \cite{CR} Definition 4.1.2. This product is defined in an analogous way 
to the definition of the orbifold cup product for compact orbifolds, see Definition \ref{ocpdefn}.
In particular the obstruction bundle is as defined in Definition \ref{obstruction for orbifold coho}.

\vspace{0.5cm}

\begin{prop}
The restriction of the orbifold cup product on $H^*_{orb}([\ti{U}])$ to the subspace 
$\left( H^*_{orb}([\ti{U}])\right)^{\rm{Aut}(\Ga_{A_n})}$ is an associative product. 
The map (\ref{orbcohowithmon}) is then a ring isomorphism.
\end{prop}

\noindent \textbf{Proof.} Let $\ul{a}=(a_1,a_2,a_3) \in (\Z_{n+1})^3$ with $a_1 +a_2 +a_3 =0$ and $\ul{a}\not= (0,0,0)$.
Then consider  the  morphism
        \begin{align}\label{nspcc}
        [\ti{U}_{-\ul{a}}]  \ra  [\ti{U}_{\ul{a}}]
        \end{align}
whose associated continuous function $\ul{U}_{-\ul{a}} \ra \ul{U}_{\ul{a}}$ is given by $\ti{y} \ra \e\cd \ti{y}$,
and all the morphisms between the local groups $\Z_{n+1}$ are $g \ra -g$, 
where $\e \in \rm{Aut}(\Ga_{A_n})$, $\e \not= 1$, and we have identified each $\ul{U}_{\ul{a}}$ with $\ti{S}$, so
$\e\cd \ti{y}$ is the monodromy action.

From the description of the obstruction bundle given in the proof of Theorem \ref{orbifold product},
we easily get that the pull-back of $[E_{\ul{a}}]$ under (\ref{nspcc}) is isomorphic to $[E_{-\ul{a}}]$,
for any $\ul{a}=(a_1,a_2,a_3) \in (\Z_{n+1})^3$ with $a_1 +a_2 +a_3 =0$. This prove the first claim.

The fact that (\ref{orbcohowithmon}) is a ring isomorphism follows from a direct computation. \qed

\section{Examples}

We give here some special examples of orbifold cohomology rings.

\vspace{0.5cm}

\begin{ex}\nf \textit{(Surface case).} We give now a description of the orbifold cohomology of a surface 
        with an $A_n$-singularity:
        \begin{align*}
        Y=\{ (x,y,z)\in \C^3 : xy - z^{n+1}=0\}.
        \end{align*}
        $Y$ is the quotient of $\C^2$ by the action of the group $\mu_{n+1}$
        given by $\e \cd (u,v)= (\e\cd u, \e^{-1} \cd v)$, $\e \in \mu_{n+1}$. \\
        As a vector space
                $$
                H^*_{orb}([Y])=H^*(Y)\op H^{*-2}(S)\lan e_1 \ran \op ...\op H^{*-2}(S)\lan e_n \ran
                $$
        where $e_i$ is a generator of $H^*(Y_{(i)})$ as $H^*(S)$-module. \\
        The product rule is given by
        \[ e_i\cup_{orb}e_j=
        \begin{cases}
        0 & \text{ if $i+j \not=0 (\mb{mod}~ n+1)$},\\
        \frac{1}{n+1}i_*[S]\in H^4(Y)& \text{ if $i+j =0 (\mb{mod}~ n+1)$}.
        \end{cases} \]
\end{ex}

\begin{ex}\label{orba1}\nf \textit{(Transversal $A_1$-case).} 
                $$
                H^*_{orb}([Y])=H^*(Y)\op H^{*-2}(S)
                $$
as vector space. Given $(\de_1, \al_1), (\de_2, \al_2)\in H^*_{orb}(Y)$, we have the following espression for the 
orbifold cup product:
                $$
                (\de_1, \al_1) \cup_{orb} (\de_2, \al_2)=(\de_1 \cup \de_2 +\frac{1}{2} i_*(\al_1 \cup \al_2),
                i^*(\de_1)\cup \al_2 +  \al_1 \cup i^*(\de_2))
                $$
Note that in this case  the obstruction bundle $[E]$ has rank zero.
\end{ex}

\begin{ex}\label{orba2}\nf \textit{(Transversal $A_2$-case).} 
                $$
                H^*_{orb}([Y])=H^*(Y)\op H^{*-2}(S) \op H^{*-2}(S)
                $$
as vector space. Given $(\de_1, \al_1, \be_1), (\de_2, \al_2, \be_2)\in H^*_{orb}(Y)$, we have the following espression for the 
        orbifold cup product:
                \begin{eqnarray*}
                (\de_1, \al_1, \be_1)&\cup_{orb}& (\de_2, \al_2, \be_2)=
                (\de_1 \cup \de_2 +\frac{1}{2} i_*(\al_1 \cup \be_2 +\be_1\cup \al_2),\\
                & & i^*(\de_1)\cup \al_2 +  \al_1 \cup i^*(\de_2)
                +\be_1 \cup \be_2 \cup c_1(L), \\
                & & i^*(\de_1)\cup \be_2  +  \be_1 \cup i^*(\de_2) +
                \al_1 \cup \al_2 \cup c_1(M)).
                \end{eqnarray*}
\end{ex}

\chapter{Crepant resolutions}

In this Chapter we show that any variety with transversal $ADE$-singularities $Y$ (see Chapter 3.2) 
has a unique crepant resolution $\rho :Z\ra Y$. Then we restrict our attention to the case of 
transversal $A_n$-singularities and trivial monodromy. In this case we describe the exceptional locus $E$
in terms of the line bundles $L$ and $M$ defined in Lemma \ref{LM}.
We compute the cohomology ring $H^*(Z)$ of $Z$ in terms of the  cohomology of $Y$ and of $E$.

In the first section we recall some facts about the resolution of rational double points.

\section{Crepant resolutions of rational double points}

Let $R\subset \C^3$ be a rational double point. Then, by Theorem \ref{surface resolution}, $R$ has a 
unique crepant resolution $\rho :\ti{R} \ra R$. $\ti{R}$ can be obtained  by blowing-up successively the
singular locus. The exceptional locus $C\subset \ti{R}$ is union of rational curves $C_l$ whose autointersection
numbers are $C_l \cd C_l = -2$.
The shape of $C$ inside $\ti{R}$ is described by the resolution graph, see Remark \ref{resolution graph}.
We explain with the next example the $A_n$ case.

\begin{ex}\nf \textit{(Resolution of $A_n$-surface singulerities).}
Let 
        $$
        R=\{ (x,y,z)\in \C^3 :~ xy - z^{n+1} =0 \}
        $$
be a surface singularity of type $A_n$. Let $r: R_1 = Bl_{0} R \ra R $ be the blow-up of $R$ at the origin.
Then $R_1$ is covered by three open affine varieties $U,V$ and $W$, where
\begin{eqnarray*}
U & = & \{ \left(x, \fr{v}{u}, \fr{w}{u}\right) \in \C^3 :~  \left(\fr{v}{u}\right) - x^{n-1}\left(\fr{w}{u}\right)^{n+1}=0 \}\\
V & = & \{ \left(y, \fr{u}{v}, \fr{w}{v}\right) \in \C^3 :~  \left(\fr{u}{v}\right) - y^{n-1}\left(\fr{w}{v}\right)^{n+1}=0 \}\\
W & = & \{ \left(z, \fr{u}{w}, \fr{v}{w}\right) \in \C^3 :~  \fr{u}{w}\fr{v}{w} - z^{n-1} =0 \}.
        \end{eqnarray*}
and the restruction of $r $ to $U,V,W$ is given as follows      
        \begin{eqnarray*}
        r_{|U}: \left(x, \fr{v}{u}, \fr{w}{u}\right) & \mapsto & (x,y,z)= \left(x, x\fr{v}{u}, x\fr{w}{u}\right) \\
        r_{|V}: \left(y, \fr{u}{v}, \fr{w}{v}\right) & \mapsto & (x,y,z)= \left(y\fr{u}{v}, y , y \fr{w}{v}\right)\\
        r_{|W}: \left(z, \fr{u}{w}, \fr{v}{w}\right) & \mapsto & (x,y,z)= \left(z\fr{u}{w}, z \fr{v}{w}, z\right).
        \end{eqnarray*}
If $n=1$, $R_1$ is smooth and  the exceptional locus is given by one rational curve $C$. A direct computation shows that $C\cd C =-2$. 
If $n\geq 2$, $R_1$ has a singularity of type $A_{n-2}$ at the origin of $W$ and  the  exceptional locus 
is given by the union of two rational curves meeting at the singular point. 
Then, after a finite number of blow-up, we get a smooth surface.

Let $\ti{R}$ be the first smooth surface obtained in this way and $\rho : \ti{R} \ra R$ the composition of the blow-up morphisms.
Let $C=C_1,...,C_n$ be the components of the exceptional locus. 
A direct computation shows that $C_l \cd C_l =-2$, for any $l=1,...,n$.
From adjunction formula  we get $K_{\ti{R}}\cd C_l =0$, for any $l=1,...,n$. We want to prove that $\rho^* K_R \cong K_{\ti{R}}$.
But 
        \begin{align*}
        \rho^* K_R \cong K_{\ti{R}} + \sum_{l=1}^n a_l C_l,
        \end{align*}
for some integers $a_1,...,a_n$. The intersection of the right and left side of the previous espression gives:
        $$
        \sum_{l=1}^n a_l C_l \cd C_k =0 \quad \mb{for any} ~ k=1,...,n.
        $$
Since the matrix with entries $(C_l\cd C_k)$ is negative definite (\cite{BPV}, Chapter III, Theorem 2.1),
it follows $a_l=0 $ for all $l=1,...,n$.

From this description, it is clear that the exceptional locus $C$ is a chain of rational curves
whose dual graph is $\Ga_{A_n}$, (see Remark \ref{McKay}).
\end{ex}

\section{Existence and unicity}

\begin{prop}\label{existence and unicity of res}
Let  $Y$ be  a variety with transversal $ADE$-singularities. Then, there exists a unique crepant resolution $\rho :Z\ra Y$.
\end{prop}

\noindent \textbf{Proof.} To prove existence, one can proceed as follows. Let $r:Bl_S Y \ra Y$ be the blow-up 
of $S\subset Y$. If $Bl_S Y$ is smooth, then define $Z:=Bl_S Y$ and $\rho = r$. Otherwise, blow-up again.
As in the surface case, after a finite number of blow-up, we will end with a smooth variety.
Define $Z$ to be the first smooth variety obtained in this way, and $\rho$ be the composition of the blow-up morphisms.
We will show that $\rho^* K_Y \cong K_Z$. In general we have
        \begin{align*}
        \rho^* K_Y \cong K_Z + \sum_{l=1}^n a_l E_l,
        \end{align*}
where $E_l$ are the components of the exceptional divisor $E$ of $\rho$ and $a_l$ are integers defined as follows.
Let $z\in E_l$ be a generic point, and $ g_l=0$ be an equation for $E_l$ in a neighbourhood of $z$. 
Let $s$ be a (local) generator of $K_Y$ in a neighbourhood of $\rho (z)$. Then $a_l$ is defined by the following equation
        $$
        \rho^*(s) = g_l^{a_l} (dz_1 \wedge ... \wedge dz_n),
        $$
where $z_1,...,z_d$ are local coordinates for $Z$ in $z$. For more details see \cite{CKM}, Lecture 6.
In our case, $Y$ is locally a product $R\times \C^k$, so $Z$ is locally isomorphic to $\ti{R}\times \C^k$.
Then, since $\ti{R} \ra R$ is crepant, $a_l =0$ for all $l=1,...,d$.

We now prove unicity. Assume that $\rho_1 :Z_1 \ra Y$ is another resolution of $Y$. By  \cite{Fujiki}, Lemma 2.10,
the exceptional locus  of $\rho_1$ is of pure codimension $1$ in $Z_1$. Let $I_{S/Y}$ be the ideal sheaf of $S$ in $Y$.
The sheaf $J:=\rho_1^{-1}(I_{S/Y})\cd \mcal{O}_{Z_1}$ is the ideal sheaf of the exceptional locus of $\rho_1$,
and the local picture shows that it is invertible. So, we get a morphism
$Z_1 \ra Bl_S Y$ which lifts $\rho_1$,   \cite{H}, Chapter II, Proposition 7.14. 
Repeating this argument we get a morphism
$f: Z_1 \ra Z$, furthermore $f^* K_Z =K_{Z_1}$ becouse $\rho_1 :Z_1 \ra Y$ is crepant. Then $f$ is an isomorphism.
Indeed, $f$ induces a morphism $T_{Z_1} \ra f^* T_Z$, taking the determinant we have
a morphism 
        \begin{equation}\label{local iso}
        \wedge^d T_{Z_1} \ra \wedge^d T_Z,
        \end{equation}
so we get a global section of $\mcal{O}_{Z_1}(K_{Z_1} - f^* K_Z)\cong \mcal{O}_{Z_1}$.
Since $Z_1$ is projective and $f$ is birational, this section is constant equal to $1$. So it is a 
local isomorphism and one to one. \qed

\section{Geometry of the exceptional divisor}

In this section we restrict our attention to varieties with transversal $A_n$-singularities and whose associated 
orbifold $[Y]$ has trivial monodromy. In this case, any component $E_l$ of the exceptional divisor has a structure of 
$\Pro^1$-bundle on $S$. We will describe $E_l$ as the projectivization of vector bundles over $S$ of rank $2$. 
These vector bundles are defined in terms of the line bundles $L$, $M$ and $K$, defined in Lemma \ref{LM}. 

\vspace{0.2cm}

\begin{notation}\label{pil}\nf We will denote by $E\subset Z$ the exceptional locus of $\rho$ and
by $E_1,...,E_n$ the irreducible components of $E$.   The restriction of $\rho$ to $E$ will be denoted  $\pi : E\ra S$ 
and the restriction of $\pi$ to $E_l$ with $\pi_l :E_l \ra S$.
\end{notation}

\vspace{0.2cm}

\begin{prop}\label{E1}
Let $Y$ be a variety with transversal $A_1$-singularities. Then $E$ is irredicible and there exists a vector bundle $F$ on $S$, 
of rank two, such that
        \begin{equation}
        E\cong \Pro (F),
        \end{equation}
where $\Pro (F)$ is the projective bundle of lines in $F$ as defined in \cite{Fu}, Appendix B.5.5, where it is denoted by $P(F)$.

Moreover, the normal bundle $N_{E/Z}$ is given as follows
        \begin{equation}\label{N}
        N_{E/Z} \cong \mcal{O}_F(-2) \ot \pi^* L
        \end{equation}
where $L$ is defined by $\wedge^2 F \ot L \cong R^1 \pi_* N_{E/Z}$, $\mcal{O}_F(-2)$ is defined in \cite{Fu}, Appendix B.5.1.
\end{prop}

\noindent \textbf{Proof.} In this case, $Z=Bl_SY$, and the normal cone $C_S Y$ of $S$ in $Y$ is a conic bundle with fiber
isomorphic to $\{(x,y,z)\in \C^3: xy-z^2=0\}$, so the projection $C_S Y \ra S$ induces to $\pi: E=\Pro (C_S Y)\ra S$ 
a structure of a $\Pro^1$ bundle over $S$. In particular it is irreducible.

Since $S$ is smooth, there exists a rank two vector bundle on $S$, say $F$, such that $E\cong \Pro (F)$.
This follows from \cite{H}, Chapter II, Exercise 7.10(c). Let us fix one of these bindles and denote it by $F$.

The normal bundle $N_{E/Z}$ is a line bundle whose restriction on each fiber $\pi^{-1}(s)$ is isomorphic to
$\mcal{O}_{\pi^{-1}(s)}(-2)$. So, equation \ref{N} follows. The description of $L$ is a consequence of the projection
formula, \cite{H}, Chapter III, Exercise 8.3. \qed

\newpage

\begin{prop}\label{En}
Let  $Y$ be  a variety with transversal $A_n$-singularities, $n\geq 2$. Then, for any component $E_l$ of $E$,
$\pi_l :E_l \ra S$ is a $\Pro^1$-bundle and it can be written as follows,
        $$
        E_l \cong \Pro (L_l \op M_l ) \quad \mb{for}~ l=1,...,n,
        $$
where $L_l$ and $M_l$ are line bundles on $S$  that satisfies the following equation
        \begin{equation}\label{L_l-M_l}
        L_l \ot M_l^{\vee} \cong M \ot (K^{\vee})^{\ot l}  \quad \mb{for}~ l=1,...,n.
        \end{equation}
Moreover the intersection $E_k \cap E_l$ has the following expression
        \[  E_k \cap E_l =
        \begin{cases}
        \emptyset & \mb{if}~ |k-l| >1, \\
        \Pro (M_{l-1}) \subset E_{l-1} & \mb{if} ~ k=l-1,\\
        \Pro (L_l) \subset E_{l}& \mb{if} ~ k=l-1.
        \end{cases} \]
\end{prop}

\noindent \textbf{Proof.} In our case the monodromy is trivial, so we can identify the local groups $G_y$, $y\in S$,
with $\Z_{n+1}$ in a natural way. From this it follows that the normal cone $C_S Y$ of $S$ in $Y$ is the union of
two components, each of these being a vector bundle of rank $2$ on $S$. So, we get two components of $E$
with the  structures of $\Pro^1$-bundles over $S$. If we blow-up again, we get other components of $E$.
More precisely, if $n=3$ we get one component, if $n\geq 4$, then we must get two other componets. Each of these components 
have clearly a structure of $\Pro^1$-bundles over $S$.

Since $S$ is smooth, we can choose rank two vector bundles $F_l$ on $S$ such that $E_l \cong \Pro (F_l)$.
But, for later use, we want $F_l$ to be of the claimed form. Notice that, if $E_l \cong \Pro (F_l)$, then for any
line bundle $L$ on $S$, $E_l\cong \Pro (F_l \ot L)$. So, if $F_l = L_l \op M_l$, $E_l$ is dermined by $L_l \ot M_l^{\vee}$.

We prove the Proposition in the following way: we think of  $Z$ as a finite number of blow-ups; at each blow-up
we compute the transition functions of the vector bundles that forms the normal cone of the singular locus.

Let $\{ U_i \}_{i\in I}$ be an open covering of a neighbourhood of $S$ in $Y$. Assume that each 
$U_i$ is given as follows
        $$
        U_i =\{ (\ul{w}_i, x_i, y_i, z_i)\in \C^k\times \C^3 : x_i y_i - z_i^{n+1}=0 \}.
        $$
For any $i\in I$, let $(\ti{U_i}, \Z_{n+1}, {\chi_i}_{|} )$ be a uniformizing system for $[Y]$ such that $\chi_i(\ti{U_i})= \ti{U_i}$.
Let us assume that 
        $$
        \ti{U_i} = \{ (\ul{w}_i, u_i, v_i)\in \C^k\times \C^2 \}.
        $$
We can suppose further that, for any $i,j \in I$, there are isomorphisms
        $$
        \vphi_{ij} :  \chi_i^{-1}(U_i \cap U_j) \ra \chi_j^{-1}(U_i \cap U_j)
        $$
which are $\Z_{n+1}$-equivariant. 

Let $\vphi =(\Phi, F, G)$, where $\Phi, F, G$ are the components of $\vphi$ with respect to the coordinates $(\ul{w}_i, u_i, v_i)$.
The fact that  $\varphi$ is $\Z_{n+1}$-equivariant imply that $\Phi$ depends
only on $\ul{w}_i$ so it is an isomorphism  between open subsets of $S$.

Using this coordinates we get bases for the normal bundles 
$N_{\ti{U}_i^{\Z_{n+1}} / \ti{U}_i} = \lan \fr{\partial}{\partial u_i}, \fr{\partial}{\partial v_i}\ran $.
Moreover we have
        \begin{eqnarray*}
        \fr{\partial}{\partial u_i} &\mapsto & \fr{\partial F}{\partial u_i}\fr{\partial}{\partial u_j}+
        \fr{\partial G}{\partial u_i} \fr{\partial}{\partial v_j} \\
        \fr{\partial}{\partial v_i} & \mapsto & \fr{\partial F}{\partial v_i}\fr{\partial}{\partial u_j}+
        \fr{\partial G}{\partial v_i} \fr{\partial}{\partial v_i}.
        \end{eqnarray*}
The condition for $T\vphi_{ij}: N_{\ti{U}_i^{\Z_{n+1}} / \ti{U}_i} \ra N_{\ti{U}_j^{\Z_{n+1}} / \ti{U}_j}$ to be 
an isomorphism of $\Z_{n+1}$-representations implies that $F$ is independent 
on $v_i$ and that $G$ does not depend on $u_i$.

The two conditions $F(\e \cd u_i)=\e \cd F(u_i)$, $G(\e^{-1} \cd v_i)=\e^{-1} \cd G(v_i)$
imply
        \begin{eqnarray*}
        F(u_i) &=& \sum_{k=0}^{\infty} u_i^{k(n+1) +1} \fr{\partial^{k(n+1) +1} F}{\partial u_i^{k(n+1) +1}} \\
        G(v_i) &=&  \sum_{k=0}^{\infty} v_i^{k(n+1) +1}  \fr{\partial^{k(n+1) +1} G}{\partial v_i^{k(n+1) +1}}.
        \end{eqnarray*}
From these expressions we see that $F^{n+1}$ and $G^{n+1}$ depends only on $u_i^{n+1}=x_i$
and $v_i^{n+1}=y_i$ respectively and $F\cd G$ depends on $(x_i, y_i, z_i)$. So $(F^{n+1},G^{n+1}, F\cd G)$
is an automorphism of $U_i\cap U_j$. Moreover we have the following change of variable expression
        \begin{eqnarray*}
        x_j &=& x_i \left(\fr{\partial F}{\partial u_i} \right)^{n+1} +\mbox{h.o.t.'s}\\
        y_j &=& y_i  \left( \fr{\partial G}{\partial v_i} \right)^{n+1} +\mbox{h.o.t.'s}\\
        z_j &=& z_i \fr{\partial F}{\partial u_i}\fr{\partial G}{\partial v_i}+\mbox{h.o.t.'s}.
        \end{eqnarray*}

From the previous calculations it is clear that the normal cone of $S$ in $Y$
is the union of two irreducible components, $C_1$ and $C_2$. $C_1$ and $C_2$ have a structure of 
vector bundles of rank $2$ over $S$ and they are given as follows
        \begin{eqnarray*}
        C_1& \cong& M \oplus K\\
        C_2& \cong& L \oplus K,
        \end{eqnarray*}
where $L$, $M$ and $K$ are defined in Lemma \ref{LM}. Moreover, the intersection $C_1 \cap C_2$ in $C_S Y$ 
is given by the line bundle $K$.

After the first blow-up, we get a variety over $Y$ with transversal $A_{n-2}$-singularities, the exceptional
divisor is $\Pro( C_S Y) = \Pro(C_1) \cup \Pro(C_2)$, the singular locus is $\Pro(C_1)\cap \Pro(C_2)$.
So, if $n>2$, we have to blow-up again. Let $(a_i, b_i, z_i)$, $(a_j, b_j, z_j)$ be coordinates
in a neighbourhood of the singular locus.  The blow-up morphism, in these coordinates, is given by:
$x_i =a_i z_i$, $y_i =b_i z_i$, $z_i =z_i$. Then the two systems of coordinates, $(a_i, b_i, z_i)$ and  
$(a_j, b_j, z_j)$ are related ss follows:
        \begin{eqnarray}\label{ij}
        a_j &=& \fr{ x_j}{z_j} = \fr{F^{n+1}(a_i z_i)}{F \cd G}\\
        b_j &=& \fr{ y_j}{z_j} = \fr{G^{n+1}(a_i z_i)}{F\cd G}\\
        z_j &=& F\cd G.
        \end{eqnarray}
Note that, in the first two equations, numerator and denominator on the right side, are both multiples of $z_i$.
So, after dividing by $z_i$, \ref{ij} becomes
        \begin{eqnarray*}
        a_j &=& a_i \fr{ \left(\fr{\partial F}{\partial u_i} \right)^{n+1}}{\left( \fr{\partial F}{\partial u_i}
        \fr{\partial G}{\partial v_i}\right)}
        + \mbox{h.o.t.'s}\\
        b_j &=& b_i \fr{ \left(\fr{\partial G}{\partial v_i} \right)^{n+1}}{\left( \fr{\partial F}{\partial u_i}
        \fr{\partial G}{\partial v_i}\right)}
        + \mbox{h.o.t.'s}\\
        z_j &=& F\cd G.
        \end{eqnarray*}

These calculations shows that the normal cone of the singular locus, after the first blow-up,
is the union of the following irreducible components 
$\Pro ((M \ot K^{\vee}) \op K)$ and $\Pro ( K\op (L \ot K^{\vee}) )$ intersecting along $ \Pro (K)$.

Under the identification of  the strict transform of $\Pro(C_S Y)$ with $\Pro (M\op K) \cup \Pro(K\op L)$,
we claim that $\Pro ((M \ot K^{\vee}) \op K) \cap \Pro (M\op K) = \Pro(K) \subset \Pro (M\op K)$,
$\Pro ((M \ot K^{\vee}) \op K) \cap \Pro (M\op K) = \Pro(M \ot K^{\vee}) \subset \Pro (M \ot K^{\vee} \op K)$
and $\Pro (M\op K) \cap  \Pro ( K\op (L \ot K^{\vee} )) = \emptyset$.
This follows from the surface case. In this case $Y= \{ xy -z^{n+1}=0\}$
and the blow-up is covered by three open sets with coordinates
$(x,  \fr{c}{a})$, $( y, \fr{c}{b})$ and $(\fr{a}{c}, \fr{b}{c}, z)$.
The strict transform of the $x-$axis is contained in the open with coordinates
$(x,  \fr{c}{a})$, $ \fr{c}{a} = \fr{z}{x}$. Under the identification of the exceptional locus
with $\Pro( \{ (x,y,z) : xy =0 \})$. \qed

\section{Cohomology ring of the crepant resolution}

Let $Y$ be a variety with transversal $A_n$-singularities, suppose that the monodromy of $[Y]$ is trivial.
In this section we describe the cohomology ring of the crepant resolution $Z$ of $Y$ in terms of the cohomology of
$Y$ and the geometry of the exceptional divisor $E$.

\vspace{0.2cm}

\begin{notation}\nf We will use the same notation of the preceding section. Moreover, we will denote
by $j:E\ra Z$ the embedding of $E$ in $Z$, $i:S \ra Y$ the embedding of $S$ in $Y$. 
The restriction of $j$ to the component $E_l$ will be denoted by $j_l : E_l \ra Z$. 
We will denote by $k_l:E_l \ra E$ the morphism defined by the equality $j\circ k_l = j_l$.
\end{notation}

\vspace{0.2cm}

\begin{prop}\label{cohoa1}
Let $Y$ be a variety with transversal $A_1$-singularities. 
Then the following map is an  isomorphism of vector spaces 
        \begin{eqnarray*}
         H^*(Y) \oplus H^{*-2}(S)\lan E \ran &\cong &H^*(Z)\\
         \de +\al E&\mapsto & \rho^*(\de) + j_*\pi^*(\al).
        \end{eqnarray*}
Under the identification of $H^*(Z)$ with $H^*(Y) \oplus H^{*-2}(S)\lan E \ran$ by means of this map, the cup product of $Z$
is given as follows
        \begin{eqnarray*}
        & & (\de_1 +\al_1 E) \cdot (\de_2 +\al_2 E) =  \de_1 \cup \de_2 -2i_*(\al_1 \cup \al_2)\\
        & & + \left( i^*(\de_1)\cup \al_2 +\al_1 \cup i^*(\de_2)  + 2c_1(R^1 \pi_{\ast} N_{E/Z})\cup \al_1 \cup \al_2 \right)E.
        \end{eqnarray*}
\end{prop}

\noindent \textbf{Proof.} From projection formula we get
        $$
        j_*\pi^*(\al)\cdot \rho^*(\de)=j_*\left(
        \pi^*(\al)\cdot j^*\rho^*(\de)\right)=j_*\pi^*(\al\cdot i^*\de).
        $$
So, $\al E \cd \de = (\al \cup i^* \de )E$. To get the product rule between $\al_1$ and $\al_2$, for $\al_1, \al_2 \in H^*(S)$, 
we write
        $$
        j_*\pi^*(\al_1)\cup j_*\pi^*(\al_2)=\rho^*(\de) + j_*\pi^*(\al)
        $$
for some $\de \in H^*(Y)$ and $\al \in H^*(S)$. But, from projection formula and the equality $\rho_* \circ j_* =(\rho \circ j)_*=
(i\circ \pi)_*= i_*\circ \pi_*$, we have
        \begin{eqnarray*}
        \de = & \rho_*(j_*\pi^*(\al_1)\cup j_*\pi^*(\al_2)) &=-2i_*(\al_1 \cup \al_2), \quad \mb{and}\\
        & j_*\pi^*(\al_1)\cup j_*\pi^*(\al_2) & =j_* \left(c_1( N_{E/Z})\cup \pi^*(\al_1 \cup \al_2) \right).
        \end{eqnarray*} 
On the other hand, $\pi^*(\al)$ is the coefficient of $c_1(\mcal{O}_F(-2))$ in $j^*(j_*\pi^*(\al_1)\cup j_*\pi^*(\al_2))$, which is
$2 \al_1 \cup \al_2 \cup c_1 \left( R^1 \pi_{\ast} N_{E/Z} \right)$. \qed

\vspace{0.2cm}

\begin{notation} \nf For any variety $X$ and line bundle $L$ on $X$, we will denote by $L$
the first Chern class $c_1 (L) \in H^2(X)$. If $\al \in H^*(X)$, then we will denote by $\al L$ 
the cup product $\al \cup c_1(L) \in H^*(X)$.
\end{notation}

\begin{conv}\nf The variety $Y$ has complex dimension $d$ and so also $Z$ is of complex dimension $d$ and $S$ has complex dimension 
$k:=d-2$.
\end{conv}

\vspace{0.2cm}

\begin{prop}\label{cohomologyAn}
The following map is an isomorphism of vector spaces
        \begin{eqnarray*}
         H^*(Y) \oplus_{l=1}^n H^{*-2}(S)\lan E_l \ran & \ra & H^*(Z) \\
        \de + \al_1 E_1+ ...+ \al_n E_n &\mapsto & \rho^*(\de) + \sum_{l=1}^n {j_l}_* \pi_l^* (\al_l).
        \end{eqnarray*}
Under this identification the cup product of $Z$ is given as follows:
        \begin{eqnarray*}\label{cupn}
        E_{i-1} \cup E_i &=&  [S] +\sum_{l=1}^n \{ [(c_n^{-1})_{il} -(c_n^{-1})_{i-1,l}]M +
        [i(c_n^{-1})_{i-1,l}-(i-1)(c_n^{-1})_{il}]K \} E_l\\
        E_i \cup E_i &=& -2[S] + \sum_{l=1}^n \{ [(c_n^{-1})_{i-1,l} -(c_n^{-1})_{i+1,l}]M  \\
        & &  +  [-(i-1)(c_n^{-1})_{i-1,l}-4(c_n^{-1})_{il}+ (i+1)(c_n^{-1})_{i+1,l}]K \} E_l\\
        E_i \cup E_j &=&0 \qquad \mbox{if} \quad |i-j|>1
        \end{eqnarray*}
where $[S] \in H^4(X)$ denotes the fundamental class $i_*([S])$ of $S$ in $X$,  and 
$(c_n)_{ij}$ is the element in the $i-$th row and $j-$th column of the following $n\times n$ matrix
        \begin{equation}\label{cn} c_n = 
        \begin{pmatrix}
        -2  & 1  & 0  & ... & ... & 0 \\
         1  & -2 & 1  & 0   &  ...& 0 \\
        ... & ...& ...& ...& ... & ... \\
         0  & ...& 0 &  1 & -2  & 1  \\
        0   & ... & ... & 0 & 1 & -2 
        \end{pmatrix}
        \end{equation}
\end{prop}

\vspace{0.7cm}

\noindent The proof will use the following lemmas.

\begin{lem}
The following sequence is  exact for any $q$,
        $$
        0\ra H^q(Y) \xrightarrow{\rho^*} H^q(Z) \xrightarrow{[j^*]} H^q(E)/\pi^*(H^q(S)) \ra 0,
        $$
where $[j^*]$ is the composition of $j^*$ with the projection $ H^q(E)\ra H^q(E)/\pi^*(H^q(S))$.
The sequence splits, so we get an isomorphism of vector spaces
        $$
         H^*(Z)\cong H^*(Y)\oplus H^*(E)/\pi^*(H^*(S)).
        $$
\end{lem}

\noindent \textbf{Proof.} The exactness follows by comparing the exact sequences of the pairs $(E,Z)$ and 
$(S,Y)$. The sequence is split since there is a push-forward morphism
$\rho_* :H^*(Z)\ra H^*(Y)$ which satisfies $\rho_*\circ \rho^* =id_{H^*(Y)}$.\qed

\vspace{0.5cm}

\begin{lem}
There is a canonical isomorphism of vector spaces
        $$
        H^*(E)/\pi^*(H^*(S)) \cong \oplus_{l=1}^n H^*(E_l)/\pi_l^*(H^*(S)).
        $$
\end{lem}

\noindent \textbf{Proof.} Let $\hat{E}_1= \overline{E-E_1}$ be the closure of the complement of $E_1$ in $E$.
Then $E=E_1 \cup \hat{E}_1$ and $E_1 \cap \hat{E}_1\cong S$. From the Mayer-Vietoris sequence with respect to
the covering $\{E_1, \hat{E}_1\}$ of $E$ it follows that the following morphism is an isomorphism
        $$
        H^q(E)/\pi^*(H^q(S)) \xrightarrow{(k_1^*, -\hat{k}_1^*)} H^q(E_1)/\pi_1^*(H^q(S))\oplus H^q(\hat{E}_1)/\hat{\pi}_1^*(H^q(S))
        $$
where $k_1 :E_1 \ra E$ and $\hat{k}_1:\hat{E}_1 \ra E$ are  inclusions, and $\hat{\pi}_1$ is the restriction of $\pi$ to $\hat{\pi}$.
The result  follows by induction. \qed

\vspace{0.5cm}

\noindent \textbf{Proof of Proposition \ref{cohomologyAn}.} As a consequence of the above 
lemmas we have that the vector spaces $H^*(Y) \oplus_{l=1}^n H^{*-2}(S)$ and $H^*(Z)$ have the same dimension,
so it is enough to show 
that the map is injective. Assume that $\rho^*(\de) + \sum_{l=1}^n {j_l}_* \pi_l^* (\al_l)=0$. Then 
$\de =\rho_* (\rho^*(\de) + \sum_{l=1}^n {j_l}_* \pi_l^* (\al_l)) =0$. Next, applying $j_k^*$ we get
        $$
        j_k^* (\rho^*(\de) + \sum_{l=1}^n {j_l}_* \pi_l^* (\al_l)) =0,
        $$
the left side of this equation, up to terms of the form $\pi_k^*(..)$ assumes the following expression:
        $$
        \pi_k^*( -\frac{1}{2} \al_{k-1} +\al_k -\frac{1}{2} \al_{k+1} )\cup e_k,
        $$
where $e_k\in H^2(E_k)$ is the class $j_k^*\left( c_1(\mcal{O}_Z(E_k))\right)$. We use the convenction $\al_{-1}=\al_{n+1}=0$.
It follows that  
        $$
        (c_n)_{kl} \al_l =0\quad \mb{for any}~k=1,...,n.
        $$ 
So, $\al_l = 0$ for all $l$, since $c_n $ is nondegenerate. 

Finally note that the following maps are isomorphisms
        \begin{eqnarray*}
        H^{q-2}(S) &\ra& H^q(E_l)/\pi_l^*(H^q(S))\\
        \al &\mapsto& [\pi_l^*(\al) \cup e_l].  
        \end{eqnarray*}

To prove \ref{cupn}, write 
        \begin{equation}\label{general}
        E_i \cup E_j = \rho^*(\de) + \sum_{l=1}^n {j_l}_* \pi_l^* (\al_l),
        \end{equation}
where $\de \in H^*(Y)$ and $\al_l \in H^*(S)$. 

We immediately get
        \begin{eqnarray*}
        \de &=& \rho_*(E_i \cup E_j) = \rho_*({j_i}_*([E_i]) \cup {j_j}_*([E_j])) \\
        & = & \rho_* {j_i}_*([E_i]) \cup j_i^*{j_j}_*([E_j])) = i_* {\pi_i}_* ([E_i]) \cup j_i^*{j_j}_*([E_j])).
        \end{eqnarray*}
So, \[ \de = 
        \begin{cases}
        0    &  \mb{if}~ |i-j|>1 \\
        [S]  &  \mb{if}~ |i-j|=1 \\
        -2[S]& \mb{if}~ |i-j|=0 .
        \end{cases} \]

We now pull-back both sides of (\ref{general}) with $j_k$. We get the following equation
        \begin{eqnarray*}
        j_k^*(E_i \cup E_j)& = &j_k^*\rho^*(\de) +\sum_{l=1}^n j_k^*{j_l}_* \pi_l^* (\al_l) \\
        & = & \pi_k^* i^*(\de) + j_k^*( \al_{k-1}E_{k-1} + \al_{k}E_{k} + \al_{k+1}E_{k+1} \\
        & = & \pi_k^* i^*(\de) +\pi_k^*(\al_{k-1})[E_{k-1} \cap E_{k} \subset E_k] +\pi_k^*(\al_{k})N_{E_k /Z} \\
        & + & \pi_k^*(\al_{k+1})[E_{k+1} \cap E_{k} \subset E_k]
        \end{eqnarray*}
where by $[E_{k-1} \cap E_{k} \subset E_k]$ (resp. $[E_{k+1} \cap E_{k} \subset E_k]$) we mean the cohomology class dual 
of the homology class of $E_{k-1} \cap E_{k}$ (resp. $E_{k+1} \cap E_{k}$) in $E_k$. Using our description of $E_l$ 
(see Proposition \ref{En}) and \cite{Fu}, Appendix B.5.6, we have the following expression for $j_k^*(E_i \cup E_j)$
up to terms which are pulled-back from $S$ with $\pi_k$:
        \begin{equation}\label{right}
        j_k^*(E_i \cup E_j) = \pi_k^*(\al_{k-1} -2 \al_k +\al_{k+1})\mcal{O}_{F_k}(1),
        \end{equation}
where, as usual, $\al_{-1}=\al_{n+1}=0$.

On the other hand, the left side of (\ref{general}) gives
        \begin{eqnarray}\label{lhsgeneral}
        j_k^*(E_i \cup E_j) & = & j_k^*({j_i}_*([E_i]) \cup {j_j}_*([E_j])) \\
        & = & [E_i \cap E_k \subset E_k] \cap [E_j \cap E_k \subset E_k].
        \end{eqnarray}
We now distingiush three cases.

\noindent \textit{Case  $|i-j|>1$.} Then clearly $E_i \cup E_j=0$.

\vspace{0.5cm}

\noindent \textit{Case  $|i-j|=1$.} Then \[ j_k^*(E_{i-1} \cup E_i) =
        \begin{cases}
        0 & \mb{for} ~ k<i-1, \\
        N_{E_{i-1}/Z} \cup [E_{i-1} \cap E_i \subset E_{i-1}] & \mb{for} ~ k=i-1, \\
        N_{E_{i}/Z} \cup [E_{i-1} \cap E_i \subset E_{i}] & \mb{for} ~ k=i,\\
        0 & \mb{for} ~ k>i.
        \end{cases} \]
In order to compute $N_{E_{l}/Z}$, let us denote by $S_{l-1}= E_{l-1} \cap E_l$.
Then notice that ${N_{E_{l}/Z}}_{|S_{l-1}} \cong N_{S_{l-1}/E_{l-1}}$ and then, from \cite{Fu} Appendix B.5.6, it follows that
        \begin{equation}\label{Ni}
        N_{E_i/Z} \cong \mcal{O}_{F_i}(-2) + \pi_i^*(K-L_i -M_i).
        \end{equation}
A direct application of \cite{Fu} Appendix B.5.6 gives
        \begin{equation}\label{ElintEl-1inEl-1}
        [E_{l-1} \cap E_l \subset E_{l-1}] = \mcal{O}_{F_{l-1}}(1) + \pi_{l-1}^*L_{l-1}.
        \end{equation}
From (\ref{Ni}) and (\ref{ElintEl-1inEl-1}) it follows that, up to terms of the form $\pi_k^*(...)$, $j_k^*(E_{i-1} \cup E_i)$
is equal to \[ j_k^*(E_{i-1} \cup E_i) =
        \begin{cases}
        0 & \mb{for} ~ k<i-1, \\
        \mcal{O}_{F_{i-1}}(1)(iK-M) & \mb{for} ~ k=i-1, \\
        \mcal{O}_{F_{i}}(1)(M-(i-1)K) & \mb{for} ~ k=i,\\
        0 & \mb{for} ~ k>i.
        \end{cases} \]
Then we get the following system of equations for the $\al_l$'s.
\[ \begin{array}{ccccccccccc}
\left( \begin{array}{c} 
        0\\ ... \\ 0 \\ iK-M \\ M-(i-1)K \\ 0 \\ ... \\ 0   \end{array} \right)
\begin{array}{c}
        \\ \\ \\ = \\ \\ \\ \\ \end{array}
\left( \begin{array}{cccccccc} 
        -2 & 1 & 0 & ... & ... &... &... & 0  \\ 
        1  & -2& 1 & 0 & ... &... &... & 0  \\ 
        0  &  1&-2 & 1& 0 & ... &... & 0   \\
        0  & 0 &  1&-2 & 1& 0 & ... & 0   \\
        0 & 0  & 0 &  1&-2 & 1& ... & 0   \\
        ...&...&...&...&...  &...&...&... \\
        ...&...&...&...&...  &...&...&... \\
        0 & ...&... & ..&....&... & 1 &-2    \end{array} \right)
\left( \begin{array}{c} 
        \al_1 \\ ...  \\ ... \\ \al_{i-1} \\ \al_i \\ ... \\ ...\\ \al_n \end{array} \right) 
\end{array} \]

\vspace{0.5cm}

\noindent \textit{Case  $|i-j|=0$.} Then \[ j_k^*(E_{i} \cup E_i) =
        \begin{cases}
        0 & \mb{for} ~ k<i-1, \\
        [E_{i-1} \cap E_i \subset E_{i-1}]^2 & \mb{for} ~ k=i-1, \\
        N_{E_{i}/Z}^2 & \mb{for} ~ k=i,\\
        [E_{i+1} \cap E_i \subset E_{i+1}]^2 & \mb{for} ~ k=i+1, \\
        0 & \mb{for} ~ k>i+1.
        \end{cases} \]
Using again (\ref{Ni}) and (\ref{ElintEl-1inEl-1}), we have that \[ j_k^*(E_{i} \cup E_i) =
        \begin{cases}
        0 & \mb{for} ~ k<i-1, \\
        \mcal{O}_{F_{i-1}}(1)(M-(i-1)K) & \mb{for} ~ k=i-1, \\
        \mcal{O}_{F_{i}}(1)(-4K)& \mb{for} ~ k=i,\\
        \mcal{O}_{F_{i+1}}(1)((i+1)K-M) & \mb{for} ~ k=i+1, \\
        0 & \mb{for} ~ k>i+1.
        \end{cases} \]
Then we get the following system of equations for the $\al_l$'s.
\[ \begin{array}{cccccccccccc}
\left( \begin{array}{c} 
        0\\ ... \\ 0 \\ M-(i-1)K \\ -4K  \\ (i+1)K-M \\ 0 \\ ... \\ 0   \end{array} \right)
\begin{array}{c}
        \\ \\ \\ \\ = \\ \\ \\ \\ \end{array}
\left( \begin{array}{ccccccccc} 
        -2 & 1 & 0 & ... & ... &... &... &... & 0  \\ 
        1  & -2& 1 & 0 & ... &...&... &... & 0  \\ 
        0  &  1&-2 & 1& 0 & ... &...&... & 0   \\
        0  & 0 &  1&-2 & 1& 0 &...& ... & 0   \\
        0 & 0  & 0 &  1&-2 & 1&...& ... & 0   \\
        ...&...&...&...&...  &...&...&...&... \\
        ...&...&...&...&...  &...&...&...&... \\
        0 & ...&... & ..&....&...&... & 1 &-2    \end{array} \right)
\left( \begin{array}{c} 
        \al_1 \\ ...  \\ ... \\ \al_{i-1} \\ \al_i \\ \al_{i+1}\\ ... \\ ...\\ \al_n \end{array} \right) 
\end{array} \]

\qed

\chapter{Quantum corrections}

We compute the \textit{quantum corrected cup product} (see Definition \ref{quantum corrected cup product})
in the  transversal $A_n$ case and with trivial monodromy. 

In the first section we  give a description of the genus zero Gromov-Witten invariants of $Z$
that are needed in order to compute the quantum corrected cup product. 
Some of these invariants are computed by a direct understanding of the \textit{virtual fundamental class}.
To compute the ramaining, we use the property that Gromov-Witten invariants are invariants under deformation of the complex structure
of $Z$, so we will impose some technical hypothesis on $Z$ which guarantee that some deformations of $Z$ are unobstructed.     
In the last Section, using also the results of  Chapter 4, 
we give a presentation of the ring $H^*_{\rho}(Z)(q_1,...,q_n)$, as defined in Theorem \ref{qccohomology}.
We review in Section 2, some basic facts about virtula fundamental classes that are used in the proof
of the results.

\section{Gromov-Witten invariants of the crepant resolution}

We  describe some of the genus zero Gromov-Witten invariants of $Z$ in order to compute the quantum corrected cup product. 
Since we are able to compute  some of the invariants in complete generality and some others under particular hypothesis,
we decided to collect these results in three different Theorems.

\vspace{0.2cm}

\begin{conv}\nf Let $X$ be a variety of dimension $d$, then we define to be zero the integral of any cohomology
class $\al \in H^*(X)$   on $X$ of degree different to $2d$.
\end{conv}

\begin{notation}\nf Through this Chapter $Y$ will be a variety with transversal $A_n$-singularities
such that the corresponding orbifold $[Y]$ has trivial monodromy. We will denote by $\rho :Z \ra Y$
the crepant resolution.
\end{notation}

\vspace{0.5cm}

We have the following isomorphism of vector spaces (see Proposition \ref{cohomologyAn}):
        \begin{eqnarray*}
         H^*(Y) \oplus_{l=1}^n H^{*-2}(S)\lan E_l \ran & \ra & H^*(Z) \\
        \de + \al_1 E_1+ ...+ \al_n E_n &\mapsto & \rho^*(\de) + \sum_{l=1}^n {j_l}_* \pi_l^* (\al_l).
        \end{eqnarray*}
So, a cohomology class $\ga \in H^*(Z)$ of $Z$ will be denoted by 
        \begin{align*}
        \ga = \de + \al_1 E_1 +...+\al_n E_n, \quad \mb{with} ~ \de \in H^*(Y), ~ \al_l \in H^{*-2}(S).
        \end{align*}

The homology group $H_2(Z,\Z)$ of $Z$ can be described as follows,
        \begin{align*}
        H_2(Z,\Z) \cong H_2(Y\Z) \op H_0(S)\lan \be_1 \ran \op ... \op  H_0(S)\lan \be_n \ran,
        \end{align*}
where $\be_l \in H_2(Z,\Z)$ is the class of a fiber of $\pi_l :E_l \ra S$, see Notation \ref{pil}.
We have that $\be_1,...,\be_n$ is an integral basis of $\mb{Ker}~\rho_*$, where $\rho_* : H_2(Z,\Q) \ra H_2(Y,\Q)$
is the group homomorphism induced by $\rho$, see Chapter 2.4.
So, any class $\Ga \in H_2(Z,\Q)$ of a rational curve that is contracted by $\rho$, i.e. $\rho_*(\Ga)=0$,
is written in a unique way as
        $$
        \Ga = a_1 \be_1 + ... +a_n \be_n \quad \mb{with $a_l$ positive integers}.
        $$

\vspace{0.5cm}

We now compute the  $3$-point, genus zero  Gromov-Witten invariants 
        \begin{eqnarray} \label{GWI}
        \Psi_{\Ga}^Z(\ga_1, \ga_2, \ga_3 )=\int_{\left[\bar{\mcal{M}}_{0,3}(Z, \Ga)\right]^{vir}}ev_3^*(\ga_1 \otimes\ga_2 \otimes\ga_3)
        \end{eqnarray}
where $\ga_i \in H^*(Z)$, $\Ga=a_1\cdot \be_1 +...+a_n\cdot \be_n \in H_2(Z,\Z)$ 
is an homology class such that $\rho_*(\Ga)=0$, $\bar{\mcal{M}}_{0,3}(Z,\Ga)$ is the moduli space of $3$-pointed stable maps
$[\mu:(C,p_1,p_2,p_3)\ra Z]$ such that $\mu_*[C]=\Ga$, the
arithmetic genus of $C$ is $0$,  and
$ev_3 :\bar{\mcal{M}}_{0,3}(Z,\Ga) \ra Z\times Z \times Z$ is the evaluation map.

\vspace{0.5cm}

\begin{teo}\label{GWA1}
Let $Y$ be a variety with transversal $A_1$-singularities and let $\rho :Z\ra Y$ be the crepant resolution.
Then, \[ \Psi_{a\be}^Z(\ga_1, \ga_2, \ga_3 ) =
        \begin{cases}
        0 & \mb{if}~ \ga_1, \ga_2 ~\mb{or}~ \ga_3~\mb{are in}~ H^*(Y); \\
        -8\int_S \al_1 \cd \al_2 \cd \al_3 \cd c_1( R^1\pi_* N_{E/Z}) & \mb{if}~ \ga_i = \al_i E ~ \mb{for}~ i=1,2,3.
        \end{cases}\]
Where $a\geq 1$ is an integer, $\pi :E\ra S$ is the restriction to $E$ of $\rho$, and the dots denote the cup product of $S$.
\end{teo}

\vspace{0.5cm}

\begin{teo}\label{GWgrado1}
Let $Y$ be a variety with transversal $A_n$-singularities, $n\geq 2$, and let $\rho :Z\ra Y$ be the crepant resolution.
Then, \[ \Psi_{\Ga}^Z(\ga_1, \ga_2, \ga_3 ) =
        \begin{cases}
        0 \quad \mb{if}~ \ga_1, \ga_2 ~\mb{or}~ \ga_3~\mb{are in}~ H^*(Y);& \\
        (E_{l_1} \cd \be_{ij})(E_{l_2} \cd \be_{ij}) (E_{l_3} \cd \be_{ij}) \int_S\al_1\cd \al_2 \cd \al_3\cd c_1(K)&
        \end{cases}\]
where the second possibility holds if $ \Ga=\be_{ij}:= \be_i +...+\be_j$ for $1\leq i \leq j \leq n$,
        and $\ga_i=\al_i \cd E_{l_i}$ for $i=1,2,3$. Here $K$ is the line bundle on $S$ defined in Lemma \ref{LM}.
\end{teo}

\vspace{0.5cm}

\begin{teo} \label{GW}
Let $Y$ be a variety with transversal $A_n$-singularities, $n\geq 2$, and let $\rho :Z\ra Y$ be the crepant resolution.
Assume further that the line bundle $K$ (defined in Lemma \ref{LM}) is sufficiently ample,  that $H^2(Z, T_Z)=0$
and $H^1(S, T_S)=0$.
Then, we have the following expression for the Gromov-Witten invariants:
\[ \Psi_{\Ga}^Z(\ga_1, \ga_2, \ga_3 ) =
        \begin{cases}
        0\quad \mb{if}~ \ga_1, \ga_2 ~\mb{or}~ \ga_3~\mb{are in}~ H^*(Y);& \\
        (E_{l_1} \cd \be_{ij})(E_{l_2} \cd \be_{ij}) (E_{l_3} \cd \be_{ij}) \int_S\al_1\cd \al_2 \cd \al_3\cd c_1(K) &\\
        0\quad \mb{in the remaining cases}.
        \end{cases}\]
where the second possibility holds if $ \Ga=a \cd\be_{ij}:= \be_i +...+\be_j$ for $a$ be a 
positive integer, $1\leq i \leq j \leq n$, and $\ga_i=\al_i \cd E_{l_i}$ for $i=1,2,3$.
\end{teo}

\vspace{1cm}

\begin{conj}\label{mia congettura}\nf Theorem \ref{GW} is true in complete generality, i.e. without the hypothesis  on $K$,
on  $H^2(Z, T_Z)$ and $H^1(S,T_S)$. We give in Section 6.3  an outline of the proof of this conjecture.
\end{conj}

\vspace{0.5cm}

\begin{rem}\label{GW0} \nf Notice that, if $[Y]$ carries a global holomorphic symplectic $2$-form $\om$, then we can identify 
$L$ with $M^{\vee}$ by means of $\om$. So, 
        \begin{align*}
        (n+1)K\cong M\ot L \cong \mcal{O}_S,
        \end{align*}
so, all the Gromov-Witten invariants vanish.
\end{rem}

\section{Virtual fundamental class for Gromov-Witten invariants}

We review here some basic properties about the virtual fundamental class $\left[\bar{\mcal{M}}_{g,n}(Z, \Ga)\right]^{vir}$.
We follow clesely \cite{BF} and \cite{B}.

\vspace{0.5cm}

We recall the definition of  \textit{obstruction theory}\index{obstruction theory ! obstruction theory} 
for a Deligne-Mumford stack $\mcal{X}$.
We will denote by $D(\mcal{O}_{\mcal{X}_{\acute{e}t}})$ the derived category
of $\mcal{O}_{\mcal{X}_{\acute{e}t}}$-modules, where $\mcal{X}_{\acute{e}t}$ denote the small \'etale site of $\mcal{X}$. 
The cotangent complex of $\mcal{X}$ will be denoted by $L^{\bullet}_{\mcal{X}_{\acute{e}t}} \in D(\mcal{O}_{\mcal{X}_{\acute{e}t}})$.

\begin{defn} Let $\mcal{X}$ be a Deligne-Mumford stack.
Let  $E^{\bullet} \in ob D(\mcal{O}_{\mcal{X}_{\acute{e}t}})$ be an object that satisfies the following conditions:
\begin{itemize}
\item $h^i(E^{\bullet})=0$, for $i>1$; \\
\item $h^i(E^{\bullet})$ is coherent, for $i=0,-1$.
\end{itemize}

Then, a morphism $\phi :E^{\bullet} \ra L^{\bullet}_{\mcal{X}_{\acute{e}t}}$ in $D(\mcal{O}_{\mcal{X}_{\acute{e}t}})$
is called an \textbf{obstruction theory}\index{obstruction theory ! obstruction theory} for $\mcal{X}$, if $h^0(\phi)$ 
is an isomorphism and $h^{-1}(\phi)$
is surjective. By abuse of language we also say that $E^{\bullet}$ is an obstruction theory for $\mcal{X}$.

The obstruction theory  $\phi :E^{\bullet} \ra L^{\bullet}_{\mcal{X}_{\acute{e}t}}$ is \textbf{perfect}\index{obstruction theory ! perfect},
if $E^{\bullet}$ is of perfect amplitude in $[-1,0]$.
\end{defn}

\vspace{0.2cm}

Let $E^{\bullet}$ be a perfect obstruction theory for $\mcal{X}$. Assume that locally $E^{\bullet}$ 
is written as a complex of vector bundles $[ E^{-1} \ra E^0]$. Then, the rank of $E^{\bullet}$ is defined to be
        $$
        \rm{rk} ~ E^{\bullet}= \mb{dim}E^0 - \mb{dim} E^{-1}.
        $$

\begin{defn}
The \textbf{virtual dimension}\index{virtual dimension} of $\mcal{X}$ with respect to the perfect obstruction theory 
$E^{\bullet}$ is defined to be the rank $\rm{rk} ~ E^{\bullet}$ of $E^{\bullet}$.
We will denote the virtual dimension by $\nu$.
\end{defn}

\vspace{0.2cm}

\begin{rem}\nf The virtual dimension is a locally constant function on $\mcal{X}$. We shall assume that the virtual dimension
of $\mcal{X}$ with respect to $E^{\bullet}$ is constant, equal to $\nu$.
\end{rem}

\vspace{0.2cm}

\begin{rem}\nf Let $E^{\bullet}$ be a perfect obstruction theory for $\mcal{X}$. Then $E^{\bullet}$ give rise to a 
vector bundle stack $\mathfrak{E}$, into which the intrinsic normal cone $\mathfrak{C}_{\mcal{X}}$ of $\mcal{X}$  
can be embedded as a closed subcone stack (\cite{BF}, page 72). 
Under the hypothesis that $E^{\bullet}$ has a global resolution (\cite{BF}, Definition 5.2), 
the \textit{virtual fundamental class}\index{virtual fundamental class} $[\mcal{X},E^{\bullet} ]$ 
of $\mcal{X}$, with respect to $E^{\bullet}$, is the class in the rational Chow group $A_{\nu}(\mcal{X})$,  
obtained by intersecting $\mathfrak{C}_{\mcal{X}}$ with the zero section of $\mathfrak{E}$, \cite{BF}, Proposition 5.3. 
\end{rem}

\vspace{0.2cm}

\begin{notation} \nf We will denote by  $\bar{\mcal{M}}_{g,n}(Z, \Ga)$ the Deligne-Mumford stack of stable maps 
of class $\Ga \in H_2(Z)$ from an $n$-marked prestable curve of genus $g$ to $Z$. 
The \textit{universal curve}\index{universal curve} on $\bar{\mcal{M}}_{g,n}(Z, \Ga)$ will be denoted by 
$p: \mcal{C}\ra \bar{\mcal{M}}_{g,n}(Z, \Ga)$, and the \textit{universal stable map}\index{universal stable map} 
by $f :\mcal{C} \ra Z$. 

The universal curve $p: \mcal{C}\ra \bar{\mcal{M}}_{g,n}(Z, \Ga)$ can also be seen in the following way.
Consider the stack $\bar{\mcal{M}}_{g,n+1}(Z, \Ga)$, then, there is a morphism 
        \begin{equation}\label{forget}
        f_{n+1,n} :\bar{\mcal{M}}_{g,n+1}(Z, \Ga) \ra \bar{\mcal{M}}_{g,n}(Z, \Ga)
        \end{equation}
which forgets the last marked point and contracts the unstable components. Then $f_{n+1,n} $
can be identified with the universal curve, the universal stable map is the evaluation morphism
of the $(n+1)$-th point.

The \textit{cotangent complex}\index{cotangent complex} of $\bar{\mcal{M}}_{g,n}(Z, \Ga)$ 
will be denoted with $L_{\bar{\mcal{M}}_{g,n}(Z, \Ga)}^{\bullet}$. It
is an element in the derived category $D(\mcal{O}_{\bar{\mcal{M}}_{g,n}(Z, \Ga)_{\acute{e}t}})$ of the category of
$ \mcal{O}_{\bar{\mcal{M}}_{g,n}(Z, \Ga)_{\acute{e}t}}$-modules, where $\bar{\mcal{M}}_{g,n}(Z, \Ga)_{\acute{e}t}$ is the small \'etale
site of $\bar{\mcal{M}}_{g,n}(Z, \Ga)$.
\end{notation}

\vspace{0.2cm}

\begin{rem}\nf The moduli stack $\bar{\mcal{M}}_{g,n}(Z, \Ga)$ 
has a natural perfect obstruction theory\index{obstruction theory ! for Gromov-Witten invariants} given by
        \begin{equation}\label{obs}
        E^{\bullet}= R^{\bullet}p_* \{ [f^* \Om_Z \ra \Om_p ] \ot \om_p \},
        \end{equation}
where $\Om_Z$ is the sheaf of relative differentials of $Z$ over $\rm{Spec}( \C)$,
$\Om_p$ is the  sheaf of relative differentials of $\mcal{C}$ over $\mcal{M}$ (\cite{H}, page 175)
and $\om_p$ is the relative dualizing sheaf in degree $-1$, see \cite{BF}, page 82.
Moreover $E^{\bullet}$ has a global resolution. The resulting virtula fundamental class will be denoted by 
$[\bar{\mcal{M}}_{g,n}(Z, \Ga)]^{vir}$. 

This is proved in \cite{BF} Proposition 6.3, \cite{B} Proposition 5, for the relative case.
\end{rem}

\vspace{0.5cm}

There are some basic properties that holds for $[\bar{\mcal{M}}_{g,n}(Z, \Ga)]^{vir}$,
they are given and proved for the virtual fundamental class given by any perfect obstruction theory on 
a given Deligne-Mumford stack in \cite{BF}, Propositions 5.5-5.10 and Propositions 7.2-7.5.
Here we report some of these properties in the special case of $[\bar{\mcal{M}}_{g,n}(Z, \Ga)]^{vir}$.

\begin{prop}\label{vdim} The virtual dimension $\nu$ of $\bar{\mcal{M}}_{g,n}(Z, \Ga)$ is constant and equal to
        $$
        \nu = (1-g)(\mbox{dim} Z -3) -(\Ga \cd K_Z) +n.
        $$
\end{prop}

\vspace{0.2cm}

\begin{prop}\label{obs for smooth} Let  $\bar{\mcal{M}}_{g,n}(Z, \Ga)$ be smooth. Then $h^1({E^{\bullet}}^{\vee})$
is locally free and the virtual fundamental class is 
        $$
        [\bar{\mcal{M}}_{g,n}(Z, \Ga)]^{\text{vir}}=c_r (h^1({E^{\bullet}}^{\vee})) \cd [\bar{\mcal{M}}_{g,n}(Z, \Ga)],
        $$
where $r=\rm{rk} h^1({E^{\bullet}}^{\vee})$.
\end{prop}

\vspace{0.2cm}

\begin{prop}\label{pullback under forgetting}
Let 
        $$
        f_{n+1,n} :\bar{\mcal{M}}_{g,n+1}(Z, \Ga) \ra \bar{\mcal{M}}_{g,n}(Z, \Ga)
        $$
be the forgetful morphism. Then $f_{n+1,n}^*( [\bar{\mcal{M}}_{g,n}(Z, \Ga)]^{vir}) =[\bar{\mcal{M}}_{g,n+1}(Z, \Ga)]^{vir}$.
\end{prop}

\vspace{0.2cm}

Consider the following diagram of Deligne-Mumford stacks,
        \begin{equation} \label{pullback}
        \begin{CD}
        \mcal{X}' @>u>> \mcal{X}\\
        @VVV @VVV\\
        B' @>v>> B,
        \end{CD}
        \end{equation}
where $B$ and $B'$ are smooth of constant dimension, $v$ has finite unramified diagonal.

\begin{prop} 
Let $E^{\bullet} \ra L_{\mcal{X}/B}$ be a perfect obstruction theory for $\mcal{X}$ over $B$.
If (\ref{pullback}) is cartesian, then $u^* E^{\bullet} $
is a perfect obstruction theory for $\mcal{X}'$ over $B'$. If $E^{\bullet}$
has a global resolution so does $u^* E^{\bullet} $ and for the induced virtual fundamental classes
we have
        $$
        v^{!} [\mcal{M}]^{\text{vir}}=  [\mcal{M}']^{\text{vir}}
        $$
at least in the following cases.
\begin{description}
\item[1.] $v$ is flat,
\item[2.] $v$ is a regular local immersion.
\end{description}
\end{prop}

\vspace{0.2cm}

\begin{prop}\label{exact sequence for obstruction} Let $E^{\bullet}$ be the obstruction theory defined in (\ref{obs}).
Then we have the following exact sequence of coherent sheaves on $\mcal{M}$,
        \begin{equation}
        0 \ra p_* \Om_p^{\vee} \ra p_*f^* T_Z \ra h^0({E^{\bullet}}^{\vee}[1])
        \ra  R^1p_*\Om_p^{\vee} \ra R^1p_*f^* T_Z \ra h^1({E^{\bullet}}^{\vee}[1])\ra 0,
        \end{equation}
where 
$\Om_p^{\vee}=
R\mcal{H}om_{\mcal{O}_{\bar{\mcal{M}}_{g,n}(Z, \Ga)_{\acute{e}t}}}(\Om_p, \mcal{O}_{\bar{\mcal{M}}_{g,n}(Z, \Ga)_{\acute{e}t}})$.
\end{prop}
\noindent \textbf{Proof.} First of all notice that ${E^{\bullet}}^{\vee}\cong Rp_* \left([f^* \Om_Z \ra \Om_p]^{\vee}\right)$.
Indeed, by duality we have
        \begin{eqnarray*}
        {E^{\bullet}}^{\vee}&:=& 
        R\mcal{H}om_{\mcal{O}_{\bar{\mcal{M}}_{g,n}(Z, \Ga)_{\acute{e}t}}}\left( Rp_*([f^* \Om_Z \ra \Om_p]\ot \om_p), 
        \mcal{O}_{\bar{\mcal{M}}_{g,n}(Z, \Ga)_{\acute{e}t}}\right) \\
        &\cong& Rp_* R\mcal{H}om_{\mcal{O}_{\bar{\mcal{M}}_{g,n}(Z, \Ga)_{\acute{e}t}}}\left( [f^* \Om_Z \ra \Om_p]\ot \om_p,\om_p \right)
        \end{eqnarray*}
(see \cite{H1}, Sect. VII.4, page 393), then, since $\om_p$ is locally free, we have
        $$
        E^{\bullet} \cong
        Rp_* R\mcal{H}om_{\mcal{O}_{\bar{\mcal{M}}_{g,n}(Z, \Ga)_{\acute{e}t}}
        }\left( [f^* \Om_Z \ra \Om_p]\ot \om_p \ot  \om_p^{\vee}, \mcal{O}_{\bar{\mcal{M}}_{g,n}(Z, \Ga)_{\acute{e}t}}\right).
        $$

Consider the following distinguished triangle
        $$
        f^* \Om_Z \ra \Om_p \ra [f^* \Om_Z \ra \Om_p] \xrightarrow{+1},
        $$
where the sheaves $f^* \Om_Z$ and $\Om_p$ here means complexes centered in degree $0$,
the complex $[f^* \Om_Z \ra \Om_p]$ is the mapping cone of the morphism $df:f^* \Om_Z \ra \Om_p$.
Taking duals in derived category, we get the following distinguished triangle
        $$
        [f^* \Om_Z \ra \Om_p]^{\vee} \ra \Om_p^{\vee} \ra f^*T_Z  \xrightarrow{+1}.
        $$
The axioms of triangulated categories implies that  the following triangle is distinguished
        $$
        \Om_p^{\vee} \ra f^*T_Z \ra ([f^* \Om_Z \ra \Om_p]^{\vee})[1] \xrightarrow{+1}.
        $$
Now apply the derived functor $Rp_*$ to get the following distinguished triangle
        $$
        Rp_*(\Om_p^{\vee}) \ra Rp_*f^*T_Z \ra Rp_*([f^* \Om_Z \ra \Om_p]^{\vee})[1] \xrightarrow{+1}.
        $$
Taking cohomology, we get the following long exact sequence
        \begin{eqnarray*}
        &0& \ra R^{-1}p_*([f^* \Om_Z \ra \Om_p]^{\vee})[1] \ra p_*(\Om_p^{\vee}) \ra p_*f^*T_Z 
        \ra  p_*([f^* \Om_Z \ra \Om_p]^{\vee})[1] \\
        & \ra&  R^1p_*(\Om_p^{\vee}) \ra R^1p_*f^*T_Z 
        \ra R^1p_*([f^* \Om_Z \ra \Om_p]^{\vee})[1]\ra 0,
        \end{eqnarray*}
notice that $[f^* \Om_Z \ra \Om_p]$ is centered in $[-1,0]$ so $([f^* \Om_Z \ra \Om_p]^{\vee})[1]$
is also centered in $[-1,0]$.

Since $f$ is a stable map, the morphism $p_*\Om_p^{\vee} \ra p_*f^*T_Z$ is injective.
This conclude the proof. \qed

\newpage

\section{Proof of the Theorems}

\begin{notation}\nf Through this Section, $R$ will denote a rational double point of type $A_n$ and
$\ti{R} \ra R$ the crepant resolution. See Definition \ref{rdp definition} and Remark \ref{resolution graph}.
\end{notation}

\subsection{Genaral considerations}

The key point in the proof of the Theorems is to notice that $S$ is canonically isomorphic to 
the moduli space $\bar{\mcal{M}}_{0,0}(Z,\be_{ij})$, for any $i\leq j$, and that $\bar{\mcal{M}}_{0,0}(Z,\Ga)$
has a fibered structure over $S$ with fibers isomorphic to $\bar{\mcal{M}}_{0,0}(\ti{R},\Ga)$.

\vspace{0.2cm}

\begin{lem} \label{moduli}
There is a morphism
        $$
        \phi : \bar{\mcal{M}}_{0,0}(Z,\Ga) \ra S
        $$
such that,  for any point $p\in S$, the fiber $\phi^{-1}(p)$ is isomorphic to 
$\bar{\mcal{M}}_{0,0}(\ti{R},\Ga)$. Moreover, there is an \'etale cover $U\ra S$ and a cartesian diagram 
\begin{equation*}
        \begin{CD}
        U\times \bar{\mcal{M}}_{0,0}(\ti{R},\Ga) @>>> \bar{\mcal{M}}_{0,0}(Z,\Ga)\\
        @V pr_1 VV                                      @VV \phi V \\
        U @>>> S.
        \end{CD}
\end{equation*}
If $\Ga = \be_{ij}:= \be_i +...+\be_j$, for $i\leq j$, then $\phi$ is an isomorphism.
\end{lem}

\noindent \textbf{Proof.} \textit{Step 1.} We first prove that for any scheme $B$ of finite type over $\C$ and any object
\begin{equation*}
        \begin{CD}
        C @>f>> Z \\
        @V p VV \\
        B
        \end{CD}
\end{equation*}
in $\bar{\mcal{M}}_{0,0}(Z,\Ga)(B)$, there is a morphism $g:C\ra E$ such that $f=j \circ g$, where 
$j:E \ra Z$ is the inclusion.

Let $\mcal{O}_Z(E)$ be the line bundle over $Z$ associated to the divisor $E$, and let $s$ be the section of $\mcal{O}_Z(E)$
defined by $E$, that is, $s=\{ s_i \}$ where $s_i$ are functions which defines the Cartier divisor $E$.
Then $f$ factors through $E$ if and only if $f^*s$ vanishes as section of $f^*\mcal{O}_Z(E)$. 
We show that $p_* f^*\mcal{O}_Z(E)$ is the zero sheaf. 

First of all we assume that $B=\mbox{Spek}(\C)$. Then 
        $$
        \rho_* f_*([C]) =0,
        $$
where $\rho_*$ and $f_*$ are the morphisms of Chow groups induced by $\rho$ and $f$ respectively
and $[C]$ is the fundamental class of $C$. It follows that the image of $\rho\circ f $ is a point $y\in Y$.
This point must belong to $S$ because, outside $S$, $\rho$ is an isomorphism. So $f(C) \subset E$ and, since $C$
is reduced, $f$ factors through $E$ (\cite{H} exercise 3.11(d) chapter II). 
Note that $f(C)$ is contained in a fiber of $\pi :E \ra S$.

Assume now that  $B$ is a scheme of finite type over $\C$ and $f:C \ra Z$ is a stable map over $B$. 
Given a point $b\in B$, let $X$ be the subvariety of $B$ 
whose generic point is $b$, namely $X= \overline{ \{b\}}$. For any closed point  $x\in X$
we have that $H^0(C_x, f^*\mcal{O}_Z(E)_x)=0$ becouse $f\mid{C_x}$ factors through a fiber of $E$ over $S$. 
By Cohomology and Base Change it follows that $(p_* f^*\mcal{O}_Z(E))_x\ot k(x) =0$, since 
$p_* f^*\mcal{O}_Z(E)$ is coherent it follows that it is zero on a neighbourhood of $x$,
so it is zero on $b$.

\vspace{0.2cm}

\noindent \textit{Step 2.} Let $\varphi := \pi \circ g :C \ra S$,  we prove that there exists a morphism 
$\phi :B \ra S$ such that $\varphi = \phi \circ p$. 

First of all we define a continuous map $\phi :B \ra S$ such that $\varphi = \phi \circ p$.
From \textit{Step 1} it follows that, if $b\in B$ is a closed point, then we can define $\phi(b)$ by:
        \begin{align*}
        \phi(b)=\pi(f(C_b)).
        \end{align*}
Now, let  $b\in B$ be any point, and
let $X$ be the subvariety whose generic point is $b$. Then we define $\phi(b)$ to be the generic point of 
$\overline{\phi(X)}$, the closure of $\phi(X)$ in  $S$. The condition $\varphi = \phi \circ p$ 
implies that $\phi$ is continuous. Indeed, $p$ is surjective and so,
for any closed subset $T \subset S$, $\phi^{-1}(T)=\pi ( \pi^{-1}( \phi^{-1}(T))) =\pi (\varphi^{-1}(T))$.
Note that $p$ is surjective since it is dominant  and proper (the flatness of $p$ implies
that  $p^{\sharp}:\mcal{O}_B \ra p_* \mcal{O}_C$ is injective).

In order to give a morphism $\phi:B \ra S$ it remains to give a morphism of sheaves
        $$
        \phi^{\sharp}: \mcal{O}_S \ra \phi_* \mcal{O}_B.
        $$

We have the following diagram
        \begin{equation} \label{mor}
        \xymatrix{ \mcal{O}_S  \ar[r]^{\varphi^{\sharp}} & \varphi_*\mcal{O}_C \\
        & \varphi_* p^{-1} \mcal{O}_B \ar[u]_{\varphi_* p^{-1}(p^{\sharp})}}
        \end{equation}
where $p^{-1}(p^{\sharp}):p^{-1} \mcal{O}_B\ra \mcal{O}_C$ denote the adjoint of
$p^{\sharp}: \mcal{O}_B\ra   p_* \mcal{O}_C$. Note that  $\varphi_* p^{-1}(p^{\sharp})$ is injective.

Since $p$ is proper and surjective, a direct analysis show that the canonical morphism $\mcal{O}_B \ra p_* p^{-1} \mcal{O}_B$
is an isomorphism. So the morphism $\phi_* \mcal{O}_B \ra \varphi_* p^{-1} \mcal{O}_B= \phi_* (p_* p^{-1} \mcal{O}_B)$
is also an isomorphism. Then we can replace in the previous diagram (\ref{mor}), $\varphi_* p^{-1} \mcal{O}_B$
with $\phi_* \mcal{O}_B$, and get
        \begin{equation*}
        \xymatrix{ \mcal{O}_S  \ar[r]^{\varphi^{\sharp}} & \varphi_*\mcal{O}_C \\
        & \phi_* \mcal{O}_B \ar[u]}
        \end{equation*} 
with the vertical arrow being an inclusion.     
It follows that we can consider $\phi_* \mcal{O}_B$ as a subsheaf of 
$\varphi_*\mcal{O}_C$ and so it is enough to show that the image of $\mcal{O}_S$
under $\varphi^{\sharp}$ is contained in $\phi_* \mcal{O}_B$.
This is equivalent to say that the morphism $\mcal{O}_S \ra \varphi_* \mcal{O}_C /\phi_* \mcal{O}_B$
induced by $\varphi^{\sharp}$ is the zero morphism. But this is true on the geometric points and so it is true everywhere. 

The last statement follows from the fact that, if $\Ga = \be_1+...+\be_n$, then we have 
an inverse of $\phi$. It is given by sending any morphism $B\ra S$
to the following stable map 
\begin{equation*}
        \begin{CD}
        B\times_S E  @>j\circ pr_2>> Z \\
        @V pr_1 VV \\
        B
        \end{CD}
\end{equation*}\qed

\vspace{1cm}

\begin{lem}
Let $\ga_1,\ga_2$ or $\ga_3$ be elements of  $H^*(Y)$. Then 
        \begin{equation*}
        \Psi_{\Ga}^Z(\ga_1, \ga_2, \ga_3 )=0
        \end{equation*}
for any $\Ga = a_1 \be_1 +...+a_n \be_n$.
\end{lem}

\noindent \textbf{Proof.} By the Equivariance Axiom for Gromov-Witten invariants, see \cite{CK} Chapter 7.3,
we can assume that $\ga_3 = \rho^*(\de_3)$. The virtual dimension of $\bar{\mcal{M}}_{0,3}(Z,\Ga)$ is equal 
to the dimension of $Z$, $\rm{dim}Z$. So, let $\ga_1, \ga_2, \ga_3$ be cohomology classes such that
        \begin{align*}
        \rm{deg}(\ga_1)+\rm{deg}(\ga_2)+\rm{deg}(\de_3) = \rm{dim} Z. 
        \end{align*}
We have the following commutative diagram
        \begin{eqnarray*}
        \begin{CD}
        \bar{\mcal{M}}_{0,3}(Z,\Ga) @>ev_3>> Z\times  Z\times Z \\
        @V f_{3,2} \times f_{3,0}VV             @VV id \times id \times \rho V \\
         \bar{\mcal{M}}_{0,2}(Z,\Ga) \times \bar{\mcal{M}}_{0,0}(Z,\Ga) @>ev_2 \times \varphi >>Z\times  Z\times Y
        \end{CD}
        \end{eqnarray*}
where  $\varphi =i \circ \phi$. Then
        \begin{eqnarray*}
        ev_3^*(\ga_1 \otimes\ga_2 \otimes \rho^*(\de_3)) &=& (f_{3,2} \times f_{3,0})^*(ev_2 \times \varphi)^*
        (\ga_1 \otimes\ga_2 \otimes \de_3)\\
        &=&\left[ f_{3,2}^*ev_2^*(\ga_1 \otimes\ga_2)\right] \cdot
        \left[f_{3,0}^*\varphi^* (\de_3)\right]\\
        &=& f_{3,2}^*\left( \left[ev_2^*(\ga_1 \otimes\ga_2)\right] \cdot \left[f_{2,0}^*\varphi^* (\de_3)\right]\right)
        \end{eqnarray*}
we have used the equality $f_{3,0}=f_{2,0}\circ  f_{3,2}$. On the other hand, the following equalities hold
        \begin{eqnarray*}
        \left[\bar{\mcal{M}}_{0,3}(Z,\Ga)\right]^{vir}&=& f_{3,0}^*\left[\bar{\mcal{M}}_{0,0}(Z,\Ga)\right]^{vir} \\
        &=& f_{3,2}^*f_{2,0}^*\left[\bar{\mcal{M}}_{0,0}(Z,\Ga)\right]^{vir}\\
        &=& f_{3,2}^*\left[\bar{\mcal{M}}_{0,2}(Z,\Ga)\right]^{vir}.
        \end{eqnarray*}
So,
        \begin{eqnarray*}
        \Psi_{\Ga}^Z(\ga_1, \ga_2,\rho^*(\de_3) ) &=& \int_{f_{3,2}^*\left[\bar{\mcal{M}}_{0,2}(Z,\Ga)\right]^{vir}}
        f_{3,2}^*\left( \left[ev_2^*(\ga_1 \otimes\ga_2)\right] \cdot \left[f_{2,0}^*\varphi^* (\de_3)\right]\right)\\
        &=&(\mb{constant}) \cdot \int_{\left[\bar{\mcal{M}}_{0,2}(Z,\Ga)\right]^{vir}}
        \left[ev_2^*(\ga_1 \otimes\ga_2)\right] \cdot \left[f_{2,0}^*\varphi^* (\de_3)\right]
        \end{eqnarray*} 
which is zero since the virtual dimension of $\mcal{M}_{0,2}(Z,\Ga)$ is the virtual dimension of $\mcal{M}_{0,3}(Z,\Ga)$
minus $1$. \qed

\vspace{1cm}

It remains the computation of the invariants of the following form 
        $$
        \Psi_{\Ga}^Z({j_{l_1}}_* \pi_{l_1}^*(\al_1), {j_{l_2}}_* \pi_{l_2}^*(\al_2), {j_{l_2}}_* \pi_{l_2}^*(\al_2)),
        $$
where $\al_1,\al_2,\al_3 \in H^*(S)$ satisfy the following equation 
        \begin{align*}
        \rm{deg}~ \al_1 +\rm{deg}~ \al_2 +\rm{deg}~ \al_3 =\rm{dim} S -1.
        \end{align*}
With the next Lemma, we reduce this computation to an integral over the following class:
        \begin{align*}
        \phi_*[\mcal{M}_{0,0}(Z,\Ga)]^{vir} \in A_*(S)
        \end{align*}

\begin{lem}
The following equality holds,
        \begin{eqnarray*}
        & &  \Psi_{\Ga}^Z({j_{l_1}}_* \pi_{l_1}^*(\al_1), {j_{l_2}}_* \pi_{l_2}^*(\al_2), {j_{l_2}}_* \pi_{l_2}^*(\al_2))=  \\
        &=&(E_{l_1} \cd \Ga)(E_{l_2} \cd \Ga) (E_{l_3} \cd \Ga) 
        \int_{\phi_* [\bar{\mcal{M}}_{0,0}(Z,\Ga)]^{vir}} (\al_1 \cd \al_2 \cd \al_3).
        \end{eqnarray*}
\end{lem}

\noindent \textbf{Proof.} Consider the following cartesian diagram, which defines $E\times_S E\times_S E$,
        \begin{equation*}
        \begin{CD}
        E\times_S E\times_S E @>>> E\times E\times E \\
        @VVV                            @VVV\\
        S @>>> S\times S\times S
        \end{CD}
        \end{equation*}
where the botom arrow is the diagonal embedding.

From the previous Lemma \ref{moduli} it follows that the evaluation morphism $ev_3 : \bar{\mcal{M}}_{0,3}(Z,\Ga) \ra Z\times Z \times Z$
factors through a morphism $\ti{ev_3} : \bar{\mcal{M}}_{0,3}(Z,\Ga) \ra E\times E \times E$ and the inclusion
$E\times E \times E \ra Z\times Z \times Z$. Moreover, since the image of any stable map over a geometric point 
is contained in a fiber of $\pi :E\ra S$, $\ti{ev_3}$ factors through a morphism 
$\ti{\ti{ev_3}}: \bar{\mcal{M}}_{0,3}(Z,\Ga) \ra E\times_S E\times_S E$ and the inclusion 
$E\times_S E\times_S E \ra E\times E \times E$. So, we have the following equalities
        \begin{eqnarray*}
        & & \Psi_{\Ga}^Z({j_{l_1}}_* \pi_{l_1}^*(\al_1), {j_{l_2}}_* \pi_{l_2}^*(\al_2), {j_{l_2}}_* \pi_{l_2}^*(\al_2))=\\
        &=& \int_{[\bar{\mcal{M}}_{0,3}(Z,\Ga)]^{vir}} \ti{ev}_3^* \left( j^*{j_{l_1}}_* \pi_{l_1}^*(\al_1) \otimes 
        j^*{j_{l_2}}_* \pi_{l_2}^*(\al_2)
        \otimes j^*{j_{l_3}}_* \pi_{l_3}^*(\al_3)\right) \\
        &=& \int_{[\bar{\mcal{M}}_{0,3}(Z,\Ga)]^{vir}} 
        \ti{ev}_3^* \left(\mcal{O}_E(E_{l_1}) \pi^*(\al_1) \ot \mcal{O}_E(E_{l_2}) \pi^*(\al_2) \ot 
        \mcal{O}_E(E_{l_3}) \pi^*(\al_3)\right).
        \end{eqnarray*}
Notice that $j^*{j_l}_*\pi_l^*(\al)= \mcal{O}_E (E_l) \pi^*(\al)$, for any $l=1,...,n$ and $\al \in H^*(S)$.
Indeed, for any $m=1,...,n$, 
        $$
        k_m^* j^*{j_l}_*\pi_l^*(\al)=j_m^* {j_l}_*\pi_l^*(\al)=[E_m\cap E_l \subset E_m]\cd \pi_m^*(\al)
        $$ 
and 
        $$
        k_m^*(\mcal{O}_E (E_l) \pi^*(\al))=j_m^*\mcal{O}_Z(E_l) \pi_m^*(\al),
        $$ 
which are equal  (see \cite{Fu}, Chapter 2.3). Using the fact that $\ti{ev}_3$
factors through $\ti{\ti{ev_3}}$ we get
        \begin{eqnarray*}
        & & \int_{[\bar{\mcal{M}}_{0,3}(Z,\Ga)]^{vir}} \ti{ev}_3^* \left(\mcal{O}_E(E_1) \pi^*(\al_1) \ot 
        \mcal{O}_E(E_2) \pi^*(\al_2) \ot 
        \mcal{O}_E(E_3) \pi^*(\al_3)\right)\\
        &=& \int_{[\bar{\mcal{M}}_{0,3}(Z,\Ga)]^{vir}}\ti{\ti{ev_3}}\left( \mcal{O}_E(E_1)\ot \mcal{O}_E(E_2)
        \ot \mcal{O}_E(E_3)\pi^*(\al_1 \cd \al_2 \cd \al_3) \right).
        \end{eqnarray*}
Now apply the divisor axion and get
        \begin{eqnarray} \label{GWI3}
        & &  \Psi_{\Ga}^Z({j_{l_1}}_* \pi_{l_1}^*(\al_1), {j_{l_2}}_* \pi_{l_2}^*(\al_2), {j_{l_2}}_* \pi_{l_2}^*(\al_2))= \nonumber \\
        &=&(E_{l_1} \cd \Ga)(E_{l_2} \cd \Ga) (E_{l_3} \cd \Ga) \int_{\bar{[\mcal{M}}_{0,0}(Z,\Ga)]^{vir}} 
        \phi^* (\al_1 \cd \al_2 \cd \al_3).
        \end{eqnarray}

Let $pt$ denote the unique morphism from any scheme over $\C$ to $\rm{Spec}(\C)$.
Then the integral (\ref{GWI}) is the degree of the homology class 
        $$
        pt_* \left( ev_3^*(\ga_1 \otimes\ga_2 \otimes\ga_3) \cap \left[\bar{\mcal{M}}_{0,3}(Z, \Ga)\right]^{vir} \right).
        $$
Using the projection formula (\cite{Massey} pag. 328) and the equality $pt= \phi \circ pt$,  (\ref{GWI3}) becomes:
        $$
        (E_{l_1} \cd \Ga)(E_{l_2} \cd \Ga) (E_{l_3} \cd \Ga) pt_* \left( \al_1 \cd \al_2 \cd \al_3 \cap \phi_*
        [\bar{\mcal{M}}_{0,0}(Z,\Ga)]^{vir} \right).
        $$
\qed

\vspace{0.5cm}

At this point, to conclude the proof, we proceed as follows. First we show that 
        \begin{equation}\label{speranza}
         \phi_*[\bar{\mcal{M}}_{0,0}(Z,\Ga)]^{vir} = 
        \begin{cases}
        \fr{1}{a^3} [\bar{\mcal{M}}_{0,0}(Z,\be_{ij})]^{vir} & \text{if $\Ga = a \be_{ij}$},\\
        0 & \text{otherwise}.
        \end{cases}
        \end{equation}
Second, that 
        \begin{equation}
        [\bar{\mcal{M}}_{0,0}(Z,\be_{ij})]^{vir} = c_1(K). 
        \end{equation}

\vspace{0.5cm}

Let us  explain  our idea to prove (\ref{speranza}), it will be motivated by our proof of 
Theorems \ref{GWA1}, \ref{GWgrado1} and \ref{GW}. Let us denote by $\mcal{M}$ the moduli stack 
$\bar{\mcal{M}}_{0,0}(Z,\Ga)$ and by 
        \begin{align*}
        \phi_{\mcal{M}}: E_{\mcal{M}}^{\bullet} \ra L_{\mcal{M}}^{\bullet}
        \end{align*} 
the obstruction theory on $\mcal{M}$ given as in $(\ref{obs})$. We will identify $\bar{\mcal{M}}_{0,0}(Z,\be_{ij})$ with $S$
and  the obstruction theory will be denoted by 
        \begin{align*}
        \phi_S : E_S^{\bullet}\ra L_S^{\bullet}.
        \end{align*} 

\vspace{0.2cm}

\begin{rem}\nf Since $S$ is smooth, $L_S^{\bullet}$ is given by the locally free sheaf $\Om_S$
in degree zero.
\end{rem}

\vspace{0.2cm}

Suppose that we have a morphism $\phi^* \left( E_S^{\bullet} \right) \ra E_{\mcal{M}}^{\bullet} $ such that the 
following diagram commutes
        \begin{equation}\label{doppiasperanza}
        \begin{CD}
        \phi^* \left( E_S^{\bullet} \right) @>>> E_{\mcal{M}}^{\bullet} \\
        @V\phi^* \left( \phi_S \right)VV @VV \phi_{\mcal{M}}V \\
        \phi^* \left( L_S^{\bullet}\right) @>>> L_{\mcal{M}}^{\bullet}
        \end{CD}
        \end{equation}
where the botom row is the morphism induced by $\phi :\mcal{M} \ra S$. Then, we can complete (\ref{doppiasperanza})
to a morphism of distinguished triangles as follows
        \begin{equation}\label{ostruzione relativa}
        \begin{CD}
        & \phi^* \left( E_S^{\bullet} \right) & @>>> E_{\mcal{M}}^{\bullet} @>>> E_{\mcal{M}/S}^{\bullet} @>+1>> \\
        & @V\phi^* \left( \phi_S \right)VV & @VV \phi_{\mcal{M}}V @VV\phi_{\mcal{M}/S}V\\
        & \phi^* \left( L_S^{\bullet}\right) & @>>> L_{\mcal{M}}^{\bullet} @>>> L_{\mcal{M}/S}^{\bullet} @>+1>>
        \end{CD}
        \end{equation}
where $L_{\mcal{M}/S}^{\bullet}$ is the relative cotangent complex of $\mcal{M}$ over $S$ and 
$E_{\mcal{M}/S}^{\bullet} $ is, a priori, an object in $D(\mcal{O}_{\mcal{M}_{\acute{e}t}})$. 
We want to find conditions such that 
        \begin{align*}
        \phi_{\mcal{M}/S}: E_{\mcal{M}/S}^{\bullet} \ra L_{\mcal{M}/S}^{\bullet}
        \end{align*}
is a relative obstruction theory for $\mcal{M}$ over $S$, see \cite{BF} Section 7.
From diagram (\ref{ostruzione relativa}) it follows that, if the map
        $$
        h^{-1}\left(\phi^* \left( E_S^{\bullet} \right)\right) \ra h^{-1}\left(E_{\mcal{M}}^{\bullet}\right)
        $$
is injective, then $\phi_{\mcal{M}/S}: E_{\mcal{M}/S}^{\bullet} \ra L_{\mcal{M}/S}^{\bullet}$
is a relative obstruction theory. 

Suppose  that this is the case. Then, we have an induced perfect obstruction theory 
on each fiber of $\phi$,  \cite{BF} Proposition 7.2. An easy dimension count shows that 
this obstruction theory has virtual dimension zero. 

The fact that $\phi :\mcal{M} \ra S$ is locally trivial implies that, under the identification of 
each fiber of $\phi$ with $\bar{\mcal{M}}_{0,0}(\ti{R},\Ga)$, the obstruction theory induced on 
the former stack does not depends on the point in $S$. 

\begin{notation}\nf Let us denote by $F$ the stack $\bar{\mcal{M}}_{0,0}(\ti{R},\Ga)$, by $E_F^{\bullet}$
the obstruction theory on $F$ induced by $E_{\mcal{M}/S}^{\bullet} $ and by $[F]^{vir}$
the virtual fundamental class.
\end{notation}

\vspace{0.2cm}

Using the previous facts, we can prove the following formula
        \begin{equation}
        \phi_*[\mcal{M}]^{vir} = \mb{degree}([F]^{vir}) \cd [S]^{vir}
        \end{equation}
So, to prove (\ref{speranza}), it remains to show that
        \begin{equation}\label{speranzaF}
         \mb{degree}([F]^{vir}) = 
        \begin{cases}
        \fr{1}{a^3}  & \text{if $\Ga = a \be_{ij}$},\\
        0 & \text{otherwise}.
        \end{cases}
        \end{equation}

We have to understand $E_F^{\bullet}$. The obstruction theory $E_F^{\bullet}$
can be obtained as follows. Let $X\ra \De$ be a generic deformation of $\ti{R}$, where $\De \subset \C$
is a small open disc centered at the origin $0\in \C$. We mean that $\ti{R}$ is isomorphic to the fiber of $X\ra \De$ in $0$.
Notice that $X$ is a Calabi-Yau threefold. The embedding $\ti{R} \ra X$ gives a group homomorphism $H_2(\ti{R}, \Z) \ra H_2(X,\Z)$
and, by abuse of notation,  we will denote by $\Ga \in H_2(X,\Z)$ the image of  $\Ga \in H_2(\ti{R}, \Z)$.
Since $X\ra \De$ is generic, we will have an isomorphism as  follows,
        \begin{align*}
        F \cong \bar{\mcal{M}}_{0,0}(X,\Ga).
        \end{align*}
This follows easily from the explicit description of the semiuniversal deformation space of $\ti{R}$ given 
in \cite{KM} Theorem 1, see also \cite{T} Section 3. Then $E_F^{\bullet}$ coincides with the obstruction 
theory (\ref{obs}) for to $\bar{\mcal{M}}_{0,0}(X,\Ga)$.

Now, (\ref{speranzaF}) follows from the computation given in \cite{BKL} Proposition 2.10.

\subsection{The $A_1$ case}

Here we prove Theorem \ref{GWA1}. This is a straightforward generalization of \cite{LQ} Theorem 3.5 (ii).

\begin{notation}\nf  The exceptional divisor $E$ is irreducible and moreover the 
following description holds: $E\cong \Pro (F)$, where $F$ is a rank two vector bundle on $S$, 
$N_{E/Z} \cong \mcal{O}_F(-2) \ot \pi^* L$, where $L$ is defined by $\wedge^2 F \ot L \cong R^1 \pi_* N_{E/Z}$.
See  Proposition \ref{E1}.

In this section $\Ga$ will be $a\be$, for a positive integer $a$.

The universal stable morphism $f: \mcal{C} \ra Z$ factors through a morphism $g:\mcal{C} \ra E$. See Lemma \ref{moduli}.
\end{notation}

\vspace{0.2cm}

\begin{lem}
The moduli stack $\bar{\mcal{M}}_{0,0}(Z,a\be)$ is smooth of dimension $\mb{dim}~S+2a -2$. The virtual fundamental class is given by
        \begin{equation}\label{obs for smooth a1}
        [\bar{\mcal{M}}_{0,0}(Z,a\be)]^{vir}=c_r (h^1({E^{\bullet}}^{\vee})) \cd [\bar{\mcal{M}}_{0,0}(Z,a \be)]
        \end{equation}
where 
        \begin{equation}\label{obsa1}
        h^1({E^{\bullet}}^{\vee}) \cong R^1 p_* (g^*N_{E/Z})
        \end{equation}
is a vector bundle of rank $r = 2a - 1$.
\end{lem}

\noindent \textbf{Proof.} The smoothness of $\bar{\mcal{M}}_{0,0}(Z,a\be)$ 
follows from the fact that the  fibers of $\phi$ are smooth, where $\phi$ is 
defined in Lemma \ref{moduli}. Indeed they are all isomorphic to  $\bar{\mcal{M}}_{0,0}(\Pro^1,a\be)$.
Moreover from Proposition \ref{vdim} it follows that they have dimension $2a -1$. 

Equation (\ref{obs for smooth a1}) follows from Proposition \ref{obs for smooth}, so, it remains to prove equation (\ref{obsa1}).
We will show that 
        $$
         R^1 p_* (f^*T_Z) ~ \mb{is a vector bundle of rank} ~ 2a-1,
        $$
and that 
        \begin{equation}\label{non so}
        R^1 p_* (g^*N_{E/Z}) \cong  R^1 p_* (f^*T_Z).
        \end{equation}
Then the Lemma will follows from Proposition \ref{exact sequence for obstruction}.

\vspace{0.2cm}

Since $H^2 (p^{-1}(u),f^* T_Z |_{p^{-1}(u)})=0$, 
$R^2p_*(f^* T_Z)$ is locally free of rank $0$, see \cite{H} Theorem 12.11 (Cohomology and Base Change). 
So, it is enough to prove that  $H^1(p^{-1}(u),f^* T_Z )|_{p^{-1}(u)})$
is independent from $u$. 

Let $u=[\mu: D \ra Z] \in \bar{\mcal{M}}_{0,0}(Z,a \be)$ be a stable map. From the 
following exact sequence of locally free sheaves on $E$
        $$
        0\ra T_E \ra T_Z|_E \ra N_{E/Z}\cong \mcal{O}_E(-2) \ra 0,
        $$
we get the following
        $$
        0\ra \mu^*T_E \ra \mu^*T_Z|_E \ra \mu^* \mcal{O}_E(-2) \ra 0.
        $$
Since $H^1(D,\mu^* T_E)=0$, we get  
        $$
        H^1(D,\mu^* T_Z) \cong H^1(D,\mu^*\mcal{O}_{\mu(D)}(-2))
        $$
which has dimension $2a -1$. 

\vspace{0.2cm}

To prove (\ref{non so}), consider the following exact sequence
        $$
        0\ra T_E \ra {T_Z}_{|E} \ra N_{E/Z} \ra 0,
        $$
then apply $R^{\bullet}p_*$. \qed

\vspace{0.5cm}

\begin{rem}\nf From the previous Lemma, it follows that
        \begin{equation*}
        [\bar{\mcal{M}}_{0,0}(Z,\be)]^{vir} = c_1( R^1 \pi_*N_{E/Z}) \cd [\bar{\mcal{M}}_{0,0}(Z, \be)].
        \end{equation*}
\end{rem}

\vspace{0.5cm}

\begin{lem}
We have the following exact sequence:
        \begin{equation*}
        0\ra \phi^*(R^1 \pi_*N_{E/Z}) \ra R^1 p_*(g^*N_{E/Z}) \ra \mcal{F} \ra 0,
        \end{equation*}
where $\mcal{F}$ is a vector bundle of rank $2a-2$ whose restriction on each fiber of $\phi$ is given as follows.
Let $p\in S$ and consider the following commutative diagram:
        \begin{equation*}
        \begin{CD}
        \mcal{C}_{\phi^{-1}(p)} @>g_{|}>> E_{p} \\
        @Vp \mid VV @VV \pi_pV \\
        \phi^{-1}(p) @>\phi \mid >> \{p\}
        \end{CD}
        \end{equation*}
where $\mcal{C}_{\phi^{-1}(p)}$ is $p^{-1}(\phi^{-1}(p))$ and $p \mid $ (resp.$g \mid $) is the restriction of $p$ (resp. $g$)
on it, $E_p$ is the fiber  $\pi^{-1}(p)$ and $\pi_p$ is the restriction of $\pi$.
Then, the restriction of $\mcal{F}$ to $p^{-1}(\phi^{-1}(p))$ is:
        \begin{align*}
        {R^1p\mid }_*\left( {g\mid}^*(\mcal{O}_{E_p}(-1) \op \mcal{O}_{E_p}(-1)) \right).
        \end{align*}
\end{lem}

\noindent \textbf{Proof.} Since $E\cong \Pro(F)$, we have the surjective morphism: $\pi^*(F^{\vee}) \ra \mcal{O}_F(1)$.
Its kernel is $(\wedge^2 \pi^*(F^{\vee}))\ot \mcal{O}_F(-1)$. So, we have the following exact sequence
        \begin{align*}
        0\ra (\wedge^2 \pi^*(F^{\vee}))\ot \mcal{O}_F(-1) \ra \pi^*(F^{\vee}) \ra \mcal{O}_F(1)\ra 0,
        \end{align*}
which, tensorized with $(\pi^*\wedge^2 F\ot L)\ot \mcal{O}_F(-1)$ gives the following one
        \begin{equation}\label{come esce?}
        0\ra N_{E/Z} \ra \pi^*(F\ot L)\ot \mcal{O}_F(-1) \ra \pi^*(R^1\pi_* N_{E/Z}) \ra 0,
        \end{equation}
see Proposition \ref{E1}.

\vspace{0.2cm}

The pull back, under $g$, of (\ref{come esce?}) on $\mcal{C}$ gives a short exact sequence
of vector bundles, and then taking $R^{\bullet}p_*$ we have the following long exact sequence:
        \begin{eqnarray*}
        0& \ra & p_* g^*  N_{E/Z}  \ra p_* \left(p^* \phi^*(F\ot L) \ot g^* \mcal{O}_F(-1) \right) \ra p_*p^* \phi^*R^1\pi_* N_{E/Z} \\
         &\ra & R^1p_*( g^*  N_{E/Z}) \ra R^1p_* \left(p^* \phi^*(F\ot L) \ot g^* \mcal{O}_F(-1) \right) \\
        & \ra & R^1p_*(p^* \phi^*R^1\pi_* N_{E/Z}) \ra 0.
        \end{eqnarray*}
Notice that 
        \begin{align*}
        p_* \left(p^* \phi^*(F\ot L) \ot g^* \mcal{O}_F(-1) \right) \cong 0
        \end{align*}
by Cohomology and Base Change. Moreover, projection formula gives 
        \begin{eqnarray*}
         p_*p^* \phi^*R^1\pi_* N_{E/Z}& \cong &\phi^*(R^1\pi_* N_{E/Z}) \quad \mb{and}\\
         R^1p_*(p^* \phi^*R^1\pi_* N_{E/Z}) & \cong & \phi^*(R^1\pi_* N_{E/Z})\ot R^1p_* \mcal{O}_{\mcal{C}} \cong 0.
        \end{eqnarray*}

So, the result follows by defining 
        \begin{align*}
        \mcal{F} := R^1p_* \left(p^* \phi^*(F\ot L) \ot g^* \mcal{O}_F(-1) \right).
        \end{align*}
\qed 

\vspace{0.5cm}

Finally, to prove (\ref{speranza}), and so to complete the proof of Theorem \ref{GWA1}, we have to show that
        \begin{align*}
        \int_{\phi^{-1}(p)}c_{2a-2}\left(R^1{p_{|}}_*\left( g_{|}^*(\mcal{O}_{E_p}(-1) \op \mcal{O}_{E_p}(-1)) \right)\right) = \fr{1}{a^3}.
        \end{align*}
This is proved in \cite{CK} Theorem 9.2.3.

\subsection{The $A_n$ case, $n\geq 2$}

We prove Theorem \ref{GWgrado1}.

\begin{notation}\nf For any $i,j \in \{1,...,n\}$ with $1\leq i \leq j \leq n$, we will denote by
$E_{ij}$ the union $E_i \cup ... \cup E_j$. The restriction of $\pi$ to $E_{ij}$ will be denoted by $\pi_{ij}$.

The moduli stack $\bar{\mcal{M}}_{0,0}(Z, \be_{ij})$ will be denoted by $S$, see Lemma \ref{moduli}.
\end{notation}

\vspace{0.2cm}

Notice that we can assume that $i=1$ and $j=n$, so that $E_{ij}=E$ and $\be_{ij}=\be_1,...,\be_n$. 
The moduli stack is smooth of virtual dimension $\mb{dim}~S -1$, so the virtual fundamental class is given by
        \begin{align*}
        [S]^{vir} = c_1( h^1({E^{\bullet}_S}^{\vee})) \cd S
        \end{align*}
where $E^{\bullet}_S$ is given as follows, see (\ref{obs}):
        \begin{align*}
        E^{\bullet}_S = R^{\bullet}\pi_*([{\Om_Z}_{|E} \ra \Om_{\pi}] \ot \om_{\pi}).
        \end{align*}

\begin{lem}
        \begin{align*}
        h^1({E^{\bullet}_S}^{\vee}) \cong R^1 \pi_* N_{E/Z}.
        \end{align*}
\end{lem}

\noindent \textbf{proof.}
The  complex of sheaves on $E$:
        \begin{align}
        {\Om_Z}_{|E} \ra \Om_{\pi},
        \end{align}
is isomorphic, in the derived category $D(\mcal{O}_E)$, to a locally free sheaf $G$ in degree $-1$. 
Indeed, the morphism ${\Om_Z}_{|E} \ra \Om_{\pi}$  is surjective and, denoting by $G$  its kernel, 
we have the following exact sequence 
        \begin{align*}
        0\ra G \ra {\Om_Z}\mid{E} \ra \Om_{\pi} \ra 0.
        \end{align*}
Since $\Om_{\pi}$ is of projective dimension one, it follows that $G$ is locally free.

So, 
        \begin{align*}
        h^1({E^{\bullet}_S}^{\vee}) \cong R^1 \pi_* (G^{\vee}).
        \end{align*}

We have the following exact sequence
        \begin{align}\label{esattag}
        0\ra {\mcal{O}_Z(-E)}\mid{E} \ra G \ra \pi^* \Om_S \ra 0.
        \end{align}
This follows from a diagram chasing in the next  diagram, 
        \begin{align*}
        \xymatrix{\ & \ & \ & 0 \ar[d] & \ \\
                  \ & \ & \ & \pi^* \Om_S \ar[d] & \ \\         
                   0 \ar[r] &  \mcal{O}_Z (-E) \mid{E} \ar[r] & \Om_Z \mid{E} \ar[d]_= \ar[r] & \Om_E \ar[d] \ar[r] & 0 \\
                   0 \ar[r] &  G  \ar[r] & \Om_Z\mid{E}  \ar[r] & \Om_{\pi} \ar[d] \ar[r] & 0 \\
                  \ & \ & \ & 0  & \ }
        \end{align*}

Then, taking the dual of (\ref{esattag}) and applying $R^{\bullet} \pi_*$ we get the following isomorphism
        \begin{align*}
        R^1 \pi_* (G^{\vee}) \cong R^1 \pi_* N_{E/Z},
        \end{align*}
which complete the proof. \qed

\vspace{0.5cm}

The proof of Theorem \ref{GWgrado1} will be completed by the following Lemma.

\begin{lem}
        \begin{align*}
        R^1{\pi_{ij}}_* N_{E_{ij}/Z}\cong K.
        \end{align*}
\end{lem}

\noindent \textbf{Proof.} Let us first assume that $i<j$.  Let $a$ be an integer which satisfies
$i\leq a < j$. Let $\ti{E}_a =E_i\cup ... \cup E_a$ and $\bar{E}_{a+1} = E_{i+1} \cup ... \cup E_j$.
Let us denote by  $S_a =\ti{E}_a \cap \bar{E}_{a+1}$. Consider the following exact sequence
        $$
        0\ra \mcal{O}_{E_{ij}} \ra \mcal{O}_{\ti{E}_a } \op \mcal{O}_{\bar{E}_{a+1}} \ra \mcal{O}_{S_a} \ra 0,
        $$
this gives the following, 
        $$
        0\ra N_{E_{ij}/Z} \ra N_{\ti{E}_a/Z}(S_a) \op N_{\bar{E}_{a+1}/Z}(S_a) \ra {N_{E_{ij}/Z}}_{|S_a},
        $$
notice that $ {N_{E_{ij}/Z}} \mid{\ti{E}_a}=( \mcal{O}_Z(\ti{E}_a+\bar{E}_{a+1})) \mid{\ti{E}_a}=N_{\ti{E}_a/Z}\ot
\mcal{O}_Z(\bar{E}_{a+1}) \mid{\ti{E}_a}$ and ${\mcal{O}_Z(\bar{E}_{a+1})} \mid{\ti{E}_a}$ is, by definition,
the line bundle associated to the intersection of $\ti{E}_a$ with the divisor $\bar{E}_{a+1}$,
so it is $\mcal{O}_{\ti{E}_a}(S_a)$.

\vspace{0.2cm}

We now claim that 
        \begin{align*}
        R^1{\pi_{ij}}_* N_{E_{ij}/Z} \cong N_{E/Z} \mid S_a.
        \end{align*}
This follows from the long exact sequence associated to the functor $R^{\bullet}\pi_*$ once we notice that
        \begin{align}\label{7-10}
        R^p \pi_*N_{\ti{E}_a/Z}(S_a)=R^p \pi_*N_{\bar{E}_{a+1}/Z}(S_a)=0 \quad \mb{for all}~ p\geq 0.
        \end{align}
To see this, let $\ti{C}_a$ be a fiber of $\ti{E}_a$,  $\bar{C}_{a+1}$ be a fiber of $\bar{E}_{a+1}$
and  $p_a = \ti{C}_a \cap \bar{C}_{a+1}$. We can assume that $\ti{C}_a, \bar{C}_{a+1}$ and $ p_a$
are embedded in the surface $\ti{R}$ as a part of the exceptional divisor of $\ti{R} \ra R$.
Then, to prove (\ref{7-10}) is equivalent to prove the following identities,
        \begin{align*}
        H^p(\ti{C}_a,N_{\ti{C}_a/\ti{R}}(p_a))=H^p(\bar{C}_{a+1},N_{\bar{C}_{a+1}/\ti{R}}(p_a))=0 \quad \mb{for all}~ p\geq 0.
        \end{align*}
To prove the  equality: $H^p(\ti{C}_a,N_{\ti{C}_a/\ti{R}}(p_a)=0$, we proceed by induction on $a$. As before we have the exact sequence 
        \begin{align*}
        0\ra N_{\ti{C}_a/\ti{R}} \ra N_{\ti{C}_{a-1}/\ti{R}}(p_{a-1}) \op N_{C_a/\ti{R}}(p_{a-1}) 
        \ra {N_{\ti{C}_a/\ti{R}}}\mid{p_{a-1}}\ra 0
        \end{align*}
which tensorized by $\mcal{O}_{\ti{C}_a}(p_a)$ gives 
        \begin{align*}
        0\ra N_{\ti{C}_a/\ti{R}}(p_a) \ra N_{\ti{C}_{a-1}/\ti{R}}(p_{a-1}) \op N_{C_a/\ti{R}}(p_{a-1}+p_a)
         \ra {N_{\ti{C}_a/\ti{R}}}\mid{p_{a-1}}. 
        \end{align*}
To prove the second isomorphism notice that
        \begin{eqnarray*}
        {N_{E/Z}}\mid{S_a} &=& ( \mcal{O}_Z({E}_a+{E}_{a+1}))\mid{S_a} \\
        &= & (N_{E_a/Z}\ot {\mcal{O}_Z({E}_{a+1})}\mid{E_a})\mid{S_a} \\
        &=& {N_{E_a/Z}}\mid{S_a} \ot {\mcal{O}_{E_a}(S_a)}\mid{S_a}.
        \end{eqnarray*}

The result follows from the explicit description of the divisors $E_i$ in terms of the line bundles
$L_i, M_i$ and  $K$ given in Proposition \ref{En}. \qed

\vspace{0.7cm}

We now prove Theorem \ref{GW}. We will use the fact that Gromov-Witten invariants are invariants
under deformation of the complex structure of $Z$. The hypothesis of the Theorem will allow us 
to deform $Z$ in a convenient way.

\vspace{0.2cm}

\begin{notation}\nf For any variety $X$, we will denote by $T_X$ the sheaf of $\C$-derivations, i.e.,
        $$
        T_X= \mcal{H}om_{\mcal{O}_X}(\Om^1_X, \mcal{O}_X),
        $$
where $\Om^1_X$ is the  sheaf  of differentials of $X$. 
\end{notation}

\vspace{0.2cm}

\begin{prop}
Under the hypothesis $H^2(Y, T_Y)=0$, we have the following exact sequence
of cohomology groups:
        \begin{align}\label{9e15}
        0\ra H^1(Y, T_Y) \ra  H^1(Z, T_Z) \ra \op_{l=1}^n H^0(S,  R^1 {\pi_l}_* N_{E_l /Z})\ra 0.
        \end{align}
\end{prop}

\noindent \textbf{Proof.} From the Leray spectral sequence, we have the following exact sequence
of cohomology groups,
        \begin{eqnarray*}
        0&\ra& H^1(Y, T_Y) \ra  H^1(Z, T_Z) \ra H^0(Y,R^1\rho_* T_Z)\\
        &\ra&  H^2(Y, T_Y)\ra  H^2(Z, T_Z),
        \end{eqnarray*}
where we have used the fact that  $\rho_* T_Z$ is isomorphic to $ T_Y$. This is true since $Y$ has quotient singularities,
see \cite{S} Lemma 1.11 for the proof, see also \cite{BW} for a proof in dimension $2$. 
The vanishing of  $H^2(Y, T_Y)$ gives the following short exact sequence:
        \begin{align*}
        0\ra H^1(Y, T_Y) \ra  H^1(Z, T_Z) \ra H^0(Y,R^1\rho_* T_Z)\ra 0.
        \end{align*}

We now claim that $R^1\rho_* T_Z $ is isomorphic to $i_* \left( \op_{l=1}^n R^1 {\pi_l}_* N_{E_l /Z}\right)$.
First of all notice that there is a morphism  
        \begin{align}\label{9}
        R^1\rho_* T_Z \ra i_* \left( \op_{l=1}^n R^1 {\pi_l}_* N_{E_l /Z}\right).
        \end{align}
It is the composition of the morphism $i^*\left( R^1\rho_* T_Z \right) \ra R^1 \pi_* {T_Z}\mid{E}$, 
as defined in \cite{H} remark 9.3.1 Chapter III, and the morphism induced by $T_Z \mid{E} \ra \op_{l=1}^n N_{E_l /Z}$, as
the sum of $T_Z \mid{E} \ra N_{E_l /Z}$. We will prove that (\ref{9}) is an isomorphism.

Since this is a local problem, let us suppose that $Z= \C^k \times \ti{R}$. Then, from 
K\"{u}nneth formula and the fact that the surface singularity is rational, we have the following isomorphism
        \begin{align*}
        H^1(\C^k \times \ti{R}, T_{\C^k \times \ti{R}}) \cong H^0(\C^k, \mcal{O}_{\C^k}) \ot H^1(\ti{R}, T_{\ti{R}}).
        \end{align*}

Let  $C=C_1\cup ... \cup C_n \subset \ti{R}$ be the exceptional locus with the  $C_l$ be the
irreducible components. Then, we have the following isomorphism:
        \begin{align*}
        H^1(\ti{R}, T_{\ti{R}})\cong \op_{l=1}^n H^1(C_l, N_{C_l/\ti{R}}),
        \end{align*}
see \cite{BW} (1.8), and this shows that (\ref{9}) is an isomorphism. \qed

\vspace{0.5cm}

\begin{rem}\label{9e20} \nf It is known that $H^1(Y, T_Y)$ is in $1-1$ correspondence with the set of locally trivial
first order deformations of $Y$ modulo isomorphism. In the same way
$ H^1(Z, T_Z)$  is in $1-1$ correspondence with the set of first order deformations of $Z$ modulo isomorphisms.
So, the  sequence (\ref{9e15}) has the following meaning in deformation theory:
to any first order locally trivial deformation of $Y$ we can associate  a  first order deformation of
$Z$, the remaining deformations of $Z$ come from $H^0(Y,R^1\rho_* T_Z)$. We are interested in understanding
the last deformations.
\end{rem}

\vspace{0.5cm}

\begin{lem}\label{che casino}
For any subset $I\subset \{1,...,n\}$, let us denote by $E_I=\cup_{l\in I} E_l$.
Let $\pi_I :E_I \ra S$ be the restriction of $\pi$. Then $H^0(S, R^1{\pi_I}_* N_{E_I /Z})$ is an obstruction space for deformations 
of $E_I$ in $Z$. This means that, there is a morphism 
        \begin{align*}
        \rm{ob} :H^0(S, \op_{l=1}^n R^1 {\pi_l}_* N_{E_l /Z})\ra H^0(S, R^1{\pi_I}_* N_{E_I /Z})
        \end{align*}
which associates to any $\si \in H^0(S, \op_{l=1}^n R^1 {\pi_l}_* N_{E_l /Z})$ 
the obstruction \rm{ob}$(\si)$ to extend $E_I$ to the first order deformation associated to $\si$, see Remark \ref{9e20}.
\end{lem}

\noindent \textbf{Proof.} There is a morphism 
        \begin{align*}
        \rm{ob}:H^1(Z,T_Z) \ra H^1(E_I, N_{E_I/Z})
        \end{align*}
which associate, to any first order deformation $\mcal{Z}$ of $Z$, the obstruction to the existence of  a 
first order deformation $\mcal{E}_I$ of $E_I$ in $\mcal{Z}$, see \cite{Sernesi} Proposition II.3.3. 
Since $H^1(Y,T_Y)$ corresponds to locally trivial
deformations  $\mcal{Y}$ of $Y$, any element $\mcal{Y}_t$ of the family has tranversal $A_n$
singularities. The deformation $\mcal{Z}$ of $Z$ induced by $\mcal{Y}$ is a simultaneous resolution,
so it contains a deformation of $E_I$. It follows that ob factors through $H^0(S, \op_{l=1}^n R^1 {\pi_l}_* N_{E_l /Z})$.
By abuse of notations we will denote by ob this induced morphism.

To conclude  the proof, we will show that 
        \begin{align*}
        H^1(E_I, N_{E_I/Z}) \cong H^0(S, R^1{\pi_I}_* N_{E_I /Z}).
        \end{align*}
We use the Leray spectral sequence, $H^p(S, R^q{\pi_I}_* N_{E_I /Z}) \Rightarrow H^{p+q}(E_I,N_{E_I /Z})$,
to get the following exact sequence:
        \begin{eqnarray*}
        0&\ra& H^1(S,{\pi_I}_* N_{E_I /Z} ) \ra  H^1(E_I,N_{E_I /Z}) \ra H^0(S, R^1{\pi_I}_* N_{E_I /Z})\\
        &\ra&  H^2(S,{\pi_I}_* N_{E_I /Z} ),
        \end{eqnarray*}
where $\pi_I :E_I \ra S$ is  the restriction of $\pi$ to $E_I$.  
Since ${\pi_I}_* N_{E_I /Z}=0$ and $H^2(S,{\pi_I}_* N_{E_I /Z} )=0$, we have
an isomorphism between $H^1(E_I,N_{E_I /Z})$ and $H^0(S, R^1{\pi_I}_* N_{E_I /Z})$.\qed

\vspace{1cm}

Let $\si \in H^0(S,\op_{l=1}^n R^1 {\pi_l}_* N_{E_l /Z})$
be a section, and let 
        \begin{align*}
        \begin{CD}
        Z @>>> \mcal{Z}_1 \\
        @VVV @VVV \\
        \rm{Spec}(\C) @>>>\rm{Spec}\left( \fr{\C [\e]}{(\e^2)}\right)
        \end{CD}
        \end{align*}
be the corresponding first order deformation of $Z$. Then, there exists a finite deformation of $Z$
which, at the first order, coincides with the one given by $\si$, we will denote this deformation as follows
        \begin{align}\label{deffinita}
        \begin{CD}
        Z @>>> \mcal{Z} \\
        @VVV @VVV \\
        \{ 0\} @>>> \De
        \end{CD}
        \end{align}
where $\De$ is a small disc in $\C$ around the origin $0\in \C$. Indeed, by the Kodaira-Nirenberg-Spencer  Theorem (1958), 
there exists a complete family of deformations of $Z$, see \cite{Sernesi} for a review.

\vspace{0.2cm}

Let $V$ be a  neighbourhood  of $E$ in $Z$ and  $r:V\ra S$ be a morphism whose restriction to $E$
coincides with $\pi$.  Now, the deformation (\ref{deffinita}) induces a deformation $\mcal{V}$
of $V$. We can choose $\mcal{V}$ such that $r$ extends to $\ti{r} : \mcal{V} \ra S$, here we use the hypothesis $H^1(S,T_S)=0$.
Moreover, for $K$  sufficiently ample,
the deformation (\ref{deffinita}) can be chosen in such a way the locus of points $p\in S$ such 
that there exists a rational curve in $\ti{r}^{-1}(p)$ whose homology class 
is $\Ga$ has codimension one, if $\Ga = a \be_{ij}$, and has codimension strictly greater then one if
$\Ga \not= a \be_{ij}$. 

\vspace{0.2cm}

The previous considerations shows that, if $\Ga \not= a \be_{ij}$, then all the Gromov-Witten invariants vanish.
Indeed the virtual dimension of $\bar{\mcal{M}}_{0,0}(Z, \Ga)$ is $\rm{dim}~S -1$. 
If $\Ga = a \be_{ij}$, then the locus of points  $p\in S$ such 
that there exists a rational curve in $\ti{r}^{-1}(p)$ whose homology class is $\Ga$ is $\{ p\in S :~ \rm{ob}(\si)(p)=0 \}$,
where $\si \in H^0(S,\op_{l=1}^n R^1 {\pi_l}_* N_{E_l /Z})$ and $\rm{ob}$ is the morphism defined 
in Lemma \ref{che casino}.

\section{Quantum corrected cohomology ring for $A_n$-singularities}

Let $Y$ be a variety with transversal $A_n$-singularities and $\rho :Z \ra Y$ be the crepant resolution. 
We  compute the quantum corrected cup product $*_{\rho}$ Defined in Chapter 2.4.

\begin{conv}\nf We will tacitly assume that uor spaces $S,Y,Z$ fulfil the required hypothesis in order to be able
to apply Theorem \ref{GWA1}, \ref{GWgrado1} or \ref{GW}. 
\end{conv}

\begin{notation}\nf Since the crepant resolution $\rho:Z\ra Y$ is unique, the quantum corrected cup product $*_{\rho}$ will 
be denoted by  $*$.
\end{notation}

\vspace{0.2cm}

Let $\ga_1,\ga_2 \in H^*(Z)$, then the product $\ga_1 * \ga_2$ will be represented by a family of cohomology 
classes in $H^*(Z)$ which depend on $n$ complex parameters $q_1,...,q_n$. More precisely, for any $\ga_1,\ga_2,\ga_3 \in H^*(Z)$,
one defines the quantum corrected triple intersection $\lan \ga_1, \ga_2, \ga_3 \ran_{\rho}(q_1,...,q_n)$ as follows:
        \begin{align*}
        \lan \ga_1, \ga_2, \ga_3 \ran_{\rho}(q_1,...,q_n)=\lan \ga_1, \ga_2, \ga_3 \ran +\lan \ga_1, \ga_2, \ga_3 \ran_{qc}(q_1,...,q_n),
        \end{align*}
where $\lan \ga_1, \ga_2, \ga_3 \ran =\int_Z \ga_1 \cup \ga_2 \cup \ga_3$ and $\lan \ga_1, \ga_2, \ga_3 \ran_{qc}(q_1,...,q_n)$
is the quantum corrected $3$-point function, see (\ref{qc3}). Note that, $\lan \ga_1, \ga_2, \ga_3 \ran_{\rho}(q_1,...,q_n)$,
is a complex valued function defined on the domain of definition of the quantum corrected $3$-point function. 
Then, the quantum corrected cup product, $\ga_1 * \ga_2$, is defined by the following equation:
        \begin{equation*}
        \lan \ga_1 \ast \ga_2, \ga \ran =\lan \ga_1 ,\ga_2, \ga \ran_{\rho}(q_1,...,q_n), \quad \mb{for all}~ \ga \in  H^*(Z).
        \end{equation*}
So, $\ga_1 \ast \ga_2$ is a cohomology class in $H^*(Z)$ which depends on the parameters $q_1,...,q_n$.
The vector space $H^*(Z)$ with the product $*$ forms a family of rings depending on the parameters $q_1,...,q_n$,
see Theorem \ref{qccohomology}. This family will be denoted by $H^*(Z)(q_1,...,q_n)$.
Clearly the parameters $q_1,...,q_n$ belongs to the domain where the quantum corrected $3$-point function is defined.

\vspace{0.2cm}

\begin{notation}\nf Let us define $\ga_1 *_{\e} \ga_2$ by the following equation
        \begin{equation}\label{edef}
        \lan \ga_1 *_{\e} \ga_2, \ga \ran = \lan \ga_1, \ga_2, \ga \ran_{qc}(q_1,...,q_n) , \quad \mb{for all}~ \ga \in  H^*(Z).
        \end{equation}
Then, if $\ga_1 \cup \ga_2 = \de +\al_1 E_1 +...+\al_n E_n$ and $\ga_1 *_{\e} \ga_2 =  \de_{\e} +\e_1 E_1 +...+\e_n E_n$,
        \begin{align*}
         \ga_1 * \ga_2 = \de +\de_{\e} +(\al_1+\e_1) E_1 +...+(\al_n+\e_n) E_n.
        \end{align*}
\end{notation}

\vspace{0.2cm}

\begin{rem}\nf If $\ga_1 \in H^*(Y)$ or $\ga_2 \in H^*(Y)$, then $\ga_1 \ast \ga_2=0$. Indeed, under these hypothesis,
all the Gromov-Witten invariants $\Psi_{\Ga}^Z(\ga_1, \ga_2, \ga_3 )$ are zero.

So, we can assume $\ga_1 =E_i$ and $\ga_2=E_j $. 
\end{rem}

\vspace{0.2cm}

\begin{rem}\nf $E_i *_{\e} E_j \in H^*(Y)^{\perp}$, where $H^*(Y)^{\perp}$ is the subspace of $H^*(Z)$ which is orthogonal 
to $H^*(Y)$ with respect to the ususal pairing $\lan,\ran=\int_Z$.
\end{rem}

\begin{lem}\label{stare}
The following expression holds:
        \begin{equation}
        E_i *_{\e} E_j = \sum_{l,m=1}^n (c_n^{-1})_{lm} R_{ijm}(\ul{q})c_1(K)E_l,
        \end{equation}
where, $c_n$ is the $n\times n$ matrix (\ref{cn}), $K$ is the line bundle on $S$ defined in Lemma \ref{LM}
if $n\geq 2$ and is $R^1\pi_*(N_{E/Z})$ if $n=1$, $(\ul{q})$
denote $(q_1,...,q_n)$, and $R_{ijm}(\ul{q})$ is defined by the  following expression:
        \begin{eqnarray*}
        R_{ijm}(\ul{q}) = \sum_{r \leq s}(E_i \cd \be_{rs})(E_j \cd \be_{rs})(E_m \cd \be_{rs})
        \fr{q_{r} \cd \cd \cd q_s}{1-q_{r} \cd \cd \cd q_s}.
        \end{eqnarray*}
As usual, $\be_{rs}= \be_r+...+\be_s$.
\end{lem}

\noindent \textbf{Proof.} Let $E_i *_{\e} E_j = \e_1 E_1 +...+\e_n E_n$. Then, the left hand side of (\ref{edef})
becomes
        \begin{eqnarray*}
        \lan E_i *_{\e} E_j,\al E_k \ran &= & \int_Z \sum_{l=1}^n {j_l}_* \pi_l^*(\e_l) \cup \al E_k \\
        & = & \sum_{l=1}^n \int_Y \rho_* ( {j_l}_* \pi_l^*(\e_l) \cup \al E_k) \\
        & =& \sum_{l=1}^n \int_Y \rho_*{j_l}_* (\pi_l^*(\e_l) \cup {j_l}^*( \al E_k))\\
        &=& \sum_{l=1}^n \int_Y i_*{\pi_l}_* (\pi_l^*(\e_l\cup  \al) \cup [E_l \cap E_k \subset E_l]) \\
        & = & \int_S (\e_{k-1} -2 \e_k + \e_{k+1})\cup \al.
        \end{eqnarray*}
On the other hand
        \begin{eqnarray*}
        \lan E_i, E_j, \al E_k\ran_{qc}(\ul{q}) &=& \sum_{a=1}^{\infty} \sum_{1\leq r \leq s \leq n}
        (q_r\cd \cd \cd q_s)^a (E_i \cd \be_{rs})(E_j \cd \be_{rs})(E_k \cd \be_{rs})\int_S \al \cup c_1(K)\\
        &=& R_{ijk}(\ul{q})\int_S \al \cup c_1(K).
        \end{eqnarray*}

So, we get the following equations for $\e_1,...,\e_n$:
        \begin{equation*}
        \int_S (\e_{k-1} -2 \e_k + \e_{k+1})\cup \al =R_{ijk}(\ul{q})\int_S \al \cup c_1(K) \quad \mb{for}~ \al \in H^*(S).
        \end{equation*}
\qed

\vspace{0.2cm}

As a consequence, we have the following result.
\begin{prop}
The following expression holds for $E_i * E_j$:
        \begin{equation}
        E_i * E_j = \sum_{l,m=1}^n (c_n^{-1})_{lm}\{ R_{ijm}(\ul{q})c_1(K) + \al_{ijm}\}E_l,
        \end{equation}
where \[ (\al_{ij1},...,\al_{ijn}) =
        \begin{cases}
        0 & \\ \mb{if}~ |i-j|>1 ;&\\
        (0,...,0, iK-M, M-(i-1)K,0,..., 0)  & \\  \mb{if}~ j=i-1;&\\
        (0,...,0,M-(i-1)K,-4K, (i+1)K-M ,0, ..., 0) & \\  \mb{if}~ j=i, &
        \end{cases}\]
where, in the second row, $iK-M$ is in the $(i-1)$-th place, and in the third $-4K$ is in the $i$-th place.
\end{prop}

\chapter{Comparisons for $A_1$, $A_2$ and conclusions}

We put together the computations of the previous chapters in order to verify Ruan's conjecture for orbifolds with transversal 
$A_n$-singularities. Actually we have a complete picture only when $n=1,2$. These cases gives informations about how things should go
in the general case.

\section{The $A_1$ case}

Form Lemma \ref{stare} we immediately get
        $$
        E*_{\e} E = 4\fr{q}{1-q}c_1(R^1 \pi_{\ast} N_{E/Z})E.
        $$
So, from Proposition \ref{cohoa1}, we get the following expression for the quantum corrected cup product:
        $$
        E*E = -2i_*([S]) + \left( 2+ 4\fr{q}{1-q}\right) c_1(R^1 \pi_{\ast} N_{E/Z})E.
        $$

The orbifold cohomology ring $H^*_{orb}([Y])$ has been computed in Example \ref{orba1}.
It is easy to see that the following morphism is a ring isomorphism
        \begin{eqnarray*}
        H^*_{orb}([Y]) &\ra & H^*(Z)(-1)\\
        (\de ,\al) &\mapsto& (\de , \frac{\sqrt{-1}}{2} \al),
        \end{eqnarray*}
where $H^*(Z)(-1)$ is the ring $H^*(Z)(q)$, defined in Theorem \ref{qccohomology}, for $q=-1$.

\section{The $A_2$ case}

We will tacitly assume that uor spaces $S,Y,Z$ fulfil the required hypothesis in order to be able
to apply Theorem \ref{GWA1}, \ref{GWgrado1} or \ref{GW}. 

\begin{notation}\nf We will use the following notation: $\de_1= \fr{q_1}{1-q_1}$, $\de_2= \fr{q_2}{1-q_2}$
and $\de_3= \fr{q_1q_2}{1-q_1q_2}$.
\end{notation}

\vspace{0.5cm}

The following expressions holds for $*_{\e}$: 
        \begin{eqnarray*}
        E_1 \ast_{\e} E_1 &=& -2 [S] + \left( (2+4\de_1 +\de_3)M +(3+4\de_1 +\de_3)L  \right)E_1\\
                & & + \left( (\de_1 +\de_3)M +(2+\de_2 +\de_3)L  \right)E_2\\
        \\
        E_1 \ast_{\e} E_2 &=&  [S] + \left( (-1-2\de_1 +\de_3)M +(-2\de_1 +\de_3)L  \right)E_1\\
                & & + \left( (-2\de_2 +\de_3)M +(-1-2\de_2 +\de_3)L  \right)E_2\\
        \\
        E_2 \ast_{\e} E_2 &=& -2 [S] + \left( (2+\de_1 +\de_3)M +(\de_1 +\de_3)L  \right)E_1\\
                & & + \left( (3+4\de_2 +\de_3)M +(2+\de_2 +\de_3)L  \right)E_2.
        \end{eqnarray*}

\vspace{0.5cm}

\begin{rem}\nf Notice that, for $q_1=q_2=-1$, $\de_3 = \infty$, so, we have to modify slightly 
the cohomological crepant resolution conjecture (see Conjecture \ref{conj}). In any case, 
\end{rem}

\vspace{0.5cm}

The orbifold cohomology ring for transversal $A_2$-singularities has been computed in Example \ref{orba2}.
The resulting ring is described as follows:
        \begin{eqnarray*}
        e_1 \cup_{orb} e_1 &=& Le_2 \\
        e_1 \cup_{orb} e_2 &=& \fr{1}{3}[S]\\
        e_2 \cup_{orb} e_2 &=& Me_1.
        \end{eqnarray*}

\vspace{0.2cm}

\begin{rem}\label{s} \nf Notice that the previous expressions, for the quantum corrected cup product $*$ and for the 
orbifold cup product, are symmetric if we exchange $E_1 \leftrightarrow E_2$, $L\leftrightarrow M$
and $e_1 \leftrightarrow e_2$.
\end{rem} 

\newpage

\begin{teo}\label{a2}
The pairs $(q_1,q_2)$ for which there exists a ring isomorphism $H^*(Z)(q_1,q_2)\ra H^*_{orb}([Y])$ which
respects the symmetry described in Remark \ref{s} are: $q_1 =q_2 =exp(\fr{2}{3}\pi i)$ and $q_1 =q_2 =exp(\fr{4}{3}\pi i)$.
\end{teo}

\noindent \textbf{Proof.} We are looking for a linear isomorphism 
        \begin{equation}\label{isoruan} \begin{array}{cccc}
        \left(\begin{array}{c}
         E_1 \\ E_2
        \end{array} \right) =
        \left(\begin{array}{cc}
          a & b \\ c & d 
        \end{array}\right)
        \left( \begin{array}{c}
        e_1 \\ e_2 \end{array}\right).
        \end{array} \end{equation}
In order for it to respects the symmetry, we must have $a=d$ and $b=c$. The condition to be a ring homomorphism
gives a lot of equations, in particular we have $\de_1 = \de_2$. Moreover we get the following two possibilities:
$a= \sqrt{3}\rm{exp}(\fr{\pi}{6}i)$ and $b= \sqrt{3}\rm{exp}(\fr{5}{6}\pi i)$, or 
$a= - \sqrt{3}\rm{exp}(\fr{5}{6}\pi i)$ and $b= \sqrt{3}\rm{exp}(\fr{\pi}{6}i)$. The first choice corresponds to 
$q_1 =q_2 =\rm{exp}(\fr{2}{3}\pi i)$, the second to $q_1 =q_2 =\rm{exp}(\fr{4}{3}\pi i)$. \qed

\vspace{0.5cm}

\begin{com}\nf  As pointed out in Remark \ref{GW0}, if the orbifold $[Y]$ carries a holomorphic symplectic $2$-form,
then our Gromov-Witten invariants are zero. It is easy to see that the isomorphisms $H^*(Z)(q_1,q_2)\ra H^*_{orb}([Y])$
founded in Theorem \ref{a2} are still isomorphisms in this case, i.e. in the holomorphic symplectic case.
So, to find an isomorphism in the $A_n$ case, one can try to find an isomorphic under the additional hypothesis
for $[Y]$ to carries an holomorphic symplectic $2$-form, and then try to prove that this is still an isomorphism
in the general case. Of course the natural candidates for the parameters $q_1,...,q_n$ is 
$q_1=...=q_n$ equal to a $n$th root of unit, such that $q_1\cd\cd\cd q_n \not=1$.
\end{com}

\printindex

\end{document}